\documentclass{msc}
\usepackage{amsmath}
\usepackage{amsthm}
\usepackage{xy}
\usepackage{stmaryrd}         
\usepackage{amsfonts}

\undef{\gtrsim}
\undef{\lesssim}

\makeatletter
\DeclareRobustCommand*\cal{\@fontswitch\relax\mathcal}
\makeatother

\usepackage{amssymb}
\usepackage{hyperref}
\usepackage{graphicx}

\input diagxy

\DeclareMathOperator{\id}{id}
\DeclareMathOperator{\Id}{Id}
\def\Prof{{\cal P}\!{\it rof}}
\def\Tense{{\cal T}\!\!{\it ense}}

\def\Sett{{\bf Set}^{{\bf A}^{op} \times {\bf B}}}

\def\Pb{\framebox{\scriptsize{\rm Pb}}}

\def\longtod{\xy\morphism(0,0)|m|<250,0>[`;\hspace{-1mm}\bullet\hspace{-1mm}]\endxy}

\newcommand{\ov}[1]{{\bar{#1}}}

\newcommand{\todd}[2]{\xymatrix@1{\ar[r]|\bb^{#1}_{#2}&}}

\newcommand{\bgr}{\mbox{\rule[-2.5mm]{0mm}{7mm}}}  

\theoremstyle{plain}
\newtheorem{xtheorem}{Theorem}[section]
\newtheorem{xproposition}{Proposition}[section]
\newtheorem{xlemma}{Lemma}[section]
\newtheorem{xcorollary}{Corollary}[section]

\theoremstyle{definition}
\newtheorem{xdefinition}{Definition}[section]

\theoremstyle{remark}

\newtheorem{xexample}{Example}[section]

\input xy
\xyoption{all}

\def\liml{\varprojlim}
\def\limr{\varinjlim}

\def\Cat{{\cal C}\!{\it at}}

\def\DelF{\Delta{[}F{]}}

\newbox\bbox
\setbox\bbox = \hbox to 0pt{\hss $\scriptstyle\bullet$\hss}
\def\bb{\usebox{\bbox}}

\newcommand{\dedication}[1]{\par\vspace{-3ex}\hfill\emph{#1}\par}

\begin{document}


\lefttitle{Multivariate functorial difference}
\righttitle{Robert Par\'e}

\papertitle{Article}

\jnlPage{1}{00}
\jnlDoiYr{2019}
\doival{10.1017/xxxxx}

\title{Multivariate functorial difference}

\begin{authgrp}
\author{Robert Par\'e}
\affiliation{Department of Mathematics and Statistics, Dalhousie University\\
Halifax, NS, Canada, B3H 4R2\\\
        \email{pare@mathstat.dal.ca}}
\end{authgrp}

\history{(Received xx xxx xxx; revised xx xxx xxx; accepted xx xxx xxx)}

\begin{abstract}

Partial difference operators for a large class of functors
between presheaf categories are introduced, extending
our previous work on the difference operator to the
multivariable case. These combine into the Jacobian
profunctor which provides the setting for a lax chain
rule. We introduce a functorial version of multivariable
Newton series whose aim is to recover a functor from
its iterated differences. Not all functors are recovered but
we get a best approximation in the form of a left adjoint,
and the induced comonad is idempotent. Its fixed points
are what we call soft analytic functors, a generalization
of the well studied multivariable analytic functors.

\end{abstract}


\maketitle

\dedication{In memory of Phil Scott, 1947--2023}

\section*{Philip Scott}

I knew Phil for most of his career, from when he was a
post-doctoral fellow at McGill in 1977, a colleague the following
year at Dalhousie, and a friend ever since. His knowledge
of the literature in category theory, logic and computer
science was phenomenal. He travelled a lot and spoke
to many people. This way, he kept up to date on the latest
developments and each time he visited Halifax, he had
some new topic he thought I should look at. This was
good advice which I wish I had been more diligent
following up. We've lost a great ambassador for our
subject as well as a friend. I dedicate this work to him.

\section*{Introduction}

This is a sequel to \citep{Par24}. Here we are interested in
the structure of functors \( {\bf Set}^{\bf A} \to {\bf Set}^{\bf B} \)
(\(\bf A\) and \(\bf B\)  small categories)
generalizing the difference calculus for endofunctors of
\( {\bf Set} \). An important example is given by the generalized
analytic functors of \citet{FioGamHylWin08}. As in that work,
profunctors are central. That is perhaps the main difference
the present work has with \citep{Par24}. This is somewhat of
a simplification like saying that multivariate calculus is just
single variable calculus plus linear algebra. The added
dimensions open up a whole array of possibilities.

The work here is a categorified version of the classical
partial difference operators for real functions
\[
{\mathbb R}^n \to {\mathbb R}^m \rlap{\,,}
\]
a discrete version of partial derivatives. The analogy is
quite fruitful.

As the paper is quite long, it may be helpful to point out the
main results, namely the lax chain rule (Theorem~\ref{Thm-ChRule})
and the Newton adjunction (Theorem~\ref{Thm-NewtonAdj})
together with the convergence theorem (Theorem~\ref{Thm-UnitIso}).
These results are proper to the categorical setting and
have no counterpart for real-valued functions. They could not
be formulated without the pivotal definitions of the (discrete)
Jacobian as a profunctor (Definition~\ref{Def-JacobianProf})
and soft analytic functor (Definition~\ref{Def-SoftAnalytic}).

Apart from the obvious \citep{FioGamHylWin08} and the
references therein, the present work was strongly influenced
by the work of the Calgary-Ottawa-Montreal consortium on
tangent categories and cartesian differential categories.
Several talks in the ATCAT seminar by regulars Geoff
Cruttwell and Marcello Lanfranchi as well as guest speakers,
notably Robin Cockett and JS Lemay, helped form my
ideas on the categorical theory of differentials. After
completion of this work, the paper ``Cartesian difference categories''
by \citet{AlvLem21} came to my attention.
This is clearly relevant as it deals with the categorical
understanding of finite difference. What is less clear
is precisely how they are related. Further work in this
direction should prove fruitful.

Thanks to Nathanael Arkor, Andreas Blass, John Bourke,
Aaron Fairbanks, Marcelo Fiore, Richard Garner,
Theo Johnson-Freyd, Tom Leinster, Mat\'{i}as Menni, Deni
Salja, and Peter Selinger for their insightful comments
and interest. A special thanks to Peter Selinger for helping
me prepare the final version in the MSC style.


\section{Profunctors}

Profunctors (a.k.a. bimodules, modules, distributors) will be
at the heart of this work. Widely viewed as categorified  relations,
for our purposes they are better viewed as categorified matrices.
They correspond to cocontinuous functors between functor
categories. Such functors are considered to be linear. This section
contains nothing new (except perhaps Definition \ref{Def-Core} and 
Proposition \ref{Prop-Core}).  It is included for completeness and 
to set notation.

\subsection{Definitions}

We have opted, not without thought, for the following
definition which is the opposite of the majority view.

\begin{xdefinition}(Lawvere, B\'enabou)
\label{Def-ProfDef}
Let \( {\bf A} \) and \( {\bf B} \) be small categories. A
{\em profunctor} \( P \colon {\bf A} \longtod {\bf B} \) is
a functor \( P \colon {\bf A}^{op} \times {\bf B} \to {\bf Set} \). A
{\em morphism of profunctors} \( t \colon P \to Q \) is a
natural transformation.

\end{xdefinition}

This gives the basic data for a bicategory, \( {\cal P}\!{\it rof} \),
of profunctors. Composition is given by ``matrix multiplication''
which takes the form of a coend. For \( P \colon {\bf A} \longtod
{\bf B} \) and \( Q \colon {\bf B} \longtod {\bf C} \), the
composite \( Q \otimes P \) is defined by
\[
Q \otimes P (A, C) = \int^{B \in {\bf B}}
Q (B, C) \times P (A, B) \rlap{\,.}
\]
The identity \( \Id_{\bf A} \colon {\bf A} \longtod {\bf A} \) is the
hom functor
\[
\Id_{\bf A} (A, A') = {\bf A} (A, A') \rlap{\,.}
\]

The reader is referred to the standard texts (see e.g.~\citet{Bor94A}) for a proof
that we do get a bicategory.

For explicit computations involving profunctors, the
following notation is useful. An element \( x \in P (A, B) \)
is denoted by a pointed arrow, sometimes called a
heteromorphism, \( x \colon A \longtod B \), or
\( x \colon A \todd{}{P} B \) if it's necessary to keep track
of the profunctor. The functoriality of \( P \) manifests itself
as a composition
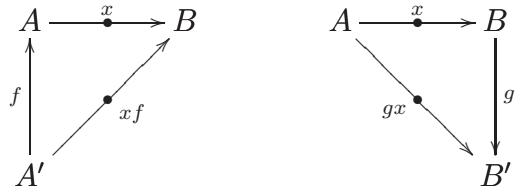
\[
\bfig
\ptriangle/@{>}|{\bb}`<-`@{<-}|{\bb}/[A`B`A';x`f`xf]

\qtriangle(1000,0)/@{>}|{\bb}`@{>}|{\bb}`>/[A`B`B';x`gx`g]

\efig
\]
which is associative (left, right, and middle) and unitary.

It is in dealing with composition that this is most useful. An
element of \( Q \otimes P (A, C) \) is an equivalence class of
pairs
\[
[A \todd{x}{P} B \todd{y}{Q} C]_B
\]
where the equivalence relation is generated by identifying
 \( [A \todd{x}{} B \todd{y}{} C] \) and\break \( [A \todd{x'}{}
 B \todd{y'}{} C] \) if we have
 \[
 \bfig
 \square/@{>}|{\bb}`=`>`@{>}|{\bb}/[A`B`A`B';x``b`x']
 
 \square(500,0)/@{>}|{\bb}``=`@{>}|{\bb}/[B`C`B'`C\rlap{\ ,};y```y']
 
 \efig
 \]
 so they are equivalent iff there exists a path of pairs
\begin{equation}\tag{*}
 \bfig
 \node a(0,0)[A]
 \node b(500,0)[B']
 \node c(1000,0)[C\rlap{\ .}]
 \node d(0,400)[\vdots]
 \node e(500,400)[\vdots]
 \node f(1000,400)[\vdots]
 \node g(0,800)[A]
 \node h(500,800)[B_2]
 \node i(1000,800)[C]
 \node j(0,1200)[A]
 \node k(500,1200)[B_1]
 \node l(1000,1200)[C]
 \node m(0,1600)[A]
 \node n(500,1600)[B]
 \node o(1000,1600)[C]

 \arrow/=/[d`a;]
 \arrow|r|/>/[b`e;b_n]
 \arrow/=/[c`f;]
 \arrow/=/[g`d;]
 \arrow|r|/>/[h`e;b_3]
 \arrow/=/[f`i;]
 \arrow/=/[g`j;]
 \arrow|r|/>/[h`k;b_2]
 \arrow/=/[i`l;]
 \arrow/=/[j`m;]
 \arrow|r|/>/[n`k;b_1]
 \arrow/=/[o`l;]
 \arrow|b|/@{>}|{\bb}/[a`b;x']
 \arrow|b|/>/[b`c;y']
 \arrow|a|/@{>}|{\bb}/[g`h;x_2]
 \arrow|a|/>/[h`i;y_2]
 \arrow|a|/@{>}|{\bb}/[j`k;x_1]
 \arrow|a|/@{>}|{\bb}/[k`l;y_1]
 \arrow|a|/@{>}|{\bb}/[m`n;x]
 \arrow|a|/@{>}|{\bb}/[n`o;y]

 \efig
 \end{equation}
We write the equivalence class \( [A \todd{x}{} B \todd{y}{} C]_B \)
as
\[
y \otimes_B x \mbox{\ \ or simply\ \ } y \otimes x \rlap{\,.}
\]
The equivalence relation is generated by
\[
y b \otimes x = y \otimes b x \rlap{\,.}
\]

Every functor \( F \colon {\bf A} \to {\bf B} \) induces two
profunctors
\[
F_* \colon {\bf A} \longtod {\bf B} \quad\quad F^* \colon {\bf B} \longtod {\bf A}
\]
\[
\mbox{and}
\]
\[F_* (A, B) = {\bf B} (FA, B) \quad\quad F^* (B, A) = {\bf B} (B, FA) .
\]
\( F^* \) is right adjoint to \( F_* \) in \( {\cal P}\!{\it rof} \).

\subsection{Biclosedness}

The bicategory \( \Prof \) is biclosed, that is \( \otimes \)
admits right adjoints in each variable giving two hom profunctors
\( \oslash \) and \( \obslash \) characterized by natural
bijections
\begin{center}
\begin{tabular}{c} 
\( P \to Q \obslash_{\bf C} R \) \\[3pt] \hline \\[-12pt]
\( Q \otimes_{\bf B} P \to R \) \\[3pt] \hline \\[-12pt]
\( Q \to R \oslash_{\bf A} P \)
\end{tabular}
\end{center}
for profunctors
\[
\bfig
\qtriangle/@{>}|{\bb}`@{>}|{\bb}`@{>}|{\bb}/<450,350>[{\bf A}`{\bf B}`{\bf C}\rlap{\,.};
P`R`Q]

\efig
\]

We use Lambek's notation for the internal homs. Inasmuch
as \( \otimes \) is a product, the right adjoints are quotients
of a sort.

An element of \( (Q \obslash_{\bf C} R) (A, B) \) is a
\( {\bf C} \)-natural transformation
\[
t \colon Q (B, -) \to R (A, -)
\]
and an element of \( (R \oslash_{\bf A} P) (B, C) \) is an
\( {\bf A} \)-natural transformation
\[
u \colon P (-, B) \to R (-, C)\rlap{\,.}
\]

\subsection{Cocontinuous functors}

Our interest is in functors between functor categories and a
profunctor will produce an adjoint pair of them. A profunctor
\( \colon {\bf 1} \longtod {\bf A} \) is a functor
\[
{\bf 1}^{op} \times {\bf A} \to {\bf Set}
\]
which we identify with a functor \( \Phi \colon {\bf A} \to {\bf Set} \).
A profunctor \( P \colon {\bf A} \longtod {\bf B} \) will then produce,
by composition, a functor
\[
P \otimes_{\bf A} (\ \ ) \colon {\bf Set}^{\bf A} \to {\bf Set}^{\bf B}
\]
with a right adjoint
\[
P \obslash_{\bf B} (\ \ ) \colon {\bf Set}^{\bf B} \to {\bf Set}^{\bf A} \rlap{\,.}
\]
It follows that \( P \otimes_{\bf A} (\ \ ) \) is cocontinuous and is
considered to be the linear functor corresponding to the matrix
\( P \).

As is well-known, we have:

\begin{xproposition}

The following categories are equivalent:
\begin{itemize}

\item[(1)] Profunctors \( {\bf A} \longtod {\bf B} \)

\item[(2)] Cocontinuous functors \( {\bf Set}^{\bf A} \to {\bf Set}^{\bf B} \)

\item[(3)] Adjoint pairs \( \bfig
\morphism(0,0)|a|/@{>}@<.5em>/<500,0>[{\bf Set}^{\bf A}`{\bf Set}^{\bf B};]
\morphism(0,0)|a|/@{<-}@<-.5em>/<500,0>[{\bf Set}^{\bf A}`{\bf Set}^{\bf B};\bot]
\efig \)
\end{itemize}

\end{xproposition}

Given a cocontinuous functor \( F \colon {\bf Set}^{\bf A} \to
{\bf Set}^{\bf B} \), the corresponding profunctor \sloppy \( P \colon {\bf A}
\longtod {\bf B} \) is given by
\[
P (A, B) = F ({\bf A} (A, -)) (B) \rlap{\,.}
\]
Note that this doesn't use cocontinuity of \( F \), which leads
to the following.

\begin{xdefinition}
\label{Def-Core}

The {\em core} of a functor \( F \colon {\bf Set}^{\bf A} \to
{\bf Set}^{\bf B} \) is the profunctor defined by
\[
\mbox{Cor} (F) (A, B) = F ({\bf A} (A, -)) (B) \rlap{\,.}
\]
The functor
\[
\mbox{Cor} (F) \otimes (\ \ ) \colon {\bf Set}^{\bf A} \to
{\bf Set}^{\bf B}
\]
is the ``linear core'' of \( F \).

\end{xdefinition}

\begin{xproposition}
\label{Prop-Core}

\( \mbox{Cor} \) is right adjoint to the functor
\( \Prof ({\bf A}, {\bf B}) \to \Cat ({\bf Set}^{\bf A}, {\bf Set}^{\bf B}) \)
which takes a profunctor \( P \) to the (cocontinuous) functor
\( P \otimes_{\bf A} (\ \ ) \).

\end{xproposition}

\begin{proof}

A profunctor \( P \colon {\bf A} \longtod {\bf B} \) can be viewed,
by exponential adjointness, as a functor \( {\bf A}^{op} \to
{\bf Set}^{\bf B} \). Then \( \mbox{Cor} \) is just restriction along
the Yoneda embedding
\[
F \colon {\bf Set}^{\bf A} \to {\bf Set}^{\bf B}\quad \longmapsto \quad
{\bf A}^{op} \to^Y {\bf Set}^{\bf A} \to^F {\bf Set}^{\bf B}
\]
and \( P \otimes_{\bf A} (\ \ ) \) is left Kan extension
\[
\bfig
\Vtriangle/>`>`>/<400,500>[{\bf A}^{op}`{\bf Set}^{\bf A}`{\bf Set}^{\bf B};
Y`P`\mbox{Lan}_Y P = P \otimes (\ \ )]

\morphism(320,300)/=>/<150,0>[`;]

\place(800,0)[.]

\efig
\] \end{proof}

Thus for \( F \colon {\bf Set}^{\bf A} \to {\bf Set}^{\bf B} \),
\( \mbox{Cor} (F) \otimes_A (\ \ ) \) is the best approximation
to \( F \) by a cocontinuous functor. As a matter of interest,
the counit of the adjunction
\[
\epsilon (F) \colon \mbox{Cor} (F) \otimes (\ \ ) \to F
\]
is given as follows. An element of \( (\mbox{Cor} (F)
\otimes \Phi )( B) \) is an equivalence class
\[
[x \in \Phi A, y \colon A \todd{}{\mbox{\scriptsize Cor}(F)} B]_A
\]
so
\[
[{\bf A} (A, -) \to^{\ov{x}} \Phi, \ y \in F ({\bf A} (A, -))(B)]
\]
giving
\[
F ({\bf A} (A, -))(B) \to^{F(\ov{x})(B)} F (\Phi) (B)
\]
\[
y \longmapsto F (\ov{x})(B) (y) \rlap{\,.}
\]

\begin{xexample}
\label{Ex-DiscProf}

If \( {\bf A} \) and \( {\bf B} \) are discrete categories, i.e.~sets
\( A \) and \( B \), then a profunctor \( P \colon {\bf A} \longtod {\bf B} \)
is just a \( A \times B \)-matrix of sets \( [P_{ab}] \) and a morphism
of profunctors \( P \to P' \) a \( A \times B \)-matrix of functions. The
identity \( \Id_{\bf A} \) is the matrix with \( 1 \)'s on the diagonal
and \( 0 \) elsewhere. If \( {\bf C} \) is another discrete category
and \( Q \colon {\bf B} \longtod {\bf C} \) a profunctor, then
\( Q \otimes_{\bf B} P \) is the \( B \times C \)-matrix

\[
\Bigl[ \sum_{b \in B} Q_{bc} \times P_{ab} \Bigr]_{\rlap{\,.}}
\]

\noindent If \( R \colon {\bf A} \longtod {\bf C} \) then

\[
R \oslash_{\bf A} P = \Bigl[\prod_{a \in A} R_{ac}^{P_{ab}} \Bigr]
\]

\noindent and

\[
Q \obslash_{\bf C} R = \Bigl[\prod_{c \in C} R_{ac}^{Q_{bc}}\Bigr]_{\rlap{\,.}} 
\]

A profunctor \( X \colon {\bf 1} \longtod {\bf A} \) is a \( 1 \times A \) matrix
of sets, i.e.~a vector \([X_a] \) and \( P \otimes_{\bf A} X \) is the vector

\[
\Bigl[ \sum_{a \in A} P_{ab} \times X_a\Bigr]_{b \rlap{\,.}} \]
On the other hand for \( Y \colon {\bf 1} \longtod {\bf B} \) a
\( B \)-vector \( P \obslash_{\bf B} Y \)
\[
\Bigl[\prod_b Y_b^{P_{ab}} \Bigr]_{a\rlap{\,.}}
\]
So \( P \otimes_{\bf A} (\ \ ) \) is a ``linear'' functor, and
\( P \obslash_{\bf B} (\ \ ) \) a ``monomial'' functor.

\end{xexample}


\section{Tense functors}

In \citet{Par24} we developed a difference calculus for taut endofunctors
of \( {\bf Set} \), functors preserving inverse images. However, the important
example of multivariable analytic functors of \citet{FioGamHylWin08} are not taut. In
fact the linear functors \( P \otimes (\ \ ) \) are not taut. They don't even
preserve monos. What we need are functors preserving complemented
subobjects and their inverse images. Of course, in \( {\bf Set} \), all
subobjects are complemented so it would make no difference, so maybe
that's what taut should be after all. But the word ``taut'' is pretty well established,
so we use ``tense'' instead.

\subsection{Complemented subobjects}

In this section we collect some useful facts about complemented
subobjects in functor categories \( {\bf Set}^{\bf A} \), most of which are
well-known from topos theory. We first list some general topos theory
results which will be useful for us. Proofs can be found in any of
the standard topos theory books (see \citet{Bor94} for an easily
accessible account).

\begin{xdefinition}
\label{Def-Compl}

A subobject \( \Psi \to/ >->/ \Phi \) is {\em complemented} if there exists
another subobject \( \Psi' \ \to/>->/ \Phi \) for which the induced
morphishm \( \Psi + \Psi' \to \Phi \) is invertible.

\end{xdefinition}

We will use the hooked arrow \( \Psi \hookrightarrow \Phi \)
as a reminder that \( \Psi \) is complemented,

Recall that every subobject \( \Psi \to/ >->/ \Phi \) has a pseudo-complement
\( \neg \Psi \to/ >->/ \Phi \), the largest subobject of \( \Phi \) whose
intersection with \( \Psi \) is \( 0 \). It can be calculated as the pullback
of the element \( \mbox{false} \colon 1 \to/ >->/ \Omega \) along the
characteristic morphism of \( \Psi \).

\begin{xproposition}
\label{Prop-ToposCompl}

1. A subobject \( \Psi \to/ >->/ \Phi \) is complemented iff its characteristic
morphism factors through \( 1 + 1 \)
\[
\bfig
\Vtriangle/>`-->`<-< /<350,400>[\Phi`\Omega`1+1;\chi_\Psi``\langle\mbox{true, false}\rangle]

\efig
\]

2. Complemented subobjects are closed under composition.

\vspace{3mm}

3. Complemented objects are stable under pullback: if
\( \Psi \hookrightarrow \Phi \) is complemented and
\( f \colon \Theta \to \Phi \), then \( \neg f^{-1} (\Psi) =
f^{-1} (\neg \Psi) \) and we have isomorphisms
\[
\bfig
\square|alra|/ >->`>`>` >->/<800,500>[f^{-1}(\Psi) + f^{-1}(\neg \Psi)`\Theta
`\Psi + (\neg \Psi)`\Psi\rlap{\ .};\cong`g + g'`f`\cong]

\efig
\]

\vspace{3mm}

4. If \( \Psi \hookrightarrow \Phi \) is complemented, its complement
is \( \neg \Psi \), so complements are unique when they exist.

\vspace{3mm}

5. Given an inverse image diagram (pullback)
\[
\bfig
\square/ >->`>`>` >->/[\Gamma`\Theta`\Psi`\Phi\rlap{\,,};
`g`f`]

\place(250,250)[\Pb]

\efig
\]
\( f \) restricts to
\[
\bfig
\square/ >->`-->`>` >->/[\neg \Gamma`\Theta`\neg \Psi`\Phi;
`\neg g`f`]

\efig
\]
and the resulting square is also a pullback.

\end{xproposition}

Complemented subobjects in functor categories \( {\bf Set}^{\bf A} \)
are better behaved than in general toposes. For example
\( \neg \Psi \to/ >->/ \Phi \) is always complemented for any subobject
\( \Psi \to/ >->/ \Phi \).

\begin{xproposition}
\label{Prop-FunCatCompl}

For \( {\bf Set}^{\bf A} \) we have

(1) \( \Psi \to/ >->/ \Phi \) is complemented iff for all
\( f \colon A \to A' \) and \( x \in \Phi A \) we have
\[
x \in \Psi A \quad \Longleftrightarrow \quad \Phi (f) (x) \in \Psi (A) \rlap{\,.}
\]
This is equivalent to saying that for all \( f \colon A \to A' \)
\[
\bfig
\square/ >->`>`>` >->/[\Psi A`\Phi A`\Psi A'`\Phi A';
`\Psi f`\Phi f`]

\efig
\]
is a pullback diagram. This in turn is equivalent to saying that
for all \( f \colon A \to A' \), every commutative square
\[
\bfig
\square/>`>`>` >->/<700,500>[{\bf A}(A',-)`{\bf A}(A,-)`\Psi`\Phi;
{\bf A}(f,-)```]

\morphism(700,500)/-->/<-700,-500>[{\bf A}(A,-)`\Psi;]

\efig
\]
has a unique fill-in making the bottom triangle commute,
i.e.~\( \Psi \to \Phi \) is orthogonal to every representable
transformation.

(2) For \( \Psi \to/ >->/ \Phi \),
\[
\neg \Psi (A) = \{a \in \Phi A\ |\ (\forall f \colon A \to A')
(\Phi (f) (a) \notin \Psi (A')\}
\]
and \( \neg \neg \Psi (A) \) consists of all elements, \( x \)
of \( \Phi (A) \) connected to an element \( x' \) of \( \Psi \)
by a zigzag of elements of \( \Phi \)
\[
A \to/<-/ A_1 \to A_2 \to/<-/ \cdots \to A_n \to/=/ A'
\]
\[
\Phi A \to/<-/ \Phi A_1 \to \Phi A_2 \to/<-/ \cdots \Phi A_n \to/<-</\  \Psi A'
\]
\[
x \  \ \to/<-|/ \    x_1\  \to/|->/ \ \  x_2 \to/<-|/ \cdots \to/|->/  x_n \to/=/ x'
\]

(3) For any \( \Psi \to/ >->/ \Phi \), \( \neg \Psi \) is
complemented and its complement is \( \neg \neg \Psi \)
which is the smallest complemented subobject of
\( \Phi \) containing \( \Psi \).

\end{xproposition}

Thus the class of complemented subobjects consists of
all transformations right orthogonal to the representable
transformations \( {\bf A} (f, -) \), suggesting that it may
be the \( {\cal M} \) part of a factorization system on
\( {\bf Set}^{\bf A} \), which is indeed the case.

For \( \Phi \) in \( {\bf Set}^{\bf A} \), let \( \sim \) be the
equivalence relation on the set of all elements of \( \Phi \)
generated by identifying \( x \in \Phi A \) with
\( \Phi f (x) \in \Phi A' \) for all \( f \colon A \to A' \). Thus
\( x \in \Phi A \sim x' \in \Phi A' \) if there exists a zigzag
path as in (2) above.
The set of equivalence classes is the set of components
of \( \Phi \), \( \pi_0 \Phi = \limr_A \Phi A \), and two
elements are equivalent if and only if they are in the
same component.

\begin{xdefinition}
\label{Def-PiZeroSurj}

A transformation \( t \colon \Psi \to \Phi \) is
\( \pi_0 \)-{\em surjective} if \( \pi_0 t \colon \pi_0 \Psi
\to \pi_0 \Phi \) is surjective.

\end{xdefinition}

Thus \( t \) is \( \pi_0 \)-surjective iff every element of
\( \Phi \) is connected by a zigzag path to an element
in the image of \( t \).

\begin{xproposition}
\label{Prop-PiZeroSurj}

\begin{itemize}

	\item[(1)] \( t , u\  \pi_0\mbox{-surjective}  \Rightarrow
	t u\  \pi_0\mbox{-surjective} \).
	
	\item[(2)] \( t u\ \pi_0\mbox{-surjective} \Rightarrow
	t\ \pi_0\mbox{-surjective} \).
	
	\item[(3)] Every \( t \) factors uniquely up to a unique
	isomorphism as a \( \pi_0\mbox{-surjective} \) followed
	by a complemented monomorphism.
	
	\item[(4)] The \( \pi_0\mbox{-surjective} \) transformations
	are left orthogonal to the complemented monos.

\end{itemize}

\end{xproposition}

\begin{proof}

(1) and (2) are obvious from the definition. For (3), let
\( t \colon \Psi \to \Phi \) be any transformation. Let
\( \Phi_0 A \subseteq \Phi A \) be the set of all
\( x \in \Phi A \) connected to an element in the image
of \( t \). \( \Phi_0 \) is easily seen to be a complemented
subfunctor of \( \Phi \), and is in fact the union of all
of the components of \( \Phi \) that contain an element
in the image of \( t \). Then \( t \) factors as
\[
\Psi \to/>/<250>^{t_0} \Phi_0\  \to/^(->/<200> \Phi
\]
and \( t_0 \) is \( \pi_0\mbox{-surjective} \) by construction.
This is our factorization. The uniqueness part will follow
from (4).

Consider a commutative square in \( {\bf Set}^{\bf A} \)
\[
\bfig
\square/>`>`>`^(->/[\Psi`\Phi`\ \Gamma\ `\ \Delta\ ;
t`r`s`m]

\efig
\]
where \( t \) is \( \pi_0\mbox{-surjective} \) and \( m \) is
a complemented mono. Any \( x \in \Phi A \) is connected
to some \( t (A') (y) \) for \( y \in \Psi A' \), so \( s (A) (x) \)
is connected to \( s (A') t(A') (y) = m (A') r (A') (y) \). As
\( m \) is complemented, this implies that \( s (A) (x) \)
is in \( \Gamma (A) \). This gives the diagonal fill-in
\( \delta \colon \Phi \to \Gamma \) such that \( m\  \delta = s \)
and \( \delta\  t = r \). \( \delta \) is unique as \( m \) is monic.
\end{proof}

These results tell us that we have a factorization system on
\( {\bf Set}^{\bf A} \) with \( {\cal E} \) the class of
\( \pi_0 \)-surjections and \( {\cal M} \) the class of
complemented monos. We call it the {\em Boolean factorization}.
Note that the class of \( \pi_0 \)-surjections is not stable
under pullback however. Consider morphisms
\( f_i \colon A_0 \to A_i , i = 1, 2 \) in \( {\bf A} \) and
consider the pullback
\[
\bfig
\square<750,500>[\Sigma`{\bf A} (A_1, -)`{\bf A} (A_2, -)`{\bf A} (A_0, -)\rlap{\ .};
``{\bf A} (f_1, -)`{\bf A} (f_2, -)]

\place(375,250)[\Pb]

\efig
\]
\( \Sigma(A) \) consists of pairs of morphisms \( (g_1, g_2) \)
such that
\[
\bfig
\square[A_0`A_1`A_2`A;f_1`f_2`g_1`g_2]

\efig
\]
commutes, which well may be empty for all \( A \). In that
case, taking \( \pi_0 \) of the above pullback gives
\[
\bfig
\square[0`1`1`1;```]

\efig
\]
showing that \( {\bf A} (g_1, -) \) is \( \pi_0 \)-surjective
but its pullback is not.

Nevertheless, it will be useful for us in Section~\ref{Sec-Newton}
where we will be particularly interested in transformations
defined on sums of representables. We record here the
following facts for use later.

A natural transformation
\[
t \colon \sum_{j \in J} {\bf A} (C_j, -) \to
\sum_{i \in I} {\bf A} (A_i, -)
\]
is determined by a function on the indices \( \alpha \colon
J \to I \) and a \( J \)-family of functions \( \langle f_j\rangle \),
\[
f_j \colon A_{\alpha (j)} \to C_j \rlap{\ .}
\]
Write \( t = \displaystyle{\sum_\alpha} {\bf A} (f_j, -) \).

\begin{xproposition}
\label{Prop-TransReps}

With \( t \), \( \alpha \), \( f_i \) as above we have
\begin{itemize}

	\item[(1)] \( t \) is a complemented mono if and only if
	\( \alpha \) is one-to-one and the \( f_j \) are isomorphims.
	
	\item[(2)] \( t \) is \( \pi_0 \)-surjective if and only if
	\( \alpha \) is onto.
	
	\item[(3)] For a general \( t \) given by \( (\alpha,
	\langle f_j \rangle) \) we get its Boolean factorization
	by  factoring \( \alpha \)
	\[
	\bfig
	\Vtriangle/>`>`<-_{)}/<250,350>[J`I`K;\alpha`\sigma`\mu]
	
	\efig
	\]
	and then taking
	\[
	\sum_{j \in J} {\bf A}(C_j, -)
	\to^{\sum_\sigma {\bf A} (f_k, -)}
	\sum_{k \in K} {\bf A}(A_k, -)
	\to^{\sum_\mu {\bf A}(1_{A_i}, -)}
	\sum_{i\in I} {\bf A}, (A_i, -) \rlap{\ .}
	\]

\end{itemize}

\end{xproposition}

It's implicit in (1), but may be worth mentioning explicitly,
that the complemented subobjects of
\( \sum_{i \in I} {\bf A} (A_i, -) \) are the subsums, i.e.~of the
form \( \sum_{k \in K} {\bf A} (A_k, -) \) for \( K \subseteq I \).
It is also clear from the fact that each hom functor
\( {\bf A} (A_i, -) \) is connected and complemented, so is one
of the components of \( \sum_{i \in I} {\bf A} (A_i, -) \), and
any complemented subfunctor is a union of components.


The following is well-known (see \citet{Bor94}, Example 7.2.4).

\begin{xproposition}

Every subobject in \( {\bf Set}^{\bf A} \) is complemented
(\( {\bf Set}^{\bf A} \)\,is boolean) if and only if \( {\bf A} \) is
a groupoid.

\end{xproposition}

We end this subsection with the following, which says that
limits and confluent colimits of complemented subobjects
are again complemented.

\begin{xproposition}
\label{Prop-LimColimCompl}

Let \( \Gamma \colon {\bf I} \to {\bf Set}^{\bf A} \) be a
diagram in \( {\bf Set}^{\bf A} \) and \( \Gamma_0 \to/ >->/ \Gamma \)
a subdiagram such that for every \( I \), \( \Gamma_0 (I)
\hookrightarrow \Gamma (I) \) is complemented, then

\vspace{3mm}

\noindent (1) \( \liml \Gamma_0 \to \liml \Gamma \) is a complemented
subobject.

\vspace{3mm}

\noindent If \( {\bf I} \) is confluent we also have that

\vspace{3mm}

\noindent (2) \( \limr \Gamma_0 \to \limr \Gamma \) is a complemented
subobject.

\end{xproposition}

\begin{proof}

(1) \( \Gamma_0 (I) \hookrightarrow \Gamma (I) \) is complemented
iff for every \( f \colon A \to A' \),
\[
\bfig
\square/ >->`>`>` >->/<700,500>[\Gamma_0 (I) (A)`\Gamma (I) (A)
`\Gamma_0 (I) (A')`\Gamma (I) (A');
`\Gamma_0 (I) (f)`\Gamma (I) (f)`]

\efig
\]
is a pullback (\ref{Prop-FunCatCompl} (1)). Limits of pullback
diagrams are pullbacks, and the result follows.

(2) Recall from \citet{Par24} that a category \( {\bf I} \) is
confluent if any span can be completed to a commutative
square, and that confluent colimits commute with inverse
image diagrams in \( {\bf Set} \). This gives (2) immediately.
\end{proof}

\begin{xcorollary}
\label{Cor-IntUnionCompl}

The intersection of an arbitrary family of complemented
subobjects in a presheaf category is again complemented. 
The same for union.

\end{xcorollary}

\begin{proof}

Let \( \Psi_i \hookrightarrow \Phi \) be a family of
complemented subobjects. Without loss of generality
we can assume that the total subobject
\( \Phi \hookrightarrow \Phi \) is contained in it so that
the indexing poset \( {\bf I} \) is connected. Then by the
previous proposition
\[
\liml \Psi_i \to \liml \Phi
\]
is a complemented mono. Because \( {\bf I} \) is
connected the limit of the constant diagram 
\( \liml \Phi \cong \Phi \), and the \( \liml \Psi_i \)
is \( \cap \Psi_i \hookrightarrow \Phi \). The lattice
of complemented subobjects of \( \Phi \) is
self-dual which implies the result for unions.
\end{proof}

Note that this result does not hold in an arbitrary
Grothendieck topos.

\subsection{Tense functors}

As mentioned above, the functors
\( P \otimes (\ \ ) \colon {\bf Set}^{\bf A} \to {\bf Set}^{\bf B} \)
arising from profunctors are not generally taut. In fact
they don't even preserve monos in general. This may not
be surprising if we consider the tensor product of modules
but one might have hoped that things would be better in
the simpler \( {\bf Set} \) case.

\begin{xexample}
\label{Ex-NotTaut}

For any epimorphism \( e \colon A \to A' \) in \( {\bf A} \),
the natural transformation \( {\bf A} (e, -) \colon
{\bf A} (A', -) \to {\bf A} (A, -) \) is a monomorphism. If, for
a profunctor \( P \colon {\bf A} \longtod {\bf B} \),
\( P \otimes (\ \ ) \colon {\bf Set}^{\bf A} \to {\bf Set}^{\bf B} \)
were to preserve monos, we would need that
\( P \otimes {\bf A} (e, -) \) be a mono, but \( P \otimes {\bf A}
(e, -) \) is
\[
P (e, -) \colon P (A', -) \to P (A, -) \rlap{\ .}
\]
So \( P (e, B) \colon P (A', B) \to P (A, B) \) would have to
be one-to-one for all \( B \), but that's hardly always the
case. The simplest example is when \( {\bf A} = {\bf 2} \)
and \( {\bf B} = {\bf 1} \). Then \( P (e, 0) \) is an {\em arbitrary}
function in \( {\bf Set} \) (\( e \) is the unique morphism
\( 0 \to 1 \), which is of course epi).

\end{xexample}

Now, the functors \( P \otimes (\ \ ) \) are  ``linear functors''
and any theory of functorial differences that doesn't apply to
them is seriously flawed. This leads to the main definition of
the section.

\begin{xdefinition}
\label{Def-Tense}

A functor \( F \colon {\bf Set}^{\bf A} \to {\bf Set}^{\bf B} \) is
{\em tense} if it preserves

\vspace{2mm}

\noindent (1) complemented subobjects, and

\vspace{2mm}

\noindent (2) inverse images (pullbacks) of complemented
subobjects.

\vspace{2mm}
A natural transformation is {\em tense} if the naturality
squares corresponding to complemented subobjects
are pullbacks.
\end{xdefinition}

Tense functors are closely related to, though incomparable
with, taut functors. For this reason we chose the word ``tense''
as an approximate synonym and homonym of ``taut''.

Any functor preserving binary coproducts is tense,
in particular \( P \otimes (\ \ ) \), which preserves all colimits,
is tense.
So Example~\ref{Ex-NotTaut} shows that tense does not
imply taut. On the other hand the functor
\[
{\bf Set} \to {\bf Set}^{\bf 2}
\]
\[
A \longmapsto (A \to 1)
\]
is taut (a right adjoint, so preserves all limits) but not tense:
any proper subset \( A \subsetneq B \) gives a non-complemented
subobject
\[
\bfig
\square/ >->`>`>`=/<400,400>[A`B`1`1\rlap{\,.};```]

\efig
\]

The following is obvious but worth stating explicitly.

\begin{xproposition}
\label{Prop-Tense2Cat}

Identities are tense and compositions of tense functors
are tense. Horizontal and vertical composition of tense
natural transformations are again tense, giving a
sub-2-category \( {\cal T}\!\!{\it ense} \) of the
\( 2 \)-category \( \Cat \) of categories.

\end{xproposition}

\begin{xproposition}
\label{Prop-TenseVsTaut}

For any functor \( F \colon {\bf Set}^{\bf A} \to {\bf Set}^{\bf B} \)
we have
\begin{itemize}

	\item[(1)] If \( {\bf Set}^{\bf A} \) is Boolean then tense implies taut

	\item[(2)] If \( {\bf Set}^{\bf B} \) is Boolean then taut implies tense

	\item[(3)] If \( F \) is taut then it is tense if and only if \( F \) applied
to the first injection \( j \colon 1 \to 1 + 1 \) is complemented.

\end{itemize}

\end{xproposition}

\begin{proof}

(1) and (2) are obvious as is the ``only if'' part of (3), so assume
\( F \) is taut and \( F(j) \) complemented. If \( \Psi \hookrightarrow \Phi \)
is complemented, its characteristic morphism factors through
\( 1 + 1 \to/ >->/ \Omega \) giving a pullback

\[
\bfig
\square/^(->`>`>`^(->/<400,400>[\ \Psi\ `\Phi`\ 1\ `1+1\rlap{\,,};```j]

\place(200,200)[\Pb]

\efig
\]
\( F \) of which is also a pullback, so \( F (\Psi)\  \to/^(->/<200>
F (\Phi) \) is complemented.
\end{proof}

Evaluation functors preserve tenseness but, contrary to tautness,
they don't jointly create it. However if we consider ``evaluating
at a morphism'' they do.

\begin{xproposition}
\label{Prop-EvalTense}

A functor \( F \colon {\bf Set}^{\bf A} \to {\bf Set}^{\bf B} \) is
tense if and only if
\begin{itemize}

	\item[(1)] for every \( B \) in \( {\bf B} \), \( ev_B F
	\colon {\bf Set}^{\bf A} \to {\bf Set} \) is tense, and
	
	\item[(2)] for every \( g \colon B \to B' \), \( ev_g F
	\colon ev_B F \to ev_{B'} F \) is a tense transformation.

\end{itemize}

Furthermore, a natural transformation \( t \colon F \to G \)
is tense if and only if \( ev_B t \) is tense for every \( B \).

\end{xproposition}

\begin{proof}

\( ev_B \colon {\bf Set}^{\bf B} \to {\bf Set} \) preserves
coproducts so is tense and thus \( ev_B F \) will be
tense if \( F \) is. To say that \( ev_g \colon ev_B \to
ev_{B'} \) is tense is to say that for every complemented
subobject \( \Psi\  \to/^(->/ \Phi \) we have a pullback
\[
\bfig
\square/^(->`>`>`^(->/[\ \Psi B\ `\Phi B`\ \Psi B'\ `\Phi B';
`\Psi g`\Phi g`]

\place(250,250)[\Pb]

\efig
\]
which Proposition~\ref{Prop-FunCatCompl} (1) says is
indeed the case. So \( ev_g F \) will be tense when \( F \)
is.

In fact, this says that being complemented is
equivalent to every \( g \) giving a pullback as above. 
So our condition (2) implies that \( F \) preserves
complemented subobjects. And the evaluation functors
\( ev_B \) jointly create pullbacks. So (1) and (2) together
imply that \( F \) is tense.

The second part is clear as the functors \( ev_B \)
jointly create pullbacks and tenseness of natural
transformations is a purely pullback condition.
\end{proof}

\begin{xcorollary}
\label{Cor-EvalTense}

The following are equivalent.

\begin{itemize}
	\item[(1)] \( F \colon {\bf Set}^{\bf A} \to {\bf Set}^{\bf B} \) is tense.
	
	\item[(2)(a)] For every complemented subobject \( \Psi \hookrightarrow \Phi \)
	and every morphism \( g \colon B \to B' \),
	\[
	\bfig
	\square<650,500>[F(\Psi)(B)`F(\Phi)(B)`F(\Psi)(B')`F(\Phi)(B');```]
	
	\efig
	\]
	is a pullback diagram, and
	
	\item[{\phantom{(2)}}(b)] For every pullback diagram of complemented
	subobjects
	\[
	\bfig
	\square/^(->`>`>`^(->/<450,450>[\ \Psi'\ `\Phi'`\ \Psi\ `\Phi;```]
	
	\place(225,225)[\Pb]
	
	\efig
	\]
	and every \( B \) in \( {\bf B} \),
	\[
	\bfig
	\square<650,500>[F(\Psi')(B)`F(\Phi')(B)`F(\Psi)(B)`F(\Phi)(B);```]
	
	\efig
	\]
	is a pullback.
	
	\item[(3)] For every pullback diagram of complemented
	subobjects
	\[
	\bfig
	\square/^(->`>`>`^(->/<450,450>[\ \Psi'\ `\Phi'`\ \Psi\ `\Phi;```]
	
	\place(225,225)[\Pb]
	
	\efig
	\]
	and every \(g \colon B \to B' \),
	\[
	\bfig
	\square<650,500>[F(\Psi')(B)`F(\Phi')(B)`F(\Psi)(B')`F(\Phi)(B');```]
	
	\efig
	\]
	is a pullback.
	
\end{itemize}

Furthermore, \( t \colon F \to \ G \) is tense if and only if for every
complemented subobject \( \Psi \hookrightarrow \Phi \) and
every object \( B \) in \( {\bf B} \),
\[
\bfig
\square<650,500>[F(\Psi)(B)`F(\Phi)(B)`G(\Psi)(B')`G(\Phi)(B);
`t(\Psi)(B)`t(\Phi)(B)`]
	
\efig
\]
is a pullback.

\end{xcorollary}

\begin{proof}

That (1) is equivalent to (2) follows immediately from the previous
proposition, the definition of tense, and Proposition~\ref{Prop-FunCatCompl},
as does the statement about tense transformations.

(2) (a) and (b) are special cases of (3) and the pullback in (3) can be
factored into two pullbacks of type (a) and (b).
\end{proof}

\subsection{Limits and colimits of tense functors}
\label{SSec-LimColimTense}

\begin{xproposition}
\label{Prop-LimColimTense}

Let \( \Gamma \colon {\bf I} \to \Cat ({\bf Set}^{\bf A}, {\bf Set}^{\bf B}) \)
be a diagram such that for every \( I \) in \( {\bf I} \), \( \Gamma (I)  \) is
tense. Then

\vspace{1mm}

\noindent (1) \( \liml \Gamma \) is tense.

\vspace{1mm}

If \( t \colon \Gamma \to \Theta \) is a natural transformation
such that for every \( I \) in \( {\bf I} \), \( t I \colon \Gamma I \to
\Theta I \) is tense, then

\vspace{1mm}

\noindent (2) the induced transformation
\[
\liml t \colon \liml \Gamma \to \liml \Theta
\]
is tense.

If \( {\bf I} \) is confluent, then under the same conditions
as above we have

\vspace{1mm}

\noindent (3) \( \limr \Gamma \) is tense, and

\vspace{1mm}

\noindent (4) \( \limr t \) is tense.

\end{xproposition}

\begin{proof}

(1) and (3). The preservation of complemented subobjects
follows immediately from Proposition~\ref{Prop-LimColimCompl}.
The preservation of pullbacks of complemented subobjects
follows from the fact that limits commute with limits for (1)
and that confluent colimits commute with inverse images
for (3).

Tenseness of natural transformations is also a pullback
condition, so (2) and (4) follow for the same reasons.
\end{proof}

This is a result about limits and colimits of tense functors
taken in \( \Cat ({\bf Set}^{\bf A}, {\bf Set}^{\bf B}) \). It is
not assumed that the transition transformations
\( \Gamma (I) \to \Gamma (J) \) are tense, and
unsurprisingly we don't get a universal property for tense
cones or cocones. Given a tense cone or cocone, the
uniquely induced natural transformation is tense but
this doesn't establish the required bijection because
neither the projections in the limit case nor the injections
in the colimit case are tense.

It's more natural to consider diagrams where the
transitions {\em are} tense, i.e.~\( \Gamma \colon
{\bf I} \to {\cal T}\!\!{\it ense} ({\bf Set}^{\bf A}, {\bf Set}^{\bf B}) \).
For such diagrams, things are better. We lose products
as the projections are not tense but that's the only
obstruction. Limits of connected tense diagrams are
created by the inclusion
\[
{\cal T}\!\!{\it ense} ({\bf Set}^{\bf A}, {\bf Set}^{\bf B})
\to/ >->/<200>\  \Cat ({\bf Set}^{\bf A}, {\bf Set}^{\bf B})
\]
as are all colimits, not just confluent ones.

First we analyze diagrams \( \Gamma \colon
{\bf I} \to {\cal T}\!\!{\it ense} ({\bf Set}^{\bf A}, {\bf Set}^{\bf B}) \).

\begin{xproposition}
\label{Prop-TenseCoten}

The bicategory \( {\cal T}\!\!{\it ense} \) is \( \Cat \)-cotensored.
The cotensor of \( {\bf Set}^{\bf B} \) by \( {\bf I} \) is
\( {\bf Set}^{{\bf B} \times {\bf I}} \), i.e.
\begin{itemize}

	\item[(1)] diagrams \( \Gamma \colon
	{\bf I} \to {\cal T}\!\!{\it ense} ({\bf Set}^{\bf A}, {\bf Set}^{\bf B}) \)
	are in bijection with tense functors
	\( \overline{\Gamma} \colon {\bf Set}^{\bf A} \to
	{\bf Set}^{{\bf B} \times {\bf I}} \), and
	
	\item[(2)] natural transformations \( t \colon \Gamma \to\  \Theta \)
	are in bijection with tense natural transformations
	\( \overline{t} \colon \overline{\Gamma} \to \overline{\Theta} \).

\end{itemize}

\end{xproposition}

\begin{proof}

Functors \( \Gamma \colon{\bf I} \to \Cat ({\bf Set}^{\bf A},
{\bf Set}^{\bf B}) \) correspond bijectively to functors
\( \overline{\Gamma} \colon {\bf Set}^{\bf A} \to
{\bf Set}^{{\bf B} \times {\bf I}} \) by exponential adjointness:
\[
\overline{\Gamma} (\Phi) (B, I) = \Gamma (I) (\Phi) (B) \rlap{\,.}
\]
If \( \Gamma \) factors through \( {\cal T}\!\!{\it ense}
({\bf Set}^{\bf A}, {\bf Set}^{\bf B}) \) then we want to show
that \( \overline{\Gamma} \) is tense.

First of all \( \overline{\Gamma} (\Psi) \to \overline{\Gamma} (\Phi) \)
must be a complemented subobject for \( \Psi \ \to/^(->/ \Phi \)
complemented, i.e.
\[
\bfig
\square<700,500>[\overline{\Gamma}(\Psi)(B, I)
`\overline{\Gamma} (\Phi) (B, I)
`\overline{\Gamma} (\Psi)(B', I')
`\overline{\Gamma} (\Phi) (B', I');
`\overline{\Gamma} (\Psi)(g, \alpha)
`\overline{\Gamma}(\Phi)(g, \alpha)
`]

\efig
\]
for \( g \colon B \to B' \) and \( \alpha \colon I \to I' \),
should be a pullback of monos. If we rewrite this in terms
of \( \Gamma \) and use functoriality on the vertical
arrows we see that it is
\[
\bfig
\square/ >->`>`>` >->/<800,500>[\Gamma(I)(\Psi)(B')
`\Gamma (I) (\Phi)(B')
`\Gamma (I')(\Psi)(B')
`\Gamma (I')(\Phi)(B');
`\Gamma(\alpha)(\Psi)(B')
`\Gamma(\alpha)(\Phi)(B')`]

\place(400,750)[\framebox{{\scriptsize(1)}}]

\square(0,500)/ >->`>`>` >->/<800,500>[\Gamma (I)(\Psi)(B)
`\Gamma (I)(\Phi)(B)
`\Gamma (I)(\Psi)(B')
`\Gamma (I)(\Phi)(B');
`\Gamma (I)(\Psi)(g)
`\Gamma (I)(\Phi)(g)`]

\place(400,250)[\framebox{{\scriptsize(2)}}]

\efig
\]
(1) is a pullback of monos because \( \Gamma (I) \) is
tense, and (2) is a pullback of monos because
\( \Gamma (\alpha) \) is a tense transformation (the mono
part because \( \Gamma (I') \) is tense).

This shows that if \( \Gamma (I) \) preserves complemented
subobjects and \( \Gamma (\alpha) \) is tense, then
\( \overline{\Gamma} \) preserves complemented
subobjects. The converse is also true as can be seen by
taking \( \alpha = \id_I \) for \( \Gamma (I) \) and
\( g = 1_B \) for \( \Gamma (\alpha) \).

Preservation of inverse images by \( \overline{\Gamma} \)
is equivalent to that of \( \Gamma (I) \) as can be seen
immediately upon writing it down. Likewise for the tenseness
of \( \overline{t} \).
\end{proof}

\begin{xtheorem}
\label{Thm-LimColimTense}

The inclusion \( {\cal T}\!\!{\it ense} ({\bf Set}^{\bf A}, 
{\bf Set}^{\bf B}) \to/ >->/\ \Cat ({\bf Set}^{\bf A}, {\bf Set}^{\bf B}) \)
creates colimits and connected limits.

\end{xtheorem}

\begin{proof}

Given a diagram \( \Gamma \colon {\bf I} \to {\cal T}\!\!{\it ense}
({\bf Set}^{\bf A}, {\bf Set}^{\bf B}) \), its colimit is given by the
composite
\[
{\bf Set}^{\bf A} \to^{\overline{\Gamma}}
{\bf Set}^{{\bf B} \times {\bf I}}
\to^{\limr_{\bf I}} {\bf Set}^{\bf B}
\]
\( \limr_{\bf I} \) is left adjoint to the diagonal functor
\( D \colon {\bf Set}^{\bf B} \to {\bf Set}^{{\bf B} \times {\bf I}} \),
so it preserves coproducts and {\em a fortiori} is tense.
And \( \overline{\Gamma} \) is tense by the previous
proposition, so \( \limr_I \Gamma (I) \) is tense.

\( D \) itself preserves coproducts being left adjoint to
\( \liml \), the limit functor. So \( D \) is tense. Natural
transformations between coproduct preserving functors
are automatically tense, so the adjunction
\( \limr_{\bf I} \dashv D \) is an adjunction in the
bicategory \( {\cal T}\!\!{\it ense} \), and this gives the
universal property of \( \limr_I \Gamma (I) \):
\[
{\cal T}\!\!{\it ense} ({\bf Set}^{\bf A}, {\bf Set}^{\bf B})^{\bf I}
\to^\cong {\cal T}\!\!{\it ense}({\bf Set}^{\bf A}, {\bf Set}^{{\bf B} \times {\bf I}})
\to^{{\cal T}\!\!{\it ense} ({\bf Set}^{\bf A}, \ \limr)}
{\cal T}\!\!{\it ense} ({\bf Set}^{\bf A}, {\bf Set}^{\bf B})
\]
is left adjoint to
\[
{\cal T}\!\!{\it ense} ({\bf Set}^{\bf A}, {\bf Set}^{\bf B})
\to^{{\cal T}\!\!{\it ense} ({\bf Set}^{\bf A}, \ D)}
{\cal T}\!\!{\it ense}({\bf Set}^{\bf A}, {\bf Set}^{{\bf B} \times {\bf I}})
\to^\cong
{\cal T}\!\!{\it ense} ({\bf Set}^{\bf A}, {\bf Set}^{\bf B})^{\bf I}
\]
which is itself the diagonal functor.

If \( {\bf I} \) is non-empty and connected, then \( \liml_{\bf I} \colon
{\bf Set}^{{\bf B} \times {\bf I}} \to {\bf Set}^{\bf B} \) preserves
coproducts, so the same argument as above shows that
\( {\bf I} \)-limits are created in this case.
\end{proof}

\subsection{Internal homs}
\label{SSec-IntHom}

Part of the motivation for introducing tense functors was that
the functors \( P \otimes (\ \ ) \), thought of as linear, were not in general
taut but preserved coproducts, so were tense. The other
side of the story is that the right adjoint to \( P \otimes (\ \ ) \),
namely \( P \obslash (\ \ ) \), is taut but not always tense.
As Example~\ref{Ex-DiscProf} suggests \( P \obslash (\ \ ) \)
is a functorial version of a monomial with the \( P \) acting
as the powers, and perhaps we shouldn't expect them
to be nice for all \( P \). After all, even for real valued
functions, fractional powers can be problematic, and for
rings the powers are taken to be integers, not elements
of the ring.

\begin{xproposition}
\label{Prop-HomTense}

For a profunctor \( P \colon {\bf A} \longtod {\bf B} \)
the internal hom functor \( P \obslash (\ \ ) \colon
{\bf Set}^{\bf B} \to {\bf Set}^{\bf A} \) is tense if and only if
for every \( f \colon A \to A' \), the function
\[
\pi_0 P (A', -) \to \pi_0 P (A, -)
\]
is onto.

\end{xproposition}

\begin{proof}

\( P \obslash (\ \ ) \) preserves limits and so is taut. Thus
by Proposition~\ref{Prop-TenseVsTaut} (3) it is only
necessary to check that
\[
1 \cong P \obslash 1 \to P \obslash (1 + 1)
\]
is complemented, and it's also sufficient. This is equivalent
to the condition, that for every \( f \colon A \to A' \)
\[
\bfig
\square<800,500>[1`{\bf Set}^{\bf B} (P (A, -), 1+1)
`1`{\bf Set}^{\bf B} (P (A', -), 1+1);```]

\efig
\]
be a pullback. This says that every natural transformation
\( t \) for which (the outside of)
\[
\bfig
\square/>`<-`<-`>/<600,500>[P(A, -)`1+1`P(A', -)`1;
t`P (f, -)`j`]

\morphism(0,500)/-->/<600,-500>[P(A, -)`1;\exists]

\efig
\]
commutes, factors through the injection \( j \). This
is in \( {\bf Set}^{\bf B} \). Using the adjunction
\( \pi_0 \dashv Const \colon {\bf Set} \to {\bf Set}^{\bf B} \),
we have, equivalently, that every function \( \overline{t} \)
for which
\[
\bfig
\square/>`<-`<-`>/<600,500>[\pi_0 P(A, -)`1+1`\pi_0 P(A', -)`1;
\overline{t}``j`]

\morphism(0,500)/-->/<600,-500>[\pi_0 P(A, -)`1;]

\efig
\]
commutes, factors through \( j \) (in \( {\bf Set} \)). This is
equivalent to
\[
\pi_0 P (A', -) \to \pi_0 P (A, -)
\]
being onto.
\end{proof}

The condition on \( P \) making \( P \obslash (\ \ ) \) tense is
a kind of lifting condition. For every element of \( P \),
\( p \colon A \longtod B \) and morphism \( f \colon A \to A' \)
there exist a \( B' \) and a \( P \)-element \( p' \colon A'
\longtod B' \) for which \( p' f \) is connected to \( P \) by
a path of \( P \)-elements
\[
\bfig
\node a(0,0)[B']
\node b(350,0)[B_1]
\node c(700,0)[B_2]
\node d(1050,0)[B_3]
\node e(1400,0)[]
\node g(1600,0)[\cdots]
\node f(1750,0)[]
\node h(2000,0)[B]

\node a'(0,400)[A']
\node b'(350,400)[A]
\node c'(700,400)[A]
\node d'(1050,400)[A]
\node e'(1400,400)[]
\node g'(1600,400)[\cdots]
\node f'(1750,400)[]
\node h'(2000,400)[A]

\arrow/<-/[a`b;]
\arrow/>/[b`c;]
\arrow/<-/[c`d;]
\arrow/>/[d`e;]
\arrow/>/[f`h;]

\arrow/<-/[a'`b';f]
\arrow/=/[b'`c';]
\arrow/=/[c'`d';]
\arrow/=/[d'`e';]
\arrow/=/[f'`h';]

\arrow|l|/@{>}|{\bb}/[a'`a;p']
\arrow|r|/@{>}|{\bb}/[b'`b;p_1]
\arrow|r|/@{>}|{\bb}/[c'`c;p_2]
\arrow|r|/@{>}|{\bb}/[d'`d;p_3]
\arrow|r|/@{>}|{\bb}/[h'`h;p]

\place(2100,0)[.]

\efig
\]

Or more fancifully and more memorably, it's a kind
of homotopy pushout condition: for every \( f \) and
\( p \) as below there exist a lifting to a \( p' \) with a
fill in ``fan''
\[
\bfig
\node a(500,0)[B]
\node b(700,200)[B_n]
\node c(750,260)[.]
\node d(780,280)[.]
\node e(800,300)[.]
\node f(850,350)[B_2]
\node g(1050,550)[B_1]
\node h(1240,750)[B']

\node i(-240,650)[A]
\node j(500,1400)[A']

\arrow/>/[b`a;]
\arrow/>/[g`f;]
\arrow/>/[g`h;]

\arrow|l|/>/[i`j;f]
\arrow|l|/@{>}|{\bb}/[i`a;p]
\arrow|r|/@{>}|{\bb}/[j`h;p']

\arrow/@{>}|{\bb}/[i`b;]
\arrow/@{>}|{\bb}/[i`f;]
\arrow/@{>}|{\bb}/[i`g;]

\efig
\]

\subsection{Multivariable analytic functors}
\label{SSec-MultVarAndFun}

Following \citet{FioGamHylWin08} we define
analytic functors of several variables \( F \colon
{\bf Set}^{\bf A} \to {\bf Set}^{\bf B} \) as follows. First,
for a category \( {\bf A} \), its {\em exponential}
\( ! {\bf A} \) (from linear logic) is the free
symmetric strict monoidal category generated by \( {\bf A} \).
In concrete terms, \( ! {\bf A} \) is the category with
objects finite sequences \( \langle A_1 \dots A_n\rangle \)
of objects of \( {\bf A} \) and morphisms finite
sequences of morphisms of \( {\bf A} \) controlled by
a permutation. There are no morphisms between
sequences unless they have the same length and then
\[
\langle A_1 \dots A_n\rangle \to \langle A'_1 \dots , A'_n\rangle
\]
is a permutation of the indices, \( \sigma \in S_n \) and a
sequence of morphisms
\[
f_i \colon A_{\sigma i} \to A'_i \rlap{\,.}
\]
Composition is as expected
\[
(\tau, \langle g_i \rangle) (\sigma, \langle f_i \rangle) =
(\sigma \tau, \langle g_i f_{\tau i} \rangle).
\]

An \( {\bf A} \)-\( {\bf B} \) {\em symmetric sequence} is a
profunctor \( P \colon ! {\bf A} \longtod {\bf B} \), which for us is a
functor \( (! {\bf A})^{op} \times {\bf B} \to {\bf Set} \).
({\em Warning}: Our definition of profunctor is the
opposite of theirs.) \( P \) encodes what are to be the
coefficients of a \( {\bf B} \)-family of multivariable
power series.

The {\em analytic functor} determined by \( P \)
\[
\widetilde{P} \colon {\bf Set}^{\bf A} \to {\bf Set}^{\bf B}
\]
is given by
\[
\widetilde{P} (\Phi) (B) = \int^{\langle A_1 \dots A_n\rangle \in\, ! {\bf A}}
P (A_1 \dots A_n; B) \times \Phi A_1 \times \Phi A_2 \times \dots
\times \Phi A_n \rlap{\,.}
\]

We'll show that \( \widetilde{P} \) is tense. Define a
profunctor \( Q \colon ! {\bf A} \longtod {\bf A} \) by
\[
Q (A_1, \dots, A_n ; A) = {\bf A} (A_1, A) + {\bf A} (A_2, A)
+ \dots + {\bf A} (A_n, A)
\]
with the obvious definition on morphisms. We may consider
\( Q \) as a functor \( (!{\bf A})^{op} \to {\bf Set}^{\bf A} \) and
\( \widetilde{P} \) is the left Kan extension of \( P \), considered
as a functor \( (!{\bf A})^{op} \to {\bf Set}^{\bf B} \), along \( Q \)
\[
\bfig
\Vtriangle<350,500>[(!{\bf A})^{op}`{\bf Set}^{\bf A}`{\bf Set}^{\bf B}\rlap{\ .};
Q`P`\widetilde{P} = {\rm Lan}_Q P]

\morphism(290,300)/=>/<150,0>[`;]

\efig
\]
For our purposes a different description of \( \widetilde{P} \)
will be useful.


\begin{xproposition}
\label{Prop-AnalComposite}

1. \( \widetilde{P} \) is the composite \( P \otimes
(Q \obslash (\ \ )) \)
\[
{\bf Set}^{\bf A} \to^{Q \obslash (\ \ )} {\bf Set}^{! {\bf A}}
\to^{P \otimes (\ \ )} {\bf Set}^{\bf B}\rlap{\,.}
\]

\noindent 2. \( Q \) satisfies the condition of Proposition~\ref{Prop-HomTense}.

\end{xproposition}

\begin{proof}

(1) Let \( \Phi \in {\bf Set}^{\bf A} \). An element of
\( (Q \obslash \Phi) (A_1, \dots, A_n) \) is a natural
transformation
\[
{\bf A} (A_1, -) + \dots + {\bf A} (A_n, -) \to \Phi
\]
which by the universal property of coproduct and the
Yoneda lemma corresponds to an element of
\[
\Phi A_1 \times \Phi A_2 \times \dots \times \Phi A_n \rlap{\,.}
\]
Now the result follows by the definition of \( P \otimes (\ \ ) \)
and \( \widetilde{P} \).

\noindent (2) \( Q (A_1, \dots , A_n ; -) =
{\bf A} (A_1, -) + \dots + {\bf A}(A_n, -) \) a sum of
representables each of which is connected. So
\[
\pi_0 Q (A_1, \dots, A_n; -) \cong n
\]
and, as \( ! {\bf A} \) has only morphisms between
sequences of the same length, we get
\[\pi_0 Q (A_1, \dots, A_n ; -) \cong \pi_0 Q
(A'_1 , \dots , A'_n; -) \rlap{\,.}
\]
\end{proof}

\begin{xcorollary}

\( \widetilde{P} \) is tense.

\end{xcorollary}

\begin{xcorollary}
\label{Cor-TildeLin}

For \( P \colon ! {\bf A} \longtod {\bf B} \) an \( {\bf A} \)-\( {\bf B} \)
symmetric sequence and \( R \colon {\bf B} \longtod {\bf C} \) a
profunctor, we have
\[
\widetilde{R \otimes P} \cong R \otimes \widetilde{P} \rlap{\,.}
\]

\end{xcorollary}

\begin{proof}
\[
\begin{array}{lll}
\widetilde{R \otimes P} & \cong & (R \otimes P) \otimes (Q \obslash (\ ))\\
     & \cong & R \otimes (P \otimes (Q \obslash (\ ))\\
     & \cong & R \otimes \widetilde{P} \rlap{\,.}
\end{array}
\]
\end{proof}


\section{Partial difference operators}

We want to think of a functor \( F \colon {\bf Set}^{\bf A} \to
{\bf Set}^{\bf B} \) as a \( {\bf B} \)-family of \( {\bf Set} \)-valued
functors in \( {\bf A} \)-variables and study its change under
small perturbations of the variables. The context is that of
tense functors and for these we get a theory that parallels the
usual calculus of differences for real-valued functions of
several variables, much as our theory for taut functors did
for single variables \citep{Par24}.

\subsection{Partial difference}

A functor \( \Phi \in {\bf Set}^{\bf A} \) is a multisorted algebra,
the sorts being the objects of \( {\bf A} \), with unary operations
corresponding to the morphisms of \( {\bf A} \). Freely adding
a single element of sort \( A \) gives
\[
\Phi \leadsto \Phi + {\bf A} (A, -) \rlap{\,.}
\]

\begin{xdefinition}
\label{Def-AShift}

The {\em \( A \)-shift} functor, for an object \( A \) in \( {\bf A} \)
is
\[
S_A \colon {\bf Set}^{\bf A} \to {\bf Set}^{\bf A}
\]
\[
S_A (\Phi) = \Phi + {\bf A} (A, -)\rlap{\,.}
\]

\end{xdefinition}

\( S_A \) is clearly tense, in fact a tense monad. Although
we won't use it here, it may be of interest to note that an
Eilenberg-Moore algebra for \( S_A \) consists of a functor
\( \Phi \in {\bf Set}^{\bf A} \) together with an element
\( x \in \Phi A \). A Kleisli morphism \( \Phi \longtod \Psi \)
is a partial natural transformation
\[
\bfig
\Atriangle/ >->`>`/<350,300>[\Phi_0`\Phi`\Psi;``]

\efig
\]
defined on a complemented subobject \( \Phi_0 \) together
with a transformation on the complement
\( \Phi'_0 \to {\bf A} (A, -) \), perhaps quantifying the degree
of undefinedness.

These monads commute with each other
\[
S_{A_1} \circ S_{A_2} \cong S_{A_2} \circ S_{A_1}
\]
and for every \( f \colon A \to A' \) there is a monad morphism
\( S_A \to S_{A'} \) which is tense.

The main definition of the paper is the following.

\begin{xdefinition}
\label{Def-APartDiff}

The {\em partial difference with respect to} \( A \), or the
\( A \)-{\em partial difference}, \( \Delta_A [F] \), of a tense
functor \( F \colon {\bf Set}^{\bf A} \to {\bf Set}^{\bf B} \)
is given by
\[
\Delta_A [F] \colon {\bf Set}^{\bf A} \to {\bf Set}^{\bf B}
\]
\[
\Delta_A [F] (\Phi) = F (\Phi + {\bf A} (A, -)) \setminus F (\Phi)\,,
\]
the complement of \( F (\Phi) \hookrightarrow F (\Phi + {\bf A} (A, -)) \).

\end{xdefinition}

\begin{xproposition}
\label{Prop-PartTense}

For a tense functor \( F \colon {\bf Set}^{\bf A} \to {\bf Set}^{\bf B} \),
\( \Delta_A [F] \) is also a tense functor. A tense natural
transformation \( t \colon F \to G \) restricts to one,
\( \Delta_A [t] \colon \Delta_A [F] \to \Delta_A [G] \), making
\( \Delta_A \) a functor
\[
\Delta_A \colon {\cal T}\!\!{\it ense} ({\bf Set}^{\bf A}, {\bf Set}^{\bf B})
\to {\cal T}\!\!{\it ense} ({\bf Set}^{\bf A}, {\bf Set}^{\bf B}) \,.\]

\end{xproposition}

\begin{proof}

Let \( \phi \colon \Psi \to \Phi \) be a natural transformation. We have
the following pullbacks
\[
\bfig
\square/^(->`>`>`^(->/<600,500>[\ \Psi\ `\Psi + {\bf A}(A, -)`\ \Phi\ `
\Phi + {\bf A}(A, -);`\phi`\phi + {\bf A}(A, -)`]

\place(300,250)[\Pb]

\square(1500,0)/^(->`>`>`^(->/<700,500>[\ \ F \Psi\ `F(\Psi + {\bf A}(A, -))`\ F\Phi\ `
F(\Phi + {\bf A}(A, -))\rlap{\,.};`F \phi`F(\phi + {\bf A}(A, -))`]

\place(1850,250)[\Pb]

\efig
\]
From the second one we get that \( F (\phi + {\bf A} (A, -)) \) restricts
to the complements and gives another pullback
\[
\bfig
\square/^(->`>`>`^(->/<850,500>[\ \Delta_A {[}F{]}(\Psi)\ `F(\Psi + {\bf A} (A, -))
`\ \Delta_A {[}F{]} (\Phi)\ `F (\Phi + {\bf A} (A, -));
`\Delta_A {[}F{]}(\phi)`F(\phi + {\bf A} (A,-))`]

\place(425,250)[\Pb]

\efig
\]
which gives functoriality and tenseness.

Suppose \( t \colon F \to G \) is a tense transformation. Then
we get a pullback for any \( \Phi \)
\[
\bfig
\square/^(->`>`>`^(->/<700,500>[\ F(\Phi)\ `F(\Phi + {\bf A} (A, -))
`\ G(\Phi)\ `G(\Phi + {\bf A} (A, -));
`t(\Phi)`t(\Phi + {\bf A}(A, -))`]

\place(350,250)[\Pb]

\efig
\]
so \( t (\Phi + {\bf A} (A,-)) \) restricts to the complements,
giving another pullback
\[
\bfig
\square/^(->`>`>`^(->/<800,500>[\ \Delta_A {[}F{]}(\Phi)\ 
`F (\Phi + {\bf A} (A,-))`\ \Delta_A{[}G{]}(\Phi)\ 
`G(\Phi + {\bf A} (A,-))\rlap{\,.};
`\Delta_A{[}t{]} (\Phi)`t(\Phi + {\bf A}(A,-))`]

\place(400,250)[\Pb]

\efig
\]
It follows immediately that \( \Delta_A [t] \) is natural. Tenseness
follows by comparing the following diagrams that we get for
any complemented subobject \( \Psi\  \to/^(->/ \Phi \).
\[
\bfig
\square/^(->`>`>`^(->/<700,500>[\ \Delta_A{[}F{]}(\Psi)\ 
`\ \Delta_A {[}F{]}(\Phi)\ 
`\ \Delta_A{[}G{]}(\Psi)\ 
`\ \Delta_A{[}G{]}(\Phi)\ ;
`\Delta_A{[}t{]}(\Psi)`\Delta_A{[}t{]}(\Phi)`]
\place(350,250)[{\scriptstyle (1)}]
\square(700,0)/^(->``>`^(->/<900,500>[\ \Delta_A{[}F{]}(\Phi)\ %
`\ F(\Phi + {\bf A} (A,-))\ %
`\ \Delta_A{[}G{]}(\Phi)\ %
`\ G(\Phi + {\bf A} (A,-))\ ;
``t(\Phi + {\bf A} (A,-))`]

\place(1150,250)[{\scriptstyle (2)}]

\efig
\]

\vspace{2mm}

\[
\bfig
\square/^(->`>`>`^(->/<900,500>[\ \Delta_A{[}F{]}(\Psi)\ 
`\ F(\Psi + {\bf A}(A,-))\ 
`\Delta_A{[}G{]}(\Psi)\ 
`\ G(\Psi + {\bf A}(A,-))\ ;
`\Delta_A{[}t{]}(\Psi)`t(\Psi + {\bf A} (A,-))`]

\place(450,250)[{\scriptstyle (3)}]

\square(900,0)/^(->``>`^(->/<1100,500>[\ F(\Psi + {\bf A} (A,-))\ 
`\ F(\Phi + {\bf A} (A,-))\ 
`\ G(\Psi + {\bf A} (A,-))\ 
`\ G(\Phi + {\bf A} (A,-))\ \rlap{\,.};
``t(\Phi + {\bf A} (A,-))`]

\place(1480,250)[{\scriptstyle (4)}]

\efig
\]
The pasted rectangles are equal, and (2), (3) and (4) are
pullbacks, so (1) is too.
\end{proof}

\begin{xcorollary}
\label{Cor-PartTense}

\( \Delta_A [F] \) is a complemented subobject of the
shifted \( F \)
\[
\Delta_A [F] \ \to/^(->/ F \circ S_A
\]
\[
F + \Delta_A [F] \to^\cong F \circ S_A
\]
where the first component is \( F \) of the unit
\( \eta_A \colon \id \to S_A \).

\end{xcorollary}

\subsection{Limit and colimit rules}
\label{SSec-LimColimRules}

\( \Delta_A \) satisfies all the same commutation
properties with respect to limits and colimits as the
\( \Delta \) of \citet{Par24}. This may be proved directly
with virtually the same proofs as in {\em loc.~cit.}
However, just as the usual properties of partial
derivatives follow from their single variable versions
by fixing all the variables but one, those of \( \Delta_A \)
follow from their \( \Delta \) counterparts.

\begin{xproposition}
\label{Prop-Aff}

Objects \( A \) in \( {\bf A} \) and \( \Phi \) in \( {\bf Set}^{\bf A} \)
give an affine functor \( {\bf Set} \to {\bf Set}^{\bf A} \)
\[
A\!f\!f_{A, \Phi} (X) = {\bf A} (A,-) \cdot X + \Phi \rlap{\,.} 
\]
For any tense functor \( F \colon {\bf Set}^{\bf A} \to {\bf Set}^{\bf B} \)
and object \( B \) in \( {\bf B} \), the translated functor
\[
F^B_{A, \Phi} = ({\bf Set} \to^{A\!f\!f_{A, \Phi}} {\bf Set}^{\bf A}
\to^F {\bf Set}^{\bf B} \to^{e\!v_B} {\bf Set})
\]
is taut and
\[
\Delta_A [F] (\Phi) (B) \cong \Delta [F^B_{A, \Phi}](0)\rlap{\,.}
\]

\end{xproposition}

\begin{proof}

The evaluation functors are tense as is \( A\!f\!f_{A,\Phi} \)
so the composite \( e\!v_B \circ F \circ A\!f\!f_{A,\Phi} \) is
too, so taut.
\[
\begin{array}{lll}
\Delta [F^B_{A, \Phi}]  (0) & = & F^B_{A,\Phi} (1) \setminus F^B_{A,\Phi} (0)\\
       & = & F ({\bf A}(A,-) \cdot 1 + \Phi) (B) \setminus F({\bf A} (A,-) \cdot 0 + \Phi)(B)\\
       & \cong & F(\Phi + {\bf A}(A,-)) (B) \setminus F (\Phi) (B)\\
      & = & \Delta_A [F] (B) \rlap{\,.}
\end{array}
\]
\end{proof}

Precomposing by any functor, in particular \( A\!f\!f_{A, \Phi} \), preserves
all limits and colimits (of the \( F \)'s), and precomposing by a
functor that preserves complemented subobjects preserves tense
transformations. The same holds for postcomposing by \( ev_B \).
Furthermore, the \( ev_B \) jointly create limits and colimits. These
considerations give the following results.

\begin{xtheorem}
\label{Thm-LimColimPart}
\begin{itemize}

	\item[(1)] If \( {\bf I} \) is confluent and \( \Gamma \colon {\bf I}
	\to {\cal T}\!\!{\it ense} ({\bf Set}^{\bf A}, {\bf Set}^{\bf B}) \)
	a diagram of tense functors (and tense transformations),
	then
	\[
	\Delta_A [\limr_I \Gamma (I)] \cong \limr_I \Delta_A [\Gamma (I)]
	\rlap{\,.}
	\]
	
	\item[(2)] If \( {\bf I} \) is non-empty and connected and
	\( \Gamma \colon {\bf I} \to {\cal T}\!\!{\it ense}
	({\bf Set}^{\bf A}, {\bf Set}^{\bf B}) \), then
	\[
	\Delta_A [\liml_I \Gamma(I)] \cong \liml_I \Delta_A [\Gamma (I)] \rlap{\,.}
	\]
	
	\item[(3)] For any set \( I \) and tense functors \( F_i \)
	(\( i \in I \)) we have
	\[
	\Delta_A \Bigl[\prod_{i \in I} F_i\Bigr] \cong \sum_{J \subsetneqq I}
	\Bigl(\prod_{j \in J} F_j\Bigr) \times \Bigl(\prod_{k \notin J}
	\Delta_A [F_k]\Bigr) \rlap{\,.}
	\]

\end{itemize}
\end{xtheorem}

\begin{xcorollary}
\label{Cor-DeltaSumProd}

\begin{itemize}

	\item[(1)] \( \Delta_A [F + G] \cong \Delta_A [F] + \Delta_A [G] \)
	
	\item[(2)] \( \Delta_A [C \cdot F] \cong C \Delta_A [F] \)
	(\( C \) a constant set)
	
	\item[(3)] \( \Delta_A [F \times G] \cong (\Delta_A [F] \times G) +
	(F \times \Delta_A [G]) + (\Delta_A[F] \times \Delta_A [G])  \)\rlap{\,.}

\end{itemize}
\end{xcorollary}

We now look at a few special cases.

\begin{xproposition}
\label{Prop-ScalarsPart}

A profunctor \( P \colon {\bf A} \longtod {\bf B} \) gives
a tense \( P \otimes (\ \ ) \colon {\bf Set}^{\bf A} \to {\bf Set}^{\bf B} \)
and \( \Delta_A [P \otimes (\ \ )] \cong P(A,-) \).

\end{xproposition}

\begin{proof}

\( P \otimes (\ \ ) \) is cocontinuous so preserves binary coproducts
\[
\begin{array}{lll}
P \otimes (\Phi + {\bf A} (A,-)) & \cong & P \otimes \Phi + P \otimes {\bf A} (A,-)\\
             & \cong & P \otimes \Phi + P (A,-) \rlap{\,.}
\end{array}
\]
\end{proof}

\begin{xcorollary}
\label{Cor-ScalarsPart}

\( \Delta_A [\id_{{\bf Set}^{\bf A}}] = {\bf A} (A,-) \).

\end{xcorollary}

All that was used in \ref{Prop-ScalarsPart} was that 
\( P \otimes (\ \ ) \) preserved binary coproducts, so we can
improve it.

\begin{xproposition}
\label{Prop-GPart}

If \( F \colon {\bf Set}^{\bf A} \to {\bf Set}^{\bf B} \) preserves
binary coproducts, then
\[
\Delta_A [F] (\Phi) = F ({\bf A} (A,-)) \rlap{\,.}
\]

\end{xproposition}

Note that \( \Delta_A [F] \) is independent of \( \Phi \),
so \( \Delta_A [F] \) is the constant functor
\( {\bf Set}^{\bf A} \to {\bf Set}^{\bf B} \) with value
\( F ({\bf A} (A,-)) \).

We can do better than (2) in the corollary~\ref{Cor-DeltaSumProd}.

\begin{xproposition}
\label{Prop-DeltaProfTimes}

Let \( F \colon {\bf Set}^{\bf A} \to {\bf Set}^{\bf B} \) be
tense and \( P \colon {\bf B} \longtod {\bf C} \) a
profunctor. Then
\[
\Delta_A [P \otimes F] \cong P \otimes \Delta_A [F] \rlap{\,.}
\]

\end{xproposition}

\begin{proof}

We have a coproduct diagram preserved by \( P \otimes (\ \ ) \)
\[
\bfig
\Dtriangle/`^{ (}->`<-_{)}/<500,350>[F(\Phi)\quad `\quad\quad F(\Phi + {\bf A}(A,-1))\ `\Delta_A{[}F{]}(\Phi)\quad\quad ;``]

\place(1300,350)[\longmapsto]

\Dtriangle(2100,0)/`^(->`<-_)/<500,350>[P\otimes F(\Phi)\quad\quad\quad `\quad\quad\quad P\otimes F(\Phi + {\bf A}(A,-1))\ 
`P\otimes (\Delta_A{[}F{]}(\Phi)\quad\quad\quad ;``]

\efig
\]
from which the result follows. \end{proof}

The notation \( P \otimes F \) may need some explanation as 
it doesn't type check. It is componentwise tensor,
\( (P \otimes F) (\Phi) = P \otimes_{\bf B} F (\Phi) \).
We can interpret \ref{Prop-DeltaProfTimes} as saying that
multiplying \( F \) by a matrix of constants is preserved by
differences. But we can generalize this result to the
following, although the interpretation of ``pulling constants out''
may be lost.

\begin{xproposition}
\label{Prop-DeltaGTimes}

If \( F \colon {\bf Set}^{\bf A} \to {\bf Set}^{\bf B} \) is tense
and \( G \colon {\bf Set}^{\bf B} \to {\bf Set}^{\bf C} \)
preserves binary coproducts, then
\[
\Delta_A [GF] = G \Delta_A [F] \rlap{\,.}
\]

\end{xproposition}

\subsection{Analytic functors}

In this section we prove that the generalized analytic
functors of \citet{FioGamHylWin08} are closed under
taking differences and, in fact,
derive an explicit formula for the symmetric sequences
so obtained.

We start with an addition formula for analytic functors
which may look obvious but is frustratingly hard to
make precise. The integral notation for coends
conveniently hides the functoriality of the arguments,
which in the case at hand is not trivial, involving
permutations as it does.

We introduce some notation, without which we run the
risk of drowning in a sea of subscripts, subsubscripts,
ellipses, and so on.

In what follows \( \vec A \) represents an arbitrary object
of \( ! {\bf A} \), \( \langle A_1, \dots , A_n \rangle \) of length
\( n \). Recall that a morphism \( (\sigma, \langle f_1, \dots ,
f_n \rangle) \colon \langle A_1, \dots , A_n \rangle \to
\langle A'_1, \dots , A'_n \rangle \) is a permutation
\( \sigma \in S_n \) and a sequence of morphisms
\[
f_i \colon A_{\sigma i} \to A'_i\rlap{\,.}
\]
We will denote that by \( (\sigma, \vec f\,) \colon \vec A \to \vec A' \).
We also use objects \( \vec X = \langle X_1, \dots , X_k \rangle \)
and \( \vec Y = \langle Y_1, \dots , Y_l \rangle \) whose
lengths are \( k \) and \( l \) respectively. By construction,
\( ! {\bf A} \) is a monoidal category whose tensor is
concatenation
\[
\vec X \otimes \vec Y = \langle X_1, \dots , X_k \, , \,Y_1, \dots , Y_l \rangle
\]
a notation which we use extensively. Of course, it also applies
to morphisms
\[
(\tau, \vec g\,) \otimes (g, \vec h\,) = (\tau + \rho, \vec g \otimes \vec h\,)
\]
where \( \tau + \rho \colon k + l \to k + l \) is the ordinal sum, and
\( \vec g \otimes \vec h \) is concatenation.

We also use the notation, and obvious variants,
\[
\prod \Phi \vec A \colon = \Phi A_1 \times \dots \times 
\Phi A_n
\]
for \( \Phi \) in \( {\bf Set}^{\bf A} \). An element
\( \langle a_1, \dots , a_n \rangle \) of \( \prod \Phi \vec A \)
is denoted \( \vec a \in \prod \Phi \vec A \).

The addition formula alluded to above is given in the
following statement.

\begin{xtheorem}
\label{Thm-AnalyticAddition}

Let \( P \colon {\bf A} \longtod {\bf B} \) be an
\( {\bf A} \)-\( {\bf B} \) symmetric sequence and
\( \widetilde{P} \colon {\bf Set}^{\bf A} \to {\bf Set}^{\bf B} \)
the analytic functor it defines. Then for \( \Phi_1 \) and
\( \Phi_2 \) in \( {\bf Set}^{\bf A} \) and \( B \) in \( {\bf B} \)
we have a natural isomorphism
\[
\widetilde{P} (\Phi_1 + \Phi_2) (B) \cong \int^{\vec X}
\int^{\vec Y} P (\vec X \otimes \vec Y ; B) \times
\prod \Phi_1 \vec X \times \prod \Phi_2 \vec Y \rlap{\,.}
\]

\end{xtheorem}

The idea of the proof is simple:
\begin{eqnarray*}
\widetilde{P} (\Phi_1 + \Phi_2) (B) & = & \int^{\vec A} P(\vec A; B) \times
\prod (\Phi_1 \vec A + \Phi_2 \vec A\,)\\
        & \cong & \int^{\vec A} P(\vec A ; B) \times \sum_{\alpha \colon n \to<70> 2}
              \prod \Phi_\alpha \vec A\\
       & \cong & \int^{\vec X, \vec Y} P (\vec X \otimes \vec Y ; B)
       \times \prod \Phi_1 \vec X \times \prod \Phi_2 \vec Y\\
       & \cong & \int^{\vec X} \int^{\vec Y} P (\vec X \otimes \vec Y; B)
       \times \prod \Phi_1 \vec X \times \prod \Phi_2 \vec Y \rlap{\,.}
       \end{eqnarray*}

The first line is just the definition of \( \widetilde{P} \), the
second line is distributivity of \( \prod \) over \( + \), and the
last line is Fubini for coends. It's in going from the second to
the third line that everything happens. The ``reason'' for the
isomorphism is that for each summand with \( \Phi_1 \) and
\( \Phi_2 \) interspersed ``at random'' in the product, there is
an isomorphism in \( ! {\bf A} \) which permutes them so that
all the \( \Phi_1 \) come first followed by the \( \Phi_2 \). And,
indeed that's the reason. The devil is in the details, as they
say.

We step back and consider how we might show that two
coends are isomorphic. Let \( \Gamma \colon {\bf I}^{op}\times {\bf I} \to
{\bf Set} \) be a functor which we might think of as a
profunctor \( \Gamma \colon {\bf I} \longtod {\bf I} \). The
coend \( \int^I \Gamma (I, I) \) consists of equivalence
classes of elements of \( \Gamma \), \( [I \todd{x}{} I] \),
the equivalence relation generated by identifying
\( x \colon I \longtod I \) with \( x' \colon I' \longtod I' \)
when there are \( f \colon I \to I' \) and \( \ov{x} \colon
I' \longtod I \) such that \( x = \ov{x} f \) and \( x' = f \ov{x} \):
\[
\bfig
\square/>`@{>}|{\bb}`@{>}|{\bb}`>/[I`I'`I`I'\rlap{\,.};f`x`x'`f]

\morphism(500,500)|l|/@{>}|{\bb}/<-500,-500>[I'`I;\ov{x}]

\efig
\]
So \( x \) is equivalent to \( x' \) if there's a zigzag of such
diagrams joining them.

But in the case at hand the equivalence relation is simpler
because both of the diagrams whose coends we're
considering are separable into a product of a contravariant
functor times a covariant one.

\begin{xdefinition}
\label{Def-Separable}

A diagram \( \Gamma \colon {\bf I}^{op} \times {\bf I} \to
{\bf Set} \) is {\em separable} if for every \( f \colon I \to I' \),
\[
\bfig
\square<650,500>[\Gamma(I', I)`\Gamma(I, I)`\Gamma(I', I')`\Gamma(I, I');
\Gamma(f, I)`\Gamma (I', f)`\Gamma(I, f)`\Gamma(f, I')]

\efig
\]
is a pullback.

\end{xdefinition}

For example, if \( \Gamma (I, I') = \Gamma_0 I \times \Gamma_1 I' \)
for \( \Gamma_0 \colon {\bf I}^{op} \to {\bf Set} \) and
\( \Gamma_1 \colon {\bf I} \to {\bf Set} \), then \( \Gamma \) is
separable. Or, if \( {\bf I} \) is a groupoid, every \( \Gamma \) is
separable.

The point of this definition is that the equivalence relation is
generated by identifying \( x \) with \( x' \) when there is an
\( f \colon I \to I' \) such that \( x' f = f x \):
\[
\bfig
\square/>`@{>}|{\bb}`@{>}|{\bb}`>/[I`I'`I`I'\rlap{\,.};f`x`x'`f]

\efig
\]
The \( \ov{x} \) is automatic. This is important because we
can compose such squares.

Let us call an \( x \in \Gamma (I, I) \) a \( \Gamma \)-{\em algebra}
and an \( f \) as above a {\em homomorphism}. Then we get a
category \( {\bf Alg} (\Gamma) \) and \( \int^I \Gamma (I, I) =
\pi_0 {\bf Alg} (\Gamma) \), the set of connected components
of \( {\bf Alg} (\Gamma) \).

Let \( \Theta \colon {\bf J}^{op} \times {\bf J} \to {\bf Set} \)
be another bivariant diagram. A morphism \( (\Xi, \xi) \colon
\Gamma \to \Theta \) is a functor \( \Xi \colon {\bf I} \to {\bf J} \)
and a natural transformation \( \xi \colon \Gamma \to
\Theta (\Xi (-), \Xi (-)) \)
\[
\bfig
\Vtriangle[{\bf I}^{op} \times {\bf I}`{\bf J}^{op} \times {\bf J}`{\bf Set};
\Xi^{op} \times \Xi`\Gamma`\Phi]

\morphism(450,270)/=>/<150,0>[`;\xi]

\morphism(0,750)/>/<1000,0>[{\bf I}`{\bf J};\Xi]

\square(1800,0)/>`@{>}|{\bb}`@{>}|{\bb}`>/[{\bf I}`{\bf J}`{\bf I}`{\bf J}\rlap{\,.};
\Xi`\Gamma`\Theta`\Xi]

\morphism(2000,250)/=>/<150,0>[`;\xi]

\efig
\]
Such a morphism induces a functor
\[
{\bf Alg} (\Xi, \xi) \colon {\bf Alg} (\Gamma) \to {\bf Alg} (\Theta)
\]
\[
\bfig
\morphism(0,500)|l|/@{>}|{\bb}/<0,-500>[I`I;x]

\morphism(600,500)|r|/@{>}|{\bb}/<0,-500>[\Xi (I)`\Xi (I)\rlap{\,.};\xi (I, I) (x)]

\place(300,250)[\longmapsto]

\efig
\]
We are now ready to apply this to our addition formula. Let
\( \Gamma \colon ! {\bf A} \times ! {\bf A} \longtod ! {\bf A} \times ! {\bf A} \)
be given by
\[
\Gamma (\vec X, \vec Y ; \vec X', \vec Y'\,) =
P (\vec X \otimes \vec Y ; B) \times \prod \Phi_1 \vec X'
\times \prod \Phi_2 \vec Y'
\]
and \( \Theta \colon ! {\bf A} \longtod ! {\bf A} \) by
\[
\Theta (\vec A; \vec A'\,) = P (\vec A ; B) \times
\sum_{\alpha \colon n \to<70> 2} \prod \Phi_\alpha \vec A'
\]
with the obvious action on morphisms. Note that \( \Gamma \)
and \( \Theta \) are both products of a covariant part
(with the primes) and a contravariant part (without primes)
so that they are separable. Thus we will be able to compute
the coends by taking connected components of their
categories of elements.

\begin{xtheorem}
\label{Thm-EquivAlg}

With the above notation, there is a morphism
\[
\bfig
\square/>`@{>}|{\bb}`@{>}|{\bb}`>/<600,500>[! {\bf A} \times ! {\bf A}
`! {\bf A}`!{\bf A} \times ! {\bf A}`! {\bf A};
\otimes`\Gamma`\Theta`\otimes]

\morphism(250,250)/=>/<150,0>[`;\xi]

\efig
\]
such that the induced functor
\[
{\bf Alg} (\otimes, \xi) \colon {\bf Alg} (\Gamma) \to {\bf Alg} (\Theta)
\]
is an equivalence of categories.

\end{xtheorem}

\begin{proof}

Throughout, \( B \) is a fixed object of \( {\bf B} \).

An element of \( \Gamma (\vec X, \vec Y ; \vec X', \vec Y'\,) \)
is a triple
\[
\large(p \in P (\vec X \otimes \vec Y ; B),\  \vec x \in \prod \Phi_1 \vec X',\ 
\vec y \in \prod \Phi_2 \vec Y '\large)\rlap{\,,}
\]
and an element of \( \Theta (\vec A, \vec A'\,) \) is a triple
\[
(p \in P (\vec A ; B), \ \alpha \colon n \to<70> 2,
\ \vec a \in \prod \Phi_\alpha \vec A') \rlap{\,,}
\]
where \( \prod \Phi_\alpha \vec A' \) is
\(  \prod^{n'}_{i = 1} \Phi_{\alpha i} A'_i \),
as expected.
\[
\xi \colon \Gamma (\vec X, \vec Y ; \vec X', \vec Y') \to
\Theta (\vec X \otimes \vec Y, \vec X' \otimes \vec Y')
\]
is given by
\[
\xi (p, \vec x, \vec y) = (p \in P (\vec X \otimes \vec Y ; B),\ 
\alpha_{k', l'} \colon k' + l' \to<70> 2,\ 
\langle \vec x, \vec y \rangle \in \prod \Phi_{\alpha_{k', l'}}
(\vec X' \otimes \vec Y')) \rlap{\,.}
\]
Here \( \alpha_{k', l'} \) is the indexing that consists of
\( 1 \)'s followed by \( 2 \)'s,
\[
\alpha_{k', l'} (i) = \left\{ \begin{array}{ll}
1 & \mbox{if\quad} l \leq i \leq k'\\
2 & \mbox{if\quad} k' < i \leq k' + l' \rlap{\,,}
\end{array}
\right.
\]
and \( \langle \vec x, \vec y \rangle \) is concatenation
\[
\langle \vec x, \vec y \rangle = \langle x_1, \dots , x_{k'},
y_1, \dots , y_{l'} \rangle \in \Phi_1 X'_1 \times \dots
\times \Phi_1 X'_{k'} \times \Phi_2 Y'_1 \times
\dots \times \Phi_2 Y'_{l'} \rlap{\,.}
\]
Naturality of \( \xi \) is a straightforward calculation.

The morphism \( ( \otimes, \xi) \) induces a functor
\[
\Xi \colon {\bf Alg} (\Gamma) \to {\bf Alg} (\Theta) \rlap{\,.}
\]
Explicitly, a \( \Gamma \)-algebra is a \( 5 \)-tuple
\[
(\vec X,\  \vec Y,\  p \in P (\vec X \otimes \vec Y ; B),\ 
\vec x \in \prod \Phi_1 \vec X ,\  \vec y \in \prod \Phi_2 \vec Y)
\]
and a \( \Theta \)-algebra is a quadruple
\[
(\vec A, \ p \in P (\vec A ; B), \ \alpha \colon n \to<70> 2,\ 
\vec a \in \prod \Phi_\alpha \vec A\,) \rlap{\,.}
\]
\( \Xi \) assigns to \( (\vec X, \vec Y, p, \vec x , \vec y\,) \) the
algebra \( (\vec X \otimes \vec Y, \ p, \ \alpha_{k, l},\ 
\langle \vec x, \vec y \rangle) \).

A homomorphism \( (\vec X, \vec Y, p, \vec x, \vec y\,) \to
(\vec X', \vec Y', p', \vec x', \vec y\,') \) is a pair of
morphisms in \( ! {\bf A} \),
\[
(\tau, \vec g) \colon \vec X \to \vec X' \mbox{\quad and \quad}
(\rho, \vec h\,) \colon \vec Y \to \vec Y'
\]
preserving everything. It is sent to \( (\tau, \vec g) \otimes
(\rho, \vec h\,) \) by \( \Xi \).

\( \otimes \) is faithful as it is just concatenation, so \( \Xi \)
is also faithful.

If \( (\sigma, \vec f\,) \) is a homomorphism
\[
(\vec X \otimes \vec Y, \ p, \ \alpha_{k, l},\  \langle \vec x, \vec y\,\rangle)
\to ( \vec X' \otimes \vec Y' , \ p', \ \alpha_{k', l'},\ 
\langle \vec x', \vec y\,'\rangle)
\]
we have
\[
\bfig
\Dtriangle/<-`>`<-/<450,300>[k + l`2`k' + l';\sigma`\alpha_{k, l}`\alpha_{k', l'}]

\efig
\]
which implies that \( \sigma \) restricts to bijections
\( \tau \colon k' \to k \) and \( \rho \colon l' \to l \) (by taking
inverse images of \( \{1\} \) and \( \{2\} \)) so \( k' = k \)
and \( l' = l \) and \( \sigma = \tau + \rho \). It follows that
\( \vec f \) consists of morphisms \( (\tau, \vec g) \colon
\vec X \to \vec X' \) and \( (\rho, \vec h) \colon \vec Y
\to \vec Y' \) and the preservation of \( \langle \vec x,
\vec y\,\rangle \) becomes preservation of \( \vec x \)
and \( \vec y\, \) separately. I.e.~\( ( \sigma, \vec f\,) \)
is \( \Xi ((\tau, \vec g), (\rho, \vec h)) \) and so
\( \Xi \) is full.

For any \( \Theta \)-algebra \( (\vec A, p, \alpha, \vec a) \),
there is a permutation of \( \sigma \in S_n \) such that
\[
n \to^\sigma n \to^\alpha 2
\]
is order-preserving, i.e.~all the \( 1 \)'s come first
and then the \( 2 \)'s, so that \( \alpha \sigma = \alpha_{k, l} \)
where \( k \) is the cardinality of \( \alpha^{-1} \{1\} \)
and \( l \) that of \( \alpha^{-1}\{2\} \). Associated to \( \sigma \)
is an isomorphism
\[
(\sigma, \vec 1\,) \colon \vec A \to \vec {A_\sigma}
\]
where \( \vec {A_\sigma} \) is \( \langle A_{\sigma 1},
\dots , A_{\sigma n} \rangle \) and \( \vec 1 =
\langle 1_{A_{\sigma 1}}, \dots , 1_{A_{\sigma n}} \rangle \).
We can transport the \( \Theta \)-algebra structure on
\( \vec A \) to one on \( \vec {A_\sigma} \) giving an
algebra isomorphism
\[
(\sigma, \vec 1\,) \colon (\vec A, p, \alpha, \vec a) \to
(\vec {A_\sigma} , \ p \cdot (\sigma^{-1}, \vec 1\,),\ 
\alpha_{k, l} ,\  \vec {a}_\sigma)
\]
where \( p \cdot (\sigma^{-1}, \vec 1\,) =
P ((\sigma^{-1}, \vec 1\,);B)(p) \) and
\( \vec {a}_\sigma = \langle a_{\sigma 1}, \dots ,
a_{\sigma n} \rangle \) in \( \prod \Phi_{\alpha \sigma}
\vec {A_\sigma} \). The \( \Theta \)-algebras
with indexing of the form \( \alpha_{k, l} \) are
precisely those in the image of \( \Xi \). Indeed,
the \( \vec X \) are the first \( k \) \( A \)'s,
\( \langle A_{\sigma 1}, \dots , A_{\sigma k} \rangle \)
in this case and \( \vec Y \) the last \( l \) of them
\( \langle A_{\sigma(k + l)}, \dots , A_{\sigma (n)} \rangle \).
Similarly \( \vec x = \langle a_{\sigma 1} , \dots ,
a_{\sigma k}\rangle \) and \( \vec y = 
\langle a_{\sigma (k+1)}, \dots a_{\sigma n} \rangle \).
Then \( \Xi (\vec X, \vec Y, p \cdot (\sigma^{-1}, \vec 1\,),
\vec x, \vec y\,) \) is \( (\vec {A_\sigma}, p \cdot (\sigma^{-1}, \vec 1\,), 
\alpha_{k, l}, \vec a) \), so \( \Xi \) is essentially
surjective, which shows it's an equivalence.
\end{proof}

If we take connected components we get
\[
\pi_0 {\bf Alg}(\Gamma) \cong \pi_0 {\bf Alg} (\Theta)
\]
so the coend of \( \Gamma \) is isomorphic to that of
\( \Theta \).

\begin{xcorollary}
\label{Cor-EquivAlg}

\[
\int^{\vec X, \vec Y} P (\vec X \otimes \vec Y ; B) \times
\prod \Phi_1 \vec X \times \prod \Phi_2 \vec Y
\cong \int^{\vec A} P (\vec A ; B) \times
\sum_{\alpha \colon n \to<70> 2} \prod \Phi_\alpha \vec A
\rlap{\,.}
\]

\end{xcorollary}

Our addition formula, Theorem~\ref{Thm-AnalyticAddition},
now follows by a simple application of the Fubini theorem
for coends, which is what we wanted, but
Theorem~\ref{Thm-EquivAlg} is a stronger result. 

Our next step in the derivation of the formula for
\( \Delta_A [\widetilde{P}] \) is to specialize our addition
formula to the case \( \Phi_1 = \Phi \) and
\( \Phi_2 = {\bf A} (A, -) \). This gives
\[
\widetilde{P} (\Phi + {\bf A} (A, -)) (B) =
\int^{\vec X} \int^{\vec Y} P (\vec X \otimes \vec Y ; B)
\times \prod \Phi \vec X \times \prod {\bf A}
(A, \vec Y)
\]
in which the expression
\[
\prod {\bf A} (A, \vec Y) =
{\bf A} (A, Y_1) \times \dots 
\times {\bf A}(A, Y_l)
\]
appears, not surprisingly, as it already appears
in the definition of \( \widetilde{P} \). It defines a
functor
\[
\prod {\bf A}(A, -) \colon !{\bf A} \to {\bf Set}
\]
closely related to the representable functor
\( ! {\bf A} (A^{\otimes n}, -) \) where \( A^{\otimes n} =
\langle A, \dots , A \rangle \), the \( n \)-fold
tensor of \( A \).

\begin{xproposition}
\label{Prop-ProdSum}

With the above notation we have
\[
\prod {\bf A} (A, -) \cong \sum^\infty_{n = 0}
! {\bf A} (A^{\otimes n}, -)/S_n \rlap{\,.}
\]

\end{xproposition}

\begin{proof}
If \( \vec Y = \langle Y_1, \dots , Y_l \rangle \), then
\( ! {\bf A} (A^{\otimes n} , \vec Y) \) is \( 0 \) unless
\( l = n \) in which case an element of \( ! {\bf A}
(A^{\otimes n}, \vec Y) \) is a morphism
\[
(\sigma, \vec f\,) \colon A^{\otimes n} \to \vec Y
\]
so that \( ! {\bf A}(A^{\otimes n}, \vec Y) \cong
S_n \times {\bf A} (A_1 Y_1) \times \dots
\times {\bf A} (A, Y_n)\) and if we quotient by
\( S_n \) we get
\[
! {\bf A} (A^{\otimes n}, \vec Y) / {S_n} \cong
\prod {\bf A} (A, \vec Y)
\]
easily seen to be natural in \( \vec Y \). The result
follows.
\end{proof}

\begin{xlemma}
\label{Lem-TensorWithPi}

Let \( W \colon !{\bf A}^{op} \to {\bf Set} \). Then
\[
\int^{\vec Y} W(\vec Y) \times \prod {\bf A} (A, \vec Y)
\cong
\sum^\infty_{n = 0} W (A^{\otimes n})/S_n \rlap{\,.}
\]

\end{xlemma}

\begin{proof}
\begin{eqnarray*}
\int^{\vec Y} W({\vec Y}) \times \prod {\bf A} (A, \vec Y) & \cong &
       \int^{\vec Y} W (\vec Y) \times \sum^\infty_{n = 0} ! {\bf A}(A^{\otimes n},
       \vec Y)/S_n\\
 & \cong &\sum^\infty_{n = 0} \int^{\vec Y} W (\vec Y) \times
 (! {\bf A} (A^{\otimes n}, \vec Y)/{S_n}\\
 & \cong & \sum^\infty_{n =0} \Bigl(\int^{\vec Y} W (\vec Y) \times
 ! {\bf A} (A^{\otimes n} , \vec Y) \Bigr)/S_n\\
 & \cong & \sum^\infty_{n = 0} W (A^{\otimes n})/S_n \rlap{\,.}
 \end{eqnarray*}
 The second isomorphism is commutation of coends and
 coproducts, the third commutation of coends with colimits
 (``modding out'' by \( S_n \) {\em is} a colimit), the last
 isomorphism comes from the fact that tensoring with a representable
 is substitution. 
 \end{proof}

 \begin{xcorollary}
 \label{Cor-PTildeWithHom}
 
 \[
 \widetilde{P} (\Phi + {\bf A} (A, -))(B) \cong
 \int^{\vec X} \sum^\infty_{n = 0} P (\vec X \otimes A^{\otimes n} ; B)/
 (\{\id_k\} \times S_n) \times
 \prod \Phi \vec X
 \]
 
 \end{xcorollary}
 
 \begin{proof}
 
 \[
 \widetilde{P} (\Phi + {\bf A} (A, -)) (B) \cong \int^{\vec X}
 \int^{\vec Y} P (\vec X \otimes \vec Y ; B) \times
 \prod \Phi \vec X \times \prod {\bf A} (A, \vec Y) \rlap{\,.}
 \]
 If we fix \( \vec X \) and consider the coend over \( \vec Y \),
 we can apply the previous lemma with
 \[
W (\vec Y) = P (\vec X \otimes \vec Y ; B) \times \prod \Phi \vec X
\]
and the result follows immediately.
\end{proof}

\begin{xcorollary}
\label{Cor-DeltaPTilde}

\[
\Delta_A [\widetilde{P}] (\Phi) (B) \cong \int^{\vec X}
\sum^\infty_{n = 1} P(\vec X \otimes A^{\otimes n} ; B)/
(\{\id_k\} \times S_n) \times
\prod \Phi \vec X \rlap{\,.}
\]

\end{xcorollary}

\begin{proof}

The inclusion \( \widetilde{P} (\Phi) \hookrightarrow \widetilde{P}
(\Phi + {\bf A} (A, -)) \) corresponds to the \( n = 0 \) summand.
\end{proof}

For any \( {\bf A} \)-\( {\bf B} \) symmetric sequence
\( P \colon ! {\bf A} \longtod {\bf B} \) and object \( A \) of
\( {\bf A} \) we define a new symmetric sequence
\( \nabla_A P \colon ! {\bf A} \longtod {\bf B} \) by the formula
\[
\nabla_A P (\vec X ; B) = \sum^\infty_{n = 1} 
P (\vec X \otimes A^{\otimes n} ; B) /(\{\id_k\} \times S_n) \rlap{\,.}
\]
Now Corollary~\ref{Cor-DeltaPTilde} can be stated in its final
form, giving the main theorem of the section.

\begin{xtheorem}
\label{Thm-DeltaAnalytic}

Analytic functors \( {\bf Set}^{\bf A} \to {\bf Set}^{\bf B} \)
are closed under taking differences. If \( P \colon !{\bf A}
\longtod {\bf B} \) is a symmetric sequence, then
\[
\Delta_A [\widetilde{P}] \cong \widetilde{\nabla_A P} \rlap{\,.}
\]

\end{xtheorem}

The definition of \( \nabla_A P \) as a coproduct of quotients
is clear but for formal manipulations a more abstract
definition is useful. Let \( {\bf S}_+ \) be the category
whose objects are positive finite cardinals, \( k > 0 \),
and whose morphisms are bijections. So \( {\bf S}_+ \)
is the coproduct
\[
{\bf S}_+ = \sum^\infty_{k = 1} {\bf S}_k
\]
where \( {\bf S}_k \) is the symmetric group \( {\bf S}_k \)
considered as a one-object category.

Given an \( {\bf A} \)-\( {\bf B} \) symmetric sequence
\( P \colon ! {\bf A} \longtod {\bf B} \) and an object
\( A \) of \( {\bf A} \) we get an \( {\bf S}_+ \) family of
\( {\bf A} \)-\( {\bf B} \) symmetric sequences
\[
P_A \colon {\bf S}^{op}_+ \to \Prof (! {\bf A}, {\bf B})
\]
\[
P_A (k) (A_1 \dots A_n ; B) = P(A_1 \dots A_n , A, A, \dots , A ; B)
\]
where there are \( k \) \( A \)'s. Functoriality and naturality
are obvious. Now \( \nabla_A P = \limr_k P_A (k) \).

\begin{xproposition}
\label{Prop-NablaLin}

For \( P \colon ! {\bf A} \longtod {\bf B} \) an \( {\bf A} \)-\( {\bf B} \)
symmetric sequence and \( Q \colon {\bf B} \longtod {\bf C} \) a
profunctor, we have
\[
\nabla_A (Q \otimes P) \cong Q \otimes \nabla_A P \rlap{\,.}
\]

\end{xproposition}

\begin{proof}

\begin{eqnarray*}
\nabla_A (Q \otimes P) (A_1 \dots A_n ; C) & \cong
     & \limr_k \int^B Q (B, C) \times P(A_1 \dots A_n, A \dots A ; B)\\
& \cong & \int^B Q (B, C) \times \limr_k P(A_1 \dots A_n , A \dots A ; B)\\
& \cong & \int^B Q (B, C) \times \nabla_A P (A_1 \dots A_n ; B)\\
& \cong & (Q \otimes \nabla_A P) (A_1 \dots A_n ; C) \rlap{\,.}
\end{eqnarray*}
\end{proof}

\begin{xcorollary}
\label{Cor-NablaLin}

For any \( {\bf A} \)-\( {\bf B} \) symmetric sequence \( P \)
we have
\[
\nabla_A P \cong P \otimes \nabla_A \Id_{!{\bf A}} \rlap{\,.}
\]

\end{xcorollary}

\begin{proof}

\[
\nabla_A P \cong \nabla_A (P \otimes \Id_{!{\bf A}}) \cong
P \otimes \nabla_A \Id_{!{\bf A}} \rlap{\,.}
\]
\end{proof}

\( \nabla_A \Id_{!{\bf A}} \) is easy to describe:
\[
\nabla_A \Id_{!{\bf A}} \colon {\bf Set}^{!{\bf A}} \to
{\bf Set}^{!{\bf A}}
\]
\[
\nabla_A \Id_{!{\bf A}} (A_1 \dots A_n ; A'_1 \dots A'_m) \cong
\sum^\infty_{k = 1} ! {\bf A} (A_1 \dots A_n, A \dots A ;
A'_1 \dots A'_m) / (\{\id_n\} \times S_k)
\]
which is \( 0 \) if \( m \leq n \) and
\[
! {\bf A} (A_1 \dots A_n , A, \dots , A ; A'_1 \dots A'_m)/
(\{\id_n\} \times S_{m-n})
\]
when \( m > n \). There are \( m - n \) \( A \)'s and the
action we're modding out by is \( S_{m - n} \) acting
on those \( A \)'s.

There is also a generic difference formula.

\begin{xcorollary}
\label{Cor-DeltaAnalytic}

\[
\Delta_A [\widetilde{P}] \cong P \otimes \Delta_A
[\Id_{!{\bf A}}] \rlap{\,.}
\]

\end{xcorollary}

\begin{proof}

\[
\begin{array}{lll}
\Delta_A \widetilde{P} & \cong \widetilde{\nabla_A P}
       &(\mbox{Thm.~\ref{Thm-DeltaAnalytic}})\\
 & \cong (P \otimes \nabla_A \Id_{! {\bf A}}){\widetilde{\ \ \ }} &(\mbox{Cor.~\ref{Cor-NablaLin}})\\
 & \cong P \otimes \widetilde{\nabla_A \Id_{! {\bf A}}}
         &(\mbox{Cor.~\ref{Cor-TildeLin}})\\
 & \cong P \otimes \Delta_A {[}\Id_{! {\bf A}}{]}
     & (\mbox{Thm.~\ref{Thm-DeltaAnalytic}})
 \end{array}
 \]

\end{proof}


\subsection{Higher differences}
\label{SSec-HigherDiff}

As \( \Delta_A [F] \) is also tense, its difference can also
be taken \( \Delta_{A'} [\Delta_A [F]] = \Delta_{A', A} [F] \)
and so on, iteratively. For any sequence \( \langle A_1 \dots
A_n \rangle \) of length \( n \) of objects of \( {\bf A} \) we
define
\[
\Delta_{\langle A_i  \rangle} {[}F{]} = \left\{ \begin{array}{ll}
 F  & \mbox{\ if\ }  n = 0\\
 \Delta_{A_1} {[}\Delta_{\langle A_2 \dots A_{n}\rangle} {[}F{]}{]} &  \mbox{\ if\ } n \geq 1 \rlap{\,.}
 \end{array}
 \right. 
\]


\begin{xdefinition}
\label{Def-NewElt}

We say that an element of \( F (\Phi + {\bf A} (A_1, -) +
\cdots + {\bf A} (A_n, -)) (B) \) is {\em new} (for \( \langle A_1, \cdots ,
A_n\rangle) \)) if it is not in any
\( F (\Phi + {\bf A}(A_{\alpha 1}, -) + \cdots +
{\bf A} (A_{\alpha k}, -)) (B) \) for any proper mono
\( \alpha \colon k  \to/ >->/<200>\ n \).

\end{xdefinition}

If an element is in \( F \) of a subsum, it's in
every bigger subsum, so it is sufficient to consider
only those subsums with one less term.Thus the
new elements are those in the set difference
\[
F (\Phi + \sum_{i = 1}^n {\bf A} (A_i, -))(B) \setminus
\bigcup_{j = 1}^n F (\Phi + \sum_{i \neq j}
{\bf A} (A_i, -)) (B) \rlap{\ .}
\]

\begin{xtheorem}
\label{Thm-HighDiff}

The higher difference \( \Delta_{\langle A_i \rangle} {[}F{]} (\Phi) \)
consists of the new elements of
\( F (\Phi + {\bf A} (A_1, -) + \cdots + {\bf A} (A_n, -)) \).

\end{xtheorem}



\begin{proof} We prove this by induction on \( n \). For
\( n = 0, 1 \) the result holds by definition. Assume the
result holds for sequences of length \( n - 1 \) and take
\( \langle A_i \rangle = \langle A_1, \dots, A_n \rangle \).
Let \( \langle A_i \rangle^+ = \langle A_2, \dots, A_n \rangle \).

An element of \( \Delta_{\langle A_i\rangle} [F](\Phi)(B) \) is
an element of \( \Delta_{\langle A_i \rangle^+}
[F] (\Phi + {\bf A} (A_1, -)) (B) \) which is not in
\( \Delta_{\langle A_I \rangle^+} [F] (\Phi) (B) \). An element of
\( \Delta_{\langle A_i \rangle^+}
[F] (\Phi + {\bf A} (A_1, -)) (B) \) is, by the induction hypothesis, an element of
\begin{equation}\tag{1}
F (\Phi + {\bf A} (A_1, -) + \sum^n_{i = 2} {\bf A} (A_i, -)) (B)
\cong F (\Phi + \sum^n_{i = 1} {\bf A} (A_i, -)) (B)
\end{equation}
not in
\begin{equation}\tag{2}
F (\Phi + {\bf A} (A_1, -) + \sum^n_{i = 2,  i \neq j}
{\bf A} (A_i, -)) (B) \cong F (\Phi +
\sum^n_{i = 1,  i \neq j} {\bf A} (A_i, -)) (B)
\end{equation}
for any \( 2 \leq j \leq n \). From this we must exclude
the elements of \( \Delta_{\langle A_i\rangle^+} [F] (\Phi)(B) \)
and these, again by the induction hypothesis, are
elements of
\begin{equation}\tag{3}
F (\Phi + \sum^n_{i = 2} {\bf A} (A_i, -)) (\Phi) (B)
\end{equation}
except for any in some
\begin{equation}\tag{4}
F (\Phi + \sum^n_{i = 2, i \neq j} {\bf A} (A_i, -)) (\Phi)(B)
\end{equation}
for \( 2 \leq j \leq n \).

To summarize,
\[
\Delta_{\langle A_i\rangle} [F] (\Phi) (B) =
((1) \setminus (2)) \setminus ((3) \setminus (4)) \rlap{\,,}
\]
but \( (4) \subseteq (2) \) so
\[
\Delta_{\langle A_i\rangle} [F] (\Phi) (B) =
(1) \setminus ((2) \cup (3)) \rlap{\,.}
\]
Now \( (2) \cup (3) \) is the union of
\[
F (\Phi + \sum_{i = 1, i \neq j} {\bf A} (A_i, -))(B)
\]
over all \( j, \) \(1 \leq j \leq n \), and the result follows.
\end{proof}

We see from this formula that \( \Delta_{\langle A_i \rangle} [F] \)
is independent of the order of the differences, a version of
Clairaut's theorem.

\begin{xcorollary}
\label{Cor-Clairaut}

Let \( \langle A_i \rangle \) be a sequence of length \( n \)
of objects of \( {\bf A} \) and \( \sigma \in S_n \) a permutation,
then
\[
\Delta_{\langle A_{\sigma i} \rangle} [F] \cong
\Delta_{\langle A_i \rangle} [F] \rlap{\,.}
\]

\end{xcorollary}

%
%
%
%
%
%

\section{The discrete Jacobian}
\label{Sec-Jacobian}

\subsection{Definitions and functoriality}
\label{SSec-JacobianDef}

Let \( F \colon {\bf Set}^{\bf A} \to {\bf Set}^{\bf B} \) be a
tense functor and let \( f \colon A \to A' \) be a morphism of
\( {\bf A} \). Then, as
\[
\bfig
\square/^(->`=`>`^(->/<600,500>[\ \Phi\ `\ \Phi + {\bf A}(A',-)\ `\ \Phi\ `\ \Phi + {\bf A}(A,-)\ ;
``\Phi + {\bf A}(f,-)`]

\efig
\]
is a pullback of complemented objects, so is
\[
\bfig
\square/^(->`=`>`^(->/<700,500>[\ F\Phi\ `\ F(\Phi + {\bf A}(A',-))\ `\ F\Phi\ `\ F(\Phi + {\bf A}(A,-))\ ;
``F(\Phi + {\bf A}(f,-))`]

\efig
\]
and it follows that \( F (\Phi + {\bf A} (f,-)) \) restricts to
complements giving another pullback
\[
\bfig
\square/^(->`>`>`^(->/<850,500>[\ \Delta_{A'}{[}F{]}(\Phi)\ 
`\ F(\Phi + {\bf A}(A',-))\ 
`\ \Delta_A{[}F{]}(\Phi)\ \ 
`\ F(\Phi + {\bf A}(A,-))\ \rlap{\,.};
`\Delta_f {[}F{]}(\Phi)`F(\Phi + {\bf A}(f,-))`]

\efig
\]
This proves the following:

\begin{xproposition}
\label{Prop-FuncInA}

For any \( \Phi \) in \( {\bf Set}^{\bf A} \),
\( \Delta_A[F](\Phi) \) is functorial in \( A \), i.e.~is the object part of a functor
\[
\Delta [F] (\Phi) \colon {\bf A}^{op} \to {\bf Set}^{\bf B} \rlap{\,.}
\]

\end{xproposition}

By exponential adjointness we get a functor
\( {\bf A}^{op} \times {\bf B} \to {\bf Set} \), i.e.~a profunctor
\( {\bf A} \longtod {\bf B} \).

\begin{xdefinition}
\label{Def-JacobianProf}

The {\em (discrete) Jacobian profunctor} of \( F \) at \( \Phi \)
\[
\Delta [F](\Phi) \colon {\bf A} \longtod {\bf B}
\]
is given by
\[
\Delta [F](\Phi) (A,B) = \Delta_A [F](\Phi) (B) \rlap{\,.}
\]

\end{xdefinition}

It's more or less clear that \( \Delta [F](\Phi) \) is functorial
in \( \Phi \), which we express in the following proposition.

\begin{xproposition}
\label{Prop-JacobianFunct}

For any tense functor \( F \colon {\bf Set}^{\bf A} \to
{\bf Set}^{\bf B} \), \( \Delta [F] (\Phi) \) is the object
part of a tense functor
\[
\Delta [F] \colon {\bf Set}^{\bf A} \to {\bf Set}^{{\bf A}^{op} \times {\bf B}}
= \Prof ({\bf A}, {\bf B}) \rlap{\,.}
\]

\end{xproposition}

\begin{proof}

For a natural transformation \( t \colon \Phi \to \Psi \) and
object \( A \) in \( {\bf A} \),
\[
\bfig
\square/^(->`>`>`^(->/<600,500>[\ \Phi\ `\ \Phi + {\bf A}(A,-)\ 
`\ \Psi\ `\ \Psi + {\bf A} (A,-)\ ;
`t`t+ {\bf A}(A,-)`]

\efig
\]
is a pullback of a complemented subobject, so
\[
\bfig
\square/^(->`>`>`^(->/<700,500>[\ F\Phi\ `\ F(\Phi + {\bf A}(A,-))\ 
`\ F\Psi\ `\ F(\Psi + {\bf A} (A,-))\ ;
`Ft`F(t+ {\bf A}(A,-))`]

\efig
\]
is too. So \( F(t + {\bf A}(A,-)) \) restricts to the complements,
giving another pullback
\begin{equation}\tag{*}
\bfig
\square/^(->`>`>`^(->/<800,500>[\ \Delta{[}F{]}\Phi\ `\ F(\Phi + {\bf A}(A,-))\ 
`\ \Delta{[}F{]}\Psi\ `\ F(\Psi + {\bf A} (A,-))\ \rlap{\,,};
`\Delta{[}F{]}t`F(t+ {\bf A}(A,-))`]

\efig
\end{equation}
hence functoriality.

We still must prove that it is tense.

Proposition~\ref{Prop-PartTense} says that for a fixed
\(  A \), \( \Delta_A[F] \colon {\bf Set}^{\bf A} \to {\bf Set}^{\bf B} \)
is tense and \( \Delta_A [F] \) is the composite
\[
{\bf Set}^{\bf A} \to^{\Delta[F]} {\bf Set}^{{\bf A}^{op} \times {\bf B}}
\to^{ev_B} {\bf Set}^{\bf B} \rlap{\,.}
\]
The \( ev_B \) are the evaluation functors which preserve
pullbacks and collectively reflect them, so that \( \Delta[F] \)
will preserve pullbacks of complemented subobjects.
However, the \( ev_B \) don't reflect complemented
subobjects, so we still must show that \( \Delta[F] \)
preserves those.

Let \( \Phi_0\  \to/^(->/ \Phi \) be a complemented
subobject. We want to show that \( \Delta [F] (\Phi_0)\ 
\to/>->/ \Delta [F] (\Phi) \) is complemented, or
equivalently, for every \( f \colon A'\to A \) and
\( g \colon B \to B' \)
\[
\bfig
\square/^(->`>`>`^(->/<1000,500>[\ \DelF(\Phi_0)(A,B)\ 
`\ \DelF(\Phi)(A,B)\ 
`\ \DelF(\Phi_0)(A',B')\ 
`\ \DelF(\Phi)(A',B')\ ;
`\DelF (\Phi_0)(f,g)`\DelF(\Phi)(f,g)`]

\efig
\]
is a pullback. We can do this separately for \( f \) and \( g \),
fixing \( B \) and then \( A \). We already know for fixed
\( A \) it's a pullback. So let's fix \( B \).

Let \( f \colon A' \to A \) and consider
\[
\bfig
\node a(0,0)[\ F (\Phi + {\bf A}(A,-))(B)\ ]
\node b(1700,0)[\ F (\Phi + {\bf A}(A',-))(B)\ ]
\node c(0,500)[\ \DelF(\Phi)(A,B)\ ]
\node d(1700,500)[\ \DelF(\Phi)(A,'B)\ ]
\node e(0,1000)[\ \DelF(\Phi_0)(A,B)\ ]
\node f(1700,1000)[\ \DelF(\Phi_0)(A',B)\ ]

\arrow|b|/>/[a`b;F(\Phi + {\bf A}(f,-))(B)]
\arrow|b|/>/[c`d;\DelF(\Phi)(f,B)]
\arrow|a|/>/[e`f;\DelF(\Phi_0)(f,B)]

\arrow/ >->/[c`a;]
\arrow/ >->/[e`c;]

\arrow/ >->/[d`b;]
\arrow/ >->/[f`d;]

\place(850,250)[\Pb]

\efig
\]
and
\[
\bfig
\node a'(0,0)[\ F (\Phi + {\bf A}(A,-))(B)\ ]
\node b'(1700,0)[\ F (\Phi + {\bf A}(A',-))(B)\ \rlap{\,.}]
\node c'(0,500)[\ F(\Phi_0 + {\bf A} (A,-))(B)\ ]
\node d'(1700,500)[\ F(\Phi_0 + {\bf A}(A',-))(B)\ ]
\node e'(0,1000)[\ \DelF(\Phi_0)(A,B)\ ]
\node f'(1700,1000)[\ \DelF(\Phi)(A',B)\ ]

\arrow|b|/>/[a'`b';F(\Phi + {\bf A}(f,-))(B)]
\arrow|b|/>/[c'`d';F(\Phi_0 + {\bf A}(f,-))(B)]
\arrow|a|/>/[e'`f';\DelF(\Phi_0)(f,B)]
\arrow/ >->/[c'`a';]
\arrow/ >->/[e'`c';]

\arrow/ >->/[d'`b';]
\arrow/ >->/[f'`d';]

\place(850,250)[\Pb]
\place(850,750)[\Pb]

\efig
\]
The second and third squares are pullbacks by \( (*) \)
and the fourth because \( F \) is tense. As the composite
of the first and second squares is equal to the composite
of the third and fourth, we get that the first square is a
pullback, which shows that \( \DelF (\Phi_0) \ \ \to/>->/
\DelF(\Phi) \) is complemented.
\end{proof}

To complete the discussion of functoriality of \( \Delta \)
note that \( \Delta_A [F] (\Phi) \) is a subfunctor of
\( F (\Phi + {\bf A} (A,-)) \) which is not only functorial
in \( \Phi \) and \( A \) but by Proposition~\ref{Prop-PartTense}
also in \( F \) but only for tense transformations.
Proposition~\ref{Prop-EvalTense} says that the evaluation
functors \( ev_B \) jointly reflect tenseness of transformations, so that
\( \Delta_A [t] \) itself will be tense. Thus we get a functor
\[
\Delta \colon {\cal T}\!\!{\it ense} ({\bf Set}^{\bf A}, {\bf Set}^{\bf B})
\to {\cal T}\!\!{\it ense} ({\bf Set}^{\bf A}, {\bf Set}^{{\bf A}^{op} \times {\bf B}})
\]
the {\em (discrete) Jacobian functor}.

There are various ways of reformulating the Jacobian which
are of independent interest.

Given a tense functor
\( F \colon {\bf Set}^{\bf A} \to {\bf Set}^{\bf B} \), we get
another tense functor analogous to the differential operator
\[
D[F] \colon {\bf Set}^{\bf A} \times {\bf Set}^{\bf A} \to
{\bf Set}^{\bf B}
\]
\[
D[F] (\Phi, \Psi) = \DelF (\Phi) \otimes_{\bf A} \Psi
\]
where \( \Psi \) is considered as a profunctor
\( {\bf 1} \longtod {\bf A} \).

\begin{xdefinition}
\label{Def-DiffOp}

\( D [F] \) is called the {\em difference operator}.

\end{xdefinition}

In \citep{Par24} we
used the finite projection \( {\bf Set} \times {\bf Set} \to
{\bf Set} \) as a tangent bundle and saw that this
supported a definition of functorial differences where
the lax chain-rule was actually a lax functor. This
generalizes to the multivariable setting. We define
\[
\bfig
\square<900,500>[{\bf Set}^{\bf A} \times {\bf Set}^{\bf A}
`{\bf Set}^{\bf B} \times {\bf Set}^{\bf B}
`{\bf Set}^{\bf A}`{\bf Set}^{\bf B};
T{[}F{]}`P_1`P_1`F]

\efig
\]
by \( T[F] (\Phi, \Psi) = (F \Phi, \DelF(\Phi) \otimes_{\bf A} \Psi) \).
We see that \( T[F] \) preserves colimits in the second variable.

\begin{xdefinition}
\label{Def-TanBun}

\( T [F] \) is called the {\em (discrete) tangent functor}.

\end{xdefinition}

Profunctors \( {\bf A} \longtod
{\bf B} \) are in bijection with profunctors \( {\bf B}^{op} \longtod
{\bf A}^{op} \):
\[
\frac{P \colon {\bf A}^{op} \times {\bf B} \to {\bf Set}}
{P^\top \colon ({\bf B}^{op}) \times {\bf A}^{op} \to {\bf Set}}
\]
i.e.~\( P^\top (B,A) = P(A,B) \), the transpose as matrices.
This gives the reverse difference operator
\[
\Delta^\top [F] \colon {\bf Set}^A \to \Prof ({\bf B}^{op}, {\bf A}^{op}) .
\]

\begin{xdefinition}
\label{Def-RevDiff}

\( \Delta^\top [F] \) is the {\em reverse difference operator}.

\end{xdefinition}

This suggests that we
take as the cotangent bundle the first projection \( {\bf Set}^{\bf A} \times {\bf Set}^{{\bf A}^{op}}
\to {\bf Set}^{\bf A} \). As the Yoneda
embedding \( Y \colon {\bf A}^{op} \ \to/>->/ {\bf Set}^{\bf A} \)
is the cocompletion of \( {\bf A} \),
the category of cocontinuous functors
\[
{\bf Set}^{\bf A} \to {\bf Set}
\]
is equivalent to the category of functors
\[
{\bf A}^{op} \to {\bf Set}
\]
i.e.~\( {\bf Set}^{{\bf A}^{op}} \). So \( {\bf Set}^{{\bf A}^{op}} \)
has a legitimate claim to be the (linear) dual of \( {\bf Set}^{\bf A} \).
Now we can extend the reverse difference to the cotangent
bundle. Given a tense functor \( F \colon {\bf Set}^{\bf A} \to
{\bf Set}^{\bf B} \), we first pull back the cotangent bundle
along \( F \)
\[
\bfig
\square<1000,500>[{\bf Set}^{\bf A} \times {\bf Set}^{{\bf B}^{op}}
`{\bf Set}^{\bf B} \times {\bf Set}^{{\bf B}^{op}}
`{\bf Set}^{\bf A}`{\bf Set}^{\bf B};
\langle F,P_2\rangle`P_1`P_1`F]

\place(500,250)[\Pb]

\efig
\]
and then take the functor \( {\rm coT} [F] \)
\[
(\Phi, \Theta) \to/|->/<250> (\Phi, \Delta^\top [F](\Phi) \otimes \Theta)
\]
\[
\bfig
\Vtriangle[{\bf Set}^{\bf A} \times {\bf Set}^{{\bf B}^{op}}
`{\bf Set}^{\bf A} \times {\bf Set}^{{\bf A}^{op}}
`{\bf Set}^{\bf A}\rlap{\,.};`P_1`P_1]

\efig
\]
In this \( \Theta \) in \( {\bf Set}^{{\bf B}^{op}} \) is considered as
a profunctor \( {\bf 1} \longtod {\bf B}^{op} \).

\begin{xdefinition}
\label{Def-CotBun}

\( {\rm coT} [F] \) is the {\em cotangent functor}.

\end{xdefinition}

A differential form is a global
section of the cotangent bundle, which in our case amounts to
a functor \( {\bf Set}^{\bf A} \to {\bf Set}^{{\bf A}^{op}} \).

For a tense \( F \colon {\bf Set}^{\bf A} \to {\bf Set}^{\bf B} \) we
get another tense functor \( \DelF \colon {\bf Set}^{\bf A} \to
{\bf Set}^{{\bf A}^{op} \times {\bf B}} \) which, upon composing
with the evaluation at \( B \), \( ev_B \colon {\bf Set}^{{\bf A}^{op} \times {\bf B}}
\to {\bf Set}^{{\bf A}^{op}} \) gives another tense functor
\( {\bf Set}^{\bf A} \to {\bf Set}^{{\bf A}^{op}} \). This way
the difference \( \DelF \) may be viewed as a \( {\bf B} \)-family
of differential forms
\[
{\bf B} \to {\cal T}\!\!{\it ense} ({\bf Set}^{\bf A}, {\bf Set}^{{\bf A}^{op}})\rlap{\,.}
\]
It is tempting to write \( \Omega^1 ({\bf Set}^{\bf A}) \) for
\( {\cal T}\!\!{\it ense} ({\bf Set}^{\bf A}, {\bf Set}^{{\bf A}^{op}})\).

\begin{xdefinition}
\label{Def-DiffForm}

\( \Delta [F] (-) (B) \) is the {\em differential form} of \( F \)
at \( B \).

\end{xdefinition}

\subsection{Product and sum rules}

The evaluation functors \( ev_A \colon {\bf Set}^{{\bf A}^{op} \times {\bf B}} \to
{\bf Set}^{\bf B} \) jointly create limits and colimits, and as the
composites \( ev_A \circ \DelF (\Phi) \) are \( \Delta_A [F] (\Phi) \),
the limit rules of Section~\ref{SSec-LimColimRules} lift to
\( \DelF \).

Theorem~\ref{Thm-LimColimPart} gives the following.

\begin{xtheorem}
\label{Thm-LimColimDiff}

\begin{itemize}

	\item[(1)] If  \( \Gamma \colon
	{\bf I} \to \Tense ({\bf Set}^{\bf A}, {\bf Set}^{\bf B}) \) then
	\( \Delta [\limr_I \Gamma I] \cong \limr_I \Delta [\Gamma I] \).
	
	\item[(2)] If \( {\bf I} \) is non-empty and connected and
	\( \Gamma \colon
	{\bf I} \to \Tense ({\bf Set}^{\bf A}, {\bf Set}^{\bf B}) \),
	then \( \Delta [\liml_I \Gamma I] \cong \liml_I \Delta
	[\Gamma I] \).
	
	\item[(3)] For any set \( I \) and tense functors \( F_i \),
	\( i \in I \), we have
	\[
	\Delta \Bigl[\prod_{i \in I} F_i\Bigr] \cong \sum_{J \subsetneqq I}
	\Bigl(\prod_{j \in J} F_j\Bigr) \times \prod_{k \notin J}
	\Delta[F_k] .
	\]

\end{itemize}
\end{xtheorem}

\begin{xcorollary}
\label{Cor-LimColimDiff}

(1) \( \Delta [F + G] \cong \DelF + \Delta [G] \).

\vspace{2mm}

(2) \( \Delta [C \cdot F] \cong C \cdot \DelF \) for any set \( C \).

\vspace{2mm}

(3) \( \Delta [F \times G] \cong (\DelF \times G) + (F \times \Delta [G])
+ (\DelF \times \Delta [G]) \).

\end{xcorollary}

Note that on the right hand side of (3) we have
\( \DelF \times G \) for example. \( \DelF \) is a functor
\( {\bf Set}^{\bf A} \to \Sett \) whereas \( G \) is a functor
\( {\bf Set}^{\bf A} \to {\bf Set}^{\bf B} \). Looking at where
this came from
\[
\Delta_A [F \times G] \cong (\Delta_A [F] \times G) +
(F \times \Delta_A [G]) + (\Delta_A [F] \times \Delta_A [G]) \ ,
\]
we see that the \( G \) is the same for all \( A \), which
means the \( G \) in (3) should be interpreted, as is often done,
to be the functor
\[
{\bf Set}^{\bf A} \to^G {\bf Set}^{\bf B} \to^{{\bf Set}^{P_2}}
\Sett
\]
for \( P_2 \colon {\bf A}^{op} \times {\bf B} \to {\bf B} \)
the second projection, i.e.~\( G \) followed by the
inclusion of \( {\bf Set}^{\bf B} \) in \( \Sett \) given by functors
\( {\bf A}^{op} \times {\bf B} \to {\bf Set} \) constant in
the first variable.

Of course similar remarks go for the \( F \) in the second
term of (3) and the \( F_i \) in (3) of Theorem~\ref{Thm-LimColimDiff}.

\begin{xproposition}
\label{Prop-ScalarsDiff}

For a profunctor \( P \colon {\bf A} \longtod {\bf B} \) we
have
\[
\Delta [P \otimes (\ \ )] (\Phi) \cong P \rlap{\,.}
\]

\end{xproposition}

This is just a restatement of Proposition~\ref{Prop-ScalarsPart}.

We think of \( P \otimes (\ \ ) \) as a linear functor with coefficients
\( P \), and its difference is the constant functor
\[
{\bf Set}^{\bf A} \to \Sett
\]
with constant value \( P \).

\begin{xcorollary}
\label{Cor-DiffOfId}

\[
\Delta [\id_{{\bf Set}^{\bf A}}] (\Phi) = \Id_{\bf A}
\]
where \( \id_{{\bf Set}^{\bf A}} \colon {\bf Set}^{\bf A}
\to {\bf Set}^{\bf A} \) is the identity functor and
\( \Id_{\bf A} \colon {\bf A} \longtod {\bf A} \) is the
identity profunctor.

\end{xcorollary}

Like we did in Proposition~\ref{Prop-GPart}, we can
generalize \ref{Prop-ScalarsDiff} to the following:

\begin{xproposition}
\label{Prop-DiffF}

If \( F \colon {\bf Set}^{\bf A} \to {\bf Set}^{\bf B} \)
preserves binary coproducts, then
\[
\DelF (A, B) = F ({\bf A} (A,-)) (B)
\]
i.e.~\( \DelF = C\!or (F) \), the core of \( F \)
(see Definition~\ref{Def-Core}).

\end{xproposition}

We can improve (2) in Corollary~\ref{Cor-LimColimDiff},
replacing the set \( C \) by a profunctor \( P \colon {\bf B}
\longtod {\bf C} \). Given a tense functor \( F \colon
{\bf Set}^{\bf A} \to {\bf Set}^{\bf B} \), we can compose it
with \( P \otimes (\ \ ) \) to get another tense functor
\[
{\bf Set}^{\bf A} \to^F {\bf Set}^{\bf B} \to^{P \otimes (\ \ )}
{\bf Set}^{\bf C}
\]
which will be called \( P \otimes F \) as its value at \( \Phi \)
is \( P \otimes (F (\Phi)) \) although it might be hard to parse.

\begin{xproposition}
\label{Prop-ConstTimesFunct}

\[
\Delta [P \otimes F] \cong P \otimes \DelF \rlap{\,.}
\]

\end{xproposition}

\begin{proof}

The \( P \otimes F \) is the composite
\[
{\bf Set}^{\bf A} \to^F {\bf Set}^{\bf B} \to^{P \otimes (\ \ )}
{\bf Set}^{\bf C}
\]
so \( \Delta [P \otimes F] \) is the functor \( {\bf Set}^{\bf A}
\to {\bf Set}^{{\bf A}^{op} \times {\bf C}} \) with values
\[
\Delta [P \otimes F](\Phi) (A, C) =
\Delta_A [P \otimes F](\Phi) (C) \rlap{\,.}
\]
On the other hand \( P \otimes \DelF \) is the composite
\[
{\bf Set}^{\bf A} \to^{\DelF} \Sett \to^{P \otimes_{\bf B} (\ \ )}
{\bf Set}^{{\bf A}^{op} \times {\bf C}}
\]
so has values
\[
(P \otimes \DelF) (\Phi) (A, C) = (P \otimes_{\bf B} (\DelF (\Phi))) (A, C) \rlap{\,.}
\]
By the definition of composition of profunctors, this is

\[
\begin{array}{lll}
&    &  \int^B P (B, C) \times \DelF (\Phi) (A, B)\\
& = & \int^B P(B, C) \times \Delta_A [F] (\Phi) (B)\vspace{2mm}\\
& = & (P \otimes \Delta_A [F] (\Phi)) (C)
\end{array}
\]
and by Proposition~\ref{Prop-DeltaProfTimes} this is isomorphic to
\( \Delta_A [P \otimes F] (\Phi) (C) \).
\end{proof}

\subsection{A natural reformulation}

It will be conceptually clearer to reformulate the definition of
\( \Delta \) in more categorical terms, that is, in terms of
natural transformations, Yoneda style. This rids us of many of the
element-based proofs, eliminating, as it does, membership
and especially non-membership. The results are cleaner and
clearer, especially in the next section where we see the
chain rule reduced to composition. This is a vast improvement
over the construction and proof of the one-variable chain rule
given in \citep{Par24} which is far from transparent.

So why not just start with this as a definition? The basic
intuition of finite differences would be lost. It is hard to imagine
why one would define a profunctor using (2) or (3) in
the proposition below, or formulate the product and sum
rules or the chain rule.

\begin{xproposition}
\label{Prop-Yoneda}

Let \( F \colon {\bf Set}^{\bf A} \to {\bf Set}^{\bf B} \) be a tense
functor and \( \Phi \) an object of \( {\bf Set}^{\bf A} \). Then there
is a natural bijection between the following:

\noindent (1) Elements \( x \in \DelF (\Phi) (A, B) \)

\noindent (2) Natural transformations \( t \colon {\bf B} (B, -) \to
F (\Phi + {\bf A}(A, -)) \) giving a pullback
\[
\bfig
\square/>`<-_)`<-_)`>/<850,500>[\ {\bf B}(B,-)\ 
`\ F(\Phi + {\bf A}(A,-))\ 
`\bgr 0\ `\bgr F(\Phi)\ \ ;```]

\place(425,250)[\Pb]

\efig
\]

\noindent (3) Natural transformations \( u \colon F(\Phi) + {\bf B}(B,-) \to
F (\Phi + {\bf A} (A,-)) \) giving a pullback
\[
\bfig
\square/>`<-_)`<-_)`=/<1000,500>[\ F(\Phi) + {\bf B}(B,-)\ 
`\ F(\Phi + {\bf A}(A,-))\ 
`\bgr F (\Phi)\ `\bgr F(\Phi)\ \rlap{\,.};```]

\place(500,250)[\Pb]

\efig
\]

\end{xproposition}

\begin{proof}

An element of \( \DelF (\Phi) (A, B) \) is an element of
\( F (\Phi + {\bf A} (A,-)) (B) \) which is not in \( F (\Phi) (B) \).
By Yoneda, this corresponds bijectively to a natural
transformation
\[
t \colon {\bf B} (B, -) \to F (\Phi + {\bf A}(A,-))
\]
for which \( t (B)(1_B) \notin F (\Phi)(B) \). As \( F (\Phi)\ 
\to/^(->/ F (\Phi + {\bf A}(A,-)) \) is complemented by tenseness,
that's equivalent to none of the values of \( t \) being in
\( F (\Phi) \), which means that
\[
\bfig
\square/>`<-_)`<-_)`>/<850,500>[\ {\bf B}(B,-)\ 
`\ F(\Phi + {\bf A}(A,-))\ 
`\bgr 0\ `\bgr F(\Phi)\ \ ;t```]

\efig
\]
is a pullback. And this, in turn, is equivalent to
\[
\bfig
\square/>`<-_)`<-_)`=/<1000,500>[\ F(\Phi) + {\bf B}(B,-)\ 
`\ F(\Phi + {\bf A}(A,-))\ 
`\bgr F (\Phi)\ `\bgr F(\Phi)\ \rlap{\,.};u```]

\efig
\]
being a pullback, where \( u \) is the inclusion on the first
summand and \( t \) on the second.
\end{proof}

As mentioned in \ref{Def-ProfDef} it is useful to think of the elements of
a profunctor as some sort of morphism but between objects of
different categories (sometimes called heteromorphisms).
Because of the representables appearing in the natural
transformations above, it's not unreasonable to think of them
as morphisms from \( A \) to \( B \), as a kind of Kleisli
morphism although \( F \) is not a monad. If \( F \) were the
identity for example, \( t \) is equivalent to a natural
transformation \( {\bf B}(B,-) \to {\bf A} (A,-) \) so to an
actual morphism \( A \to B \). This is just another way of
saying that \( \Delta [1_{{\bf Set}^{\bf A}}] (\Phi) = \Id_{\bf A} \),
the identity profunctor on \( {\bf A} \), i.e.~the hom functor.

More generally, if \( F \) preserves binary coproducts, a \( t \)
as above corresponds to a natural transformation
\[
{\bf B}(B,-) \to F({\bf A} (A,-))\rlap{\,,}
\]
another way of viewing the identity
\[
\DelF (\Phi) = C\!or (F)
\]
of Proposition~\ref{Prop-DiffF}.

With the natural transformation version of \( \Delta \) it is
easy to see how \( \DelF (\Phi) (A, B) \) is functorial in \( A \)
and \( B \). Given a \( t \) as in (2) and morphisms
\( f \colon A' \to A \) and \( g \colon B \to B' \) we get pullbacks
\[
\bfig
\node a(0,0)[\bgr 0]
\node b(800,0)[\bgr 0]
\node c(1800,0)[\bgr F(\Phi)]
\node d(3100,0)[\bgr F(\Phi)\rlap{\,,}]

\node e(0,500)[{\bf B} (B', -)]
\node f(800,500)[{\bf B}(B,-)]
\node g(1800,500)[F(\Phi + {\bf A}(A, -))]
\node h(3100,500)[F(\Phi + {\bf A}(A',-))]

\arrow|a|/>/[e`f;{\bf B}(g,-)]
\arrow|a|/>/[f`g;t]
\arrow|a|/>/[g`h;F(\Phi + {\bf A}(f, -))]

\arrow/=/[a`b;]
\arrow/>/[b`c;]
\arrow/=/[c`d;]

\arrow/<-_)/[e`a;]
\arrow/<-_)/[f`b;]
\arrow/<-_)/[g`c;]
\arrow/<-_)/[h`d;]

\place(400,250)[\Pb]
\place(1200,250)[\Pb]
\place(2450,250)[\Pb]

\efig
\]
the third one because \( F \) is tense.

Similarly, functoriality in \( \Phi \) is clear. For \( \phi \colon
\Phi \to \Psi \) we get pullbacks
\[
\bfig
\node a(0,0)[\bgr 0]
\node b(900,0)[\bgr F(\Phi)]
\node c(2300,0)[\bgr F(\Psi)]

\node d(0,500)[{\bf B}(B,-)]
\node e(900,500)[F(\Phi + {\bf A} (A,-))]
\node f(2300,500)[F(\Psi + {\bf A}(A,-))]

\arrow/>/[a`b;]
\arrow|b|/>/[b`c;F(\phi)]
\arrow|a|/>/[d`e;t]
\arrow|a|/>/[e`f;F(\phi + {\bf A}(A,-))]

\arrow/<-_)/[d`a;]
\arrow/<-_)/[e`b;]
\arrow/<-_)/[f`c;]

\place(450,250)[\Pb]
\place(1600,250)[\Pb]

\efig
\]
again using tenseness of \( F \).

The same goes for the functoriality in \( F \). If \( \alpha \colon
F \to G \) is a tense transformation, we get pullbacks
\[
\bfig
\node a(0,0)[\bgr 0]
\node b(900,0)[\bgr F(\Phi)]
\node c(2300,0)[\bgr G(\Phi)\rlap{\,.}]

\node d(0,500)[{\bf B}(B,-)]
\node e(900,500)[\alpha (\Phi + {\bf A} (A,-))]
\node f(2300,500)[G (\Phi + {\bf A}(A,-))]

\arrow/>/[a`b;]
\arrow|b|/>/[b`c;\alpha (\Phi)]
\arrow|a|/>/[d`e;t]
\arrow|a|/>/[e`f;\alpha (\Phi + {\bf A}(A,-))]

\arrow/<-_)/[d`a;]
\arrow/<-_)/[e`b;]
\arrow/<-_)/[f`c;]

\place(450,250)[\Pb]
\place(1600,250)[\Pb]

\efig
\]

Showing that \( \DelF \colon {\bf Set}^{\bf A} \to \Sett \) is
tense in this context is probably no easier than the
element-wise proof given for Proposition~\ref{Prop-JacobianFunct}
but it may be more conceptual. It is a result we need if
we want to iterate \( \Delta \), as we do. So we reprove it.

The proof that \( \DelF \) preserves the pullbacks of
complemented subobjects is basically the same as in
\ref{Prop-JacobianFunct} but we reproduce it here without
reference to partical differences or evaluation functors.

Let
\[
\bfig
\square/^(->`>`>`^(->/[\ \Phi_0\ `\ \Phi\ `\ \Psi_0\ `\ \Psi\ ;
``\phi`]

\place(250,250)[\Pb]

\efig
\]
be a pullback of complemented subobjects in \( {\bf Set}^{\bf A} \)
and \( A \) an object of \( {\bf A} \). Consider the four squares in
\( {\bf Set}^{\bf B} \)
\[
\bfig
\node a(0,0)[\ F(\Psi_0 + {\bf A}(A,-))\ ]
\node b(1000,0)[\ F(\Psi + {\bf A}(A,-))\ ]
\node c(0,500)[\ F(\Phi_0 + {\bf A}(A,-))\ ]
\node d(1000,500)[\ F(\Phi + {\bf A}(A,-))\ ]
\node e(0,1000)[\ \Delta{[}\Phi_0{]}(A,-)\ ]
\node f(1000,1000)[\ \Delta {[}\Phi{]}(A,-)\ ]

\arrow/^(->/[a`b;]
\arrow/^(->/[c`d;]
\arrow/^(->/[e`f;]

\arrow/>/[c`a;]
\arrow/>/[e`c;]
\arrow/>/[f`d;]
\arrow/>/[d`b;]

\place(500,250)[(2)]
\place(500,750)[(1)]

\node a'(2000,0)[\ F(\Psi_0 + {\bf A}(A,-))\ ]
\node b'(3000,0)[\ F(\Psi + {\bf A}(A,-))\ \rlap{\,.}]
\node c'(2000,500)[\Delta{[}\Psi_0{]}(A, -)\ ]
\node d'(3000,500)[\ \Delta {[}\Psi{]}(A, -)\ ]
\node e'(2000,1000)[\ \Delta{[}\Phi_0{]}(A,-)\ ]
\node f'(3000,1000)[\ \Delta {[}\Phi{]}(A,-)\ ]

\arrow/^(->/[a'`b';]
\arrow/^(->/[c'`d';]
\arrow/^(->/[e'`f';]

\arrow/>/[c'`a';]
\arrow/>/[e'`c';]
\arrow/>/[f'`d';]
\arrow/>/[d'`b';]

\place(2500,750)[(3)]
\place(2500,250)[(4)]

\efig
\]
(1) and (4) are pullbacks by definition of \( \Delta \) and
(2) because \( F \) is tense. As the pasted rectangle
(1) + (2) is equal to (3) + (4), we get that (3) is also a
pullback.

As \( \DelF \) preserves pullbacks of complemented
subobjects, it will take a complemented subobject
\( \Phi_0 \ \to/^(->/ \Phi \) to a mono, but we still have to
prove that it's complemented. We have to prove that
for any \( f \colon A' \to A \) and \( g \colon B \to B' \),
\[
\bfig
\square/ >->`>`>` >->/<850,500>[\Delta{[}\Phi_0{]}(A,B)
`\Delta{[}\Phi{]}(A,B)
`\Delta{[}\Phi_0{]} (A',B')
`\Delta{[}\Phi{]}(A',B');
`\Delta{[}\Phi_0{]}(f,g)
`\Delta{[}\Phi{]}(f,g)`]

\efig
\]
is a pullback.

An element of \( \Delta[\Phi] (A,B) \) is a
natural transformation
\( t \colon {\bf B}(B,-) \to
F (\Phi + {\bf A}(A,-))
\).
To be in \( \Delta[\Phi_0] (A,B) \)
means that it factors through \( F (\Phi_0 + {\bf A}(A,-))
\ \to/^(->/ F (\Phi + {\bf A} (A, -)) \). Referring to the following
diagram
\[
\bfig
\node a(0,0)[{\bf B} (B',-)]
\node b(800,500)[\bgr F  \Phi_0 + {\bf A}(A,-))]
\node c(2500,500)[\bgr F (\Phi_0 + {\bf A} (A',-))]
\node d(0,1050)[{\bf B}(B,-)]
\node e(800,1050)[F(\Phi + {\bf A}(A,-))]
\node f(2500,1050)[F(\Phi + {\bf A} (A',-))]

\arrow|b|/>/[a`c;u]
\arrow|l|/>/[a`d;{\bf B}(g,-)]
\arrow|a|/>/[a`b;u']
\arrow|a|/>/[d`b;u'']
\arrow|b|/>/[b`c;F(\Phi_0 + {\bf A}(f,-))]
\arrow|a|/>/[d`e;t]
\arrow|a|/>/[e`f;F(\Phi + {\bf A}(f,-))]
\arrow/^(->/[b`e;]
\arrow/^(->/[c`f;]

\place(1650,800)[\Pb]

\efig
\]
\( \Delta [\Phi] (f,g) (t) \) is the composite of the left arrow
with the two top arrows, and to say that it is in
\( F(\Phi_0 + {\bf A}(A',-)) \) means that there is a \( u \)
making the outside boundary commute. The square in a
pullback because \( F \) is tense so there exists a unique
\( u' \) as shown and as \( F (\Phi_0 + {\bf A}(A,-)) \) is
complemented there exists a \( u'' \) by
Proposition~\ref{Prop-FunCatCompl}. So \( t \) factors
through \( F(\Phi_0 + {\bf A}(A,-)) \) which is what we wanted.

\subsection{Lax chain rule}

We saw in \citep{Par24} that the chain rule for the single
variable functorial difference was expressed as a laxity
morphism rather than an isomorphism, and the same applies
in the multivariable case. For tense functors
\( F \colon {\bf Set}^{\bf A} \to {\bf Set}^{\bf B} \) and
\( G \colon {\bf Set}^{\bf B} \to {\bf Set}^{\bf C} \) we will
construct a comparison transformation
\[
\gamma (\Phi) \colon \Delta [G] (F(\Phi)) \otimes_{\bf B}
\DelF (\Phi) \to \Delta[GF](\Phi)
\]
and establish associativity and unit laws for it. In fact,
considering \( \DelF (\Phi) \) as a profunctor may not
mean much unless it composes like a profunctor.
Otherwise it is just an object of \( \Sett \).

The construction of \( \gamma \) in \citep{Par24} is perhaps
a bit opaque and the profunctor interpretation clarifies
this. We'll see that it is, in a sense, just composition as it
should be.

In the previous section we described the functoriality of
\( \Delta \) in terms of the characterization (2) of
Proposition~\ref{Prop-Yoneda}, but for the chain rule
the characterization (3) is better, so we reformulate the
functorialities in this context. As we will refer to it a lot,
let us call a natural transformation \( t \) such that
\[
\bfig
\square/>`<-_)`<-_)`=/<1000,500>[F(\Phi) + {\bf B}(B,-)
`F(\Phi + {\bf A}(A,-))
`\bgr F(\Phi)`\bgr F(\Phi);t```]

\efig
\]
is a pullback, a PPI transformation (for pullback preserves
injections).

Functoriality of \( \DelF (\Phi) (A,B) \), considered as a set
of PPI transformations, is easy. It's just composition with
\( F (\Phi + {\bf A} (f, -)) \) and \( F (\Phi) + {\bf B}(g,-) \)
respectively.

Functoriality in \( \Phi \) and \( F \) are a {\em bit} more
complicated as the \( \Phi \) and \( F \) appear in both the
domain and codomain of \( t \). The following characterization
will be useful, although it is nothing but a reformulation.

\begin{xproposition}
\label{Prop-ReformFun}

Let \( t \colon F (\Phi) + {\bf B} (B,-) \to F (\Phi + {\bf A}(A,-)) \)
be a PPI transformation.

\noindent (1) If \( \phi \colon \Phi \to \Psi \) is a natural
transformation, then \( \DelF (\phi) (A,B) (t) \) is the unique
PPI transformation \( t' \) such that 
\[
\bfig
\square<1000,500>[F(\Phi) + {\bf B}(B,-)
`F(\Phi + {\bf A} (A,-))
`F(\Psi) + {\bf B}(B,-)
`F(\Psi + {\bf A} (A,-))\rlap{\,.};
t`F(\phi) + {\bf B}(B,-)`F(\phi + {\bf A}(A,-))`t']

\efig
\]

\noindent (2) If \( \alpha \colon F \to G \) is a tense transformation,
then \( \Delta [\alpha] (\Phi) (A,B) (t) \) is the unique PPI
transformation \( t'' \) such that
\[
\bfig
\square<1000,500>[F(\Phi) + {\bf B}(B,-)
`F(\Phi + {\bf A} (A,-))
`F(\Psi) + {\bf B}(B,-)
`F(\Psi + {\bf A} (A,-))\rlap{\,.};
t`F(\alpha) + {\bf B}(B,-)`F(\alpha + {\bf A}(A,-))`t'']

\efig
\]
 
\end{xproposition}

\begin{xtheorem}
\label{Thm-ChRule}

For tense functors \( F \colon {\bf Set}^{\bf A} \to {\bf Set}^{\bf B} \) and
\( G \colon{\bf Set}^{\bf B} \to {\bf Set}^{\bf C} \) there is a natural
transformation
\[
\gamma \colon (\Delta [G] \circ F) \otimes \DelF \to \Delta [GF]
\]
which is:
\begin{itemize}

	\item[(1)] natural in \( F \) and \( G \)
	
	\item[(2)] associative
	
	\item[(3)] normal (invertible unitors)

\end{itemize}
\end{xtheorem}

\begin{proof}

\( \gamma \) is to be understood pointwise, i.e.~as a profunctor morphism
\[
\gamma (\Phi) \colon \Delta [G] (F (\Phi)) \otimes_{\bf B}
\DelF (\Phi) \to \Delta [GF] (\Phi)
\]
\[
\bfig
\Atriangle/@{<-}|{\bb}`@{>}|{\bb}`@{>}|{\bb}/[{\bf B}
`{\bf A}`{\bf C};
\DelF(\Phi)`\Delta {[}G{]} (F (\Phi))
`\Delta{[}FG{]}(\Phi)]

\morphism(500,300)|r|/=>/<0,-150>[`;\gamma(\Phi)]

\efig
\]
for each \( \Phi \in {\bf Set}^{\bf A} \), and furthermore natural
in that \( \Phi \).

Let \( A \in {\bf A} \) and \( C \in {\bf C} \). An element of
\[
\big(\Delta [G] (F(\Phi)) \otimes_{\bf B} \DelF (\Phi)\big) (A,C)
\]
is an equivalence class
\[
u \otimes_B t = [A \todd{t}{} B \todd{u}{} C]
\]
where \( u \) and \( t \) are PPI transformations. Let
\( \gamma (\Phi)(A,C) (u \otimes_B t) = G t \cdot u \) which is
indeed PPI:
\[
\bfig
\node a(0,0)[\bgr GF(\Phi)]
\node b(1200,0)[\bgr GF(\Phi)]
\node c(2400,0)[\bgr GF(\Phi)\rlap{\,.}]
\node d(0,500)[GF(\Phi) + {\bf C} (C,-)]
\node e(1200,500)[G(F(\Phi) + {\bf B}(B,-))]
\node f(2400,500)[GF(\Phi + {\bf A}(A,-))]

\arrow/<-_)/[d`a;]
\arrow/<-_)/[e`b;]
\arrow/<-_)/[f`c;]

\arrow/=/[a`b;]
\arrow/=/[b`c;]
\arrow|a|/>/[d`e;u]
\arrow|a|/>/[e`f;Gt]

\place(600,250)[\Pb]
\place(1800,250)[\Pb]

\efig
\]

We must show that \( \gamma (\Phi) (A,C) \) is well-defined.
Suppose we have another pair of transformation related by
a single morphism
\[
\bfig
\square/@{>}|{\bb}`=`>`@{>}|{\bb}/[A`B`A`B';t``g`t']

\square(500,0)/@{>}|{\bb}``=`@{>}|{\bb}/[B`C`B'`C\rlap{\,.};
u```u']

\efig
\]
This means that we have commutative squares
\[
\bfig
\square/>`=`<-`>/<1000,500>[GF(\Phi) + {\bf C}(C,-)
`G(F(\Phi) + {\bf B}(B,-))
`GF(\Phi) + {\bf C}(C,-)
`G(F(\Phi) + {\bf B}(B',-));
u``G(F(\Phi) + {\bf B}(g,-))`u']

\square(2100,0)/>`<-`=`>/<1000,500>[F(\Phi) + {\bf B}(B,-)
`F(\Phi + {\bf A}(A,-))
`F(\Phi) + {\bf B}(B',-)
`F(\Phi + {\bf A}(A,-))\rlap{\ .};
t`F(\Phi) + {\bf B}(g,-)``t']

\efig
\]
If we apply \( G \) to the second and paste it to the first we
get a commutative diagram which shows that
\( G t \cdot u = G t' \cdot u' \). It follows that
\( \gamma (\Phi) (A,C) \) is well-defined.

Naturality in \( A \) and \( C \) is clear as it is just composition
with \( F (\Phi + {\bf A} (f, -)) \) and \( F (\Phi) + {\bf C}(h,-) \)
respectively and has nothing to do with the equivalence
relation, which is localized at \( B \). So we get a profunctor
morphism \( \gamma (\Phi) \).

To show that \( \gamma \) is natural in \( \Phi \), let
\( \phi \colon \Phi \to \Psi \) be a natural transformation and
consider
\[
\bfig
\square<1500,500>[\Delta{[}G{]} (F(\Phi)) \otimes_{\bf B} \DelF (\Phi)(A,C)
`\Delta{[}GF{]}(\Phi)(A,C)
`\Delta{[}G{]}(F(\Psi))\otimes_{\bf B} \DelF (\Psi)(A,C)
`\Delta {[}GF{]} (\Psi)(A,C);
\gamma(\Phi)```\gamma(\Psi)]

\efig
\]
where the vertical arrows are induced by \( \phi \). If we chase
an element \( u \otimes_B t \) in the domain, first around the
left-bottom we get \( u' \otimes_B t' \) and then \( G t' \cdot u' \)
where \( u' \) and \( t' \) are the unique PPI's such that
\[
\bfig
\square<1100,500>[GF(\Phi) + {\bf C}(C,-)
`G(F(\Phi) + {\bf B} (B,-))
`GF(\Psi) + {\bf C}(C,-)
`G(F(\Psi) + {\bf B} (B,-));u`GF(\phi) + {\bf C}(C,-)
`G(F(\phi) + {\bf B}(B,-))`u']

\efig
\]
\[
\bfig
\square<1100,500>[F(\Phi) + {\bf B}(B,-)
`F(\Phi + {\bf A} (A,-))
`F(\Psi) + {\bf B}(B,-)
`F(\Psi + {\bf A}(A,-))\rlap{\,.};
t`F(\phi) + {\bf B}(B,-)
`F(\phi + {\bf A}(A,-))`t']

\efig
\]
On the other hand, going around the top-right we get \( G t \cdot u \) and then
\( v' \) the unique PPI such that
\[
\bfig
\square<1150,500>[GF(\Phi) + {\bf C}(C,-)
`GF(\Phi + {\bf A}(A,-))
`GF(\Psi) + {\bf C} (C,-)
`GF(\Psi + {\bf A}(A,-))\rlap{\,.};
G t \cdot u
`GF(\phi) + {\bf C}(C,-)
`GF(\phi + {\bf A}(A,-))
`v']

\efig
\]
If we apply \( G \) to the diagram for \( t' \) above and
paste it to the one for \( u' \), we see that \( G t' \cdot u' \)
is such a \( v' \), and so \( v' = G t' \cdot u' \). This gives
naturality in \( \Phi \).

We can check naturality in \( F \) and \( G \) separately.
First, let \( \alpha \colon F \to F' \) be a tense natural
transformation. We wish to show that
\[
\bfig
\square<1500,500>[\Delta{[}G{]} (F(\Phi)) \otimes_{\bf B} \DelF (\Phi)(A,C)
`\Delta{[}GF{]}(\Phi)(A,C)
`\Delta{[}G{]}(F'(\Phi))\otimes_{\bf B} \Delta{[}F'{]} (\Phi)(A,C)
`\Delta {[}GF'{]} (\Phi)(A,C);
\gamma(\Phi)```\gamma(\Phi)]

\efig
\]
commutes.
\( \alpha \) acting on an element \( u \otimes_B t \) of the domain
gives \( u' \otimes_B t' \) which gets sent to \( G t' \cdot u' \), where
\[
\bfig
\square<1100,500>[GF(\Phi) + {\bf C}(C,-)
`G(F(\Phi) + {\bf B}(B,-))
`GF'(\Phi) + {\bf C}(C,-)
`G(F'(\Phi) + {\bf B}(B,-));
u`GF(\alpha) + {\bf C}(C,-)
`G(\alpha(\Phi) + {\bf B}(B,-))`u']

\efig
\]

\vspace{2mm}

\[
\bfig
\square<1100,500>[F(\Phi) + {\bf B}(B,-)
`F(\Phi + {\bf A}(A,-))
`F'(\Phi) + {\bf B}(B,-)
`F'(\Phi + {\bf A} (A,-))\rlap{\,.};
t`\alpha(\Phi) + {\bf B}(B,-)
`F(\alpha + {\bf A} (A,-))`t']

\efig
\]
On the other hand we first get \( G t \cdot u \) and
then \( v' \) such that
\[
\bfig
\square<1100,500>[GF(\Phi) + {\bf C}(C,-)
`GF(\Phi + {\bf A}(A,-))
`GF'(\Phi) + {\bf C}(C,-)
`GF(\Phi + {\bf A}(A,-))\rlap{\,.};
G t \cdot u`G\alpha(\Phi) + {\bf C}(C,-)
`GF(\alpha + {\bf A}(A,-))`v']

\efig
\]
Again, applying \( G \) to the square for \( t' \) and
pasting to the one for \( u' \), we see that \( v' = G u' \cdot t' \),
i.e.~naturality in \( F \).

For naturality in \( G \), let \( \beta \colon G \to G' \) be a
tense natural transformation. We'll show that
\[
\bfig
\square<1500,500>[\Delta{[}G{]} (F(\Phi))\otimes_{\bf B}\DelF(\Phi)(A,C)
`\Delta{[}GF{]}(\Phi) (A,C)
`\Delta{[}G'{]} (F(\Phi)) \otimes_{\bf B} \DelF(\Phi)(A,C)
`\Delta{[}G'F{]}(\Phi)(A,C);
\gamma(\Phi)```\gamma(\Phi)]

\efig
\]
commutes.
An element \( u \otimes t \) of the domain, goes down to
\( u' \otimes t \) and then \( G' u' \cdot t \) for \( u' \) such that
\[
\bfig
\square<1200,500>[GF(\Phi) + {\bf C} (C,-)
`G(F(\Phi) + {\bf B} (B,-))
`G'F(\Phi) + {\bf C}(C,-)
`G'(F(\Phi) + {\bf B} (B,-))\rlap{\,.};
u`\beta F(\Phi) + {\bf C}(C,-)`\beta (F(\Phi) + {\bf B}(B,-))`u']

\efig
\]
\( u \otimes t \) goes across to \( G t \cdot u \) and then down
to \( v' \) such that 
\[
\bfig
\square<1200,500>[GF(\Phi) + {\bf C} (C,-)
`GF(\Phi + {\bf A} (A,-))
`G'F(\Phi) + {\bf C}(C,-)
`G'F(\Phi + {\bf A} (A,-))\rlap{\,.};
G t \cdot u`\beta F(\Phi) + {\bf C}(C,-)`\beta (F(\Phi) + {\bf A}(A,-))`v']

\efig
\]
If we paste the diagram for \( u' \) with the naturality square
\[
\bfig
\square<1200,500>[G(F(\Phi) + {\bf B} (B,-))
`GF(\Phi + {\bf A}(A,-))
`G'(F(\Phi) + {\bf B}(B,-))
`G'F(\Phi + {\bf A} (A,-));
Gt`\beta (F \Phi + {\bf B}(B,-))`\beta F(\Phi + {\bf A}(A,-))
`G't]

\efig
\]
and compare with the diagram for \( v' \) we see that
\( v' = G' t \cdot u' \), which gives naturality in \( G \).

Let
\[
{\bf Set}^{\bf A} \to^F {\bf Set}^{\bf B} \to^G {\bf Set}^{\bf C}
\to^H {\bf Set}^{\bf D}
\]
be tense functors. Associativity involves taking an element
\( v \otimes u \otimes t \) of
\[
\Delta [H] (GF(\Phi)) \otimes_{\bf B} \Delta [G] (F\Phi)
\otimes_{\bf C} \DelF (\Phi)
\]
at \( (A,D) \) and applying \( \gamma \) in two different
ways to reduce it to elements of \( \Delta [HGF](\Phi) \),
and seeing that they are equal. This is for any PPI
transformations
\begin{align*}
t \colon F (\Phi) + {\bf B} (B,-) & \to F (\Phi + {\bf A} (A,-))\\
u \colon GF (\Phi) + {\bf C} (C,-) & \to G(F(\Phi) + {\bf B} (B,-))\\
v \colon HGF (\Phi) + {\bf D} (D,-) &\to H (GF(\Phi) + {\bf C} (C,-))\rlap{\,.}
\end{align*}
And indeed, we get
\[
\bfig
\square/|->`|->``|->/<1000,700>[v \otimes u \otimes t
`v \otimes (G t \cdot u)
`(H u \cdot v) \otimes t
`H G t \cdot H u \cdot v\rlap{\,.};```]

\morphism(1000,700)/|->/<0,-350>[v \otimes (G t \cdot u)`
H (G t \cdot u) \cdot v;]

\morphism(1000,270)/=/<0,-150>[`;]

\efig
\]

For the unit laws, first assume that \( {\bf B} = {\bf A} \)
and that \( F = \id_{{\bf Set}^{\bf A}} \). Then \( \gamma (\Phi) \)
takes the form
\[
\gamma (\Phi) \colon \Delta [G] \otimes_{\bf A}
\Delta [\id_{{\bf Set}^{\bf A}}] (\Phi) \to \Delta [G] (\Phi)
\]
and an element of the domain is an equivalence class
\( u \otimes t \) for PPI's
\[
\Phi + {\bf A} (A',-) \to^t \Phi + {\bf A} (A,-) \quad
G(\Phi) + {\bf C} (C,-) \to^u G(\Phi + {\bf A}(A',-))\rlap{\,.}
\]
For \( t \) to be a PPI it must be of the form
\[
\Phi + {\bf A}(A',-) \to^{\Phi + {\bf A} (f,-)} \Phi + {\bf A}(A,-)
\]
and every equivalence class has a unique representative
where \( f \) is \( 1_A \). Then \( \gamma (\Phi)(u \otimes 1) = u \)
gives our bijective right unitor.

For the left unitor, let \( {\bf B} = {\bf C} \) and \( G = \id_{{\bf Set}^{\bf C}} \).
Then \( \gamma \) takes the form
\[
\gamma (\Phi) \colon \Delta [\id_{{\bf Set}^{\bf C}}] (F (\Phi)
\otimes \DelF (\Phi)) \to \DelF (\Phi)
\]
and an element of its domain is an equivalence class
\( u \otimes t \) with PPI's
\[
F(\Phi) + {\bf C}(C',-) \to^t F(\Phi + {\bf A} (A,-))\quad
F(\Phi) + {\bf C} (C,-) \to^u F(\Phi) + {\bf C} (C',-)
\]
For \( u \) to be a PPI it must be of the form \( F(\Phi) +
{\bf C}(g,-) \). Again every equivalence class contains a
unique representative with \( g = 1_C \). Then
\[
\gamma (\Phi) (1 \otimes t) = t
\]
gives the bijective unitor.
\end{proof}

As stated, the lax chain rule is called lax just because what
might have been hoped to be an isomorphism is merely a
comparison morphism reducing a more complicated expression
to a simpler one. But, if we reformulate it in terms of the
tangent bundle of Section~\ref{SSec-JacobianDef}, we get
an actual lax normal functor.

Recall that the tangent functor \( T [F] \)
\[
\bfig
\square<1000,500>[{\bf Set}^{\bf A} \times {\bf Set}^{\bf A}
`{\bf Set}^{\bf B} \times {\bf Set}^{\bf B}
`{\bf Set}^{\bf A}`{\bf Set}^{\bf A};
T{[}F{]}`P_1`P_1`F]

\efig
\]
is given by
\[
T [F] (\Phi, \Psi) = (F(\Phi), \DelF (\Phi) \otimes_{\bf A} \Psi) \rlap{\,.}
\]
If \( G \colon {\bf Set}^{\bf B} \to {\bf Set}^{\bf C} \) is another
tense functor, then the composite
\[
T [G] \circ T[F] = (GF(\Phi), \Delta[G] (F(\Phi)) \otimes_{\bf B}
\DelF (\Phi) \otimes_{\bf A} \Phi)
\]
and
\[
(1_{GF(\Phi)}, \gamma (\Phi) \otimes_{\bf A} \Psi) \colon
T [G] \circ T[F] \to T[GF]
\]
makes \( T \colon \Tense \to \Tense \) into a lax normal
functor. We omit the details which only involve the
rearrangement of the facts proved in Theorem~\ref{Thm-ChRule}.

\section{Newton series}
\label{Sec-Newton}

\subsection{Multivariable Newton series}
\label{SSec-Newton}

The Newton series of a function of a real variable
\( f \colon {\mathbb R} \to {\mathbb R} \) is a discrete version
of Taylor series. Its aim is to recover \( f \) from its iterated
differences, or to approximate \( f \) by polynomials. The
formula is well-known
\begin{align*}
 & \sum_{n = 0}^\infty \frac{\Delta^n {[}f{]} (0)}{n!} x^{\downarrow n}\\
 =\  & \sum_{n = 0}^\infty \Delta^n {[}f{]} (0) \binom{x}{n}
\end{align*}
when \( x^{\downarrow n} \) is the falling power
\( x (x - 1) \dots (x - n + 1) \) and \( \binom{x}{n} \)
is the ``binomial coefficient'' 
\( \frac{x (x - 1) \dots (x - n + 1)}{n!} \).

Although not so well-known, a recursive argument
produces a multivariable version: for \( f \colon {\mathbb R^n}
\to {\mathbb R} \) we have
\begin{align*}
 & \sum_{k_1, k_2, \ldots,  k_n = 0}^\infty
 \frac{\Delta^{k_1}_{x_1}  \Delta^{k_2}_{x_2} \cdots \Delta^{k_n}_{x_n}
 {[}f{]} (0,\ldots, 0)}{k_1 ! k_2 ! \cdots k_n !}
 x_1^{\downarrow k_1} x_2^{\downarrow k_2} \cdots x_n^{\downarrow k_n}\\
 = \ & \sum_{k_1, k_2, \ldots,k_n = 0}^\infty \Delta^{k_1}_{x_1}  \Delta^{k_2}_{x_2} \cdots \Delta^{k_n}_{x_n}
 {[}f{]} (0,\ldots, 0)
 \binom{x_1}{k_1} \binom{x_2}{k_2} \cdots \binom{x_n}{k_n}\rlap{\ .}
\end{align*}

In \citep{Par24} we gave a categorified version for taut
endofunctors of \( {\bf Set} \) and showed that for analytic
functors their Newton series converge to them. In fact
this holds for a larger class of taut functors, which we
call soft analytic. Not only that, the approximation alluded
to above manifests itself as a categorical adjointness.
In this section we develop multivariable versions of these
results.

\subsection{Soft multivariable analytic functors}
\label{SSec-SoftAnalytic}

In order to categorify multivariable Newton series we
must modify the notion of \( {\bf A} \)-\( {\bf B} \) symmetric
sequence to take into account the extra structure that
the iterated differences have. We replace the category
\( ! {\bf A} \) of \citep{FioGamHylWin08} by the larger category \( {\downarrow}{\bf A} \)
with the same objects, finite sequences \( \langle A_1, \dots , A_n \rangle \)
of objects of \( {\bf A} \), but where the morphisms
\[
\langle A_1, \dots , A_n \rangle \to \langle C_1, \dots C_m \rangle
\]
are pairs \( (\sigma, \langle f_j \rangle) \) such that \( \sigma \colon 
m \to n \) is a {\em surjection} and \( \langle f_j \rangle \) is a family
of morphisms indexed by \( m \)
\[
f_j \colon A_{\sigma j} \to C_j \rlap{\ .}
\]
Composition is formally the same as for \( ! {\bf A} \)
\[
(\tau, \langle g_k \rangle) (\sigma, \langle f_j \rangle) =
(\sigma \tau, \langle g_k f_{\tau k} \rangle)\rlap{\ .}
\]

Whereas \( ! {\bf A} \) is the free symmetric strict monoidal
category generated by \( {\bf A} \), \( {\downarrow}{\bf A} \) is the
free symmetric monoidal category in which every object
has a canonical cocommutative coassociative
comultiplication.

\begin{xdefinition}
\label{Def-SoftSS}

A {\em soft} \( {\bf A} \)-\( {\bf B} \)-{\em symmetric sequence}
is a profunctor \( P \colon{\downarrow}{\bf A} \longtod {\bf B} \).

\end{xdefinition}

Given a soft \( {\bf A} \)-\( {\bf B} \)-symmetric sequence
\(P\colon {\downarrow}{\bf A} \longtod {\bf B} \) we define
the functor \( \widetilde{P} \colon {\bf Set}^{\bf A} \to
{\bf Set}^{\bf B} \) by the formula
\[
\widetilde{P} (\Phi)(B) = \int^{\langle A_1 \dots  A_n\rangle \in {\downarrow}{\bf A}}
P(A_1, \dots , A_n; B) \times \Phi A_1 \times {\dots} \times \Phi A_n \rlap{\ .}
\]
Of course, for this to make sense \( \Phi A_1 \times {\dots} \times
\Phi A_n \) must come from a functor \( {\downarrow}{\bf A} \to {\bf Set} \),
which is indeed the case. For a morphism
\[
(\sigma, \langle f_1 \dots f_m\rangle) \colon
\langle A_1, \dots , A_n \rangle \to
\langle C_1, \dots , C_m \rangle
\]
we have a unique morphism making
\[
\bfig
\square<1100,500>[\Phi A_1 \times {\dots} \times \Phi A_n
`\Phi C_1 \times {\dots} \times \Phi C_m
`\Phi A_{\sigma j}`\Phi C_j;
`{\rm proj}_{\sigma j}`{\rm proj}_j`\Phi f_j]

\efig
\]
commute for all \( j \in m \).

A more conceptual description of \( \widetilde{P} \) is
in terms of Kan extensions. Let \( Q \colon ({\downarrow}{\bf A})^{op}
\to {\bf Set}^{\bf A} \) be the functor defined by
\[
Q \langle A_1 {\dots} A_n \rangle = {\bf A} (A_1, -) +
{\dots} + {\bf A} (A_n, -) \rlap{\ .}
\]
It is indeed a functor, its value on a morphism
\[ (\sigma, \langle f_1, {\dots} , f_m \rangle) \colon
\langle A_1, {\dots} , A_n\rangle \to
\langle C_1, {\dots} , C_m \rangle \] 
being the unique morphism making all the squares
\[
\bfig
\square<1500,500>[{\bf A}(C_j, -)`{\bf A} (A_{\sigma j}, -)
`{\bf A} (C_1, -) + {\dots} + {\bf A} (C_m, -)
`{\bf A} (A_1, -) + {\dots} + {\bf A} (A_n, -);
{\bf A} (f_j, -)`{\rm inj}_j`{\rm inj}_{\sigma j}`]

\efig
\]
commute. A profunctor \( P \colon {\downarrow}{\bf A}
\longtod {\bf B} \) is a functor \( P \colon ({\downarrow}{\bf A})^{op}
\times {\bf B} \to {\bf Set} \) which may be alternately
described as a functor \( ({\downarrow}{\bf A})^{op} \to
{\bf Set}^{\bf B} \) (which we denote by the same letter).
Then \( \widetilde{P} \) is the left Kan extension of
\( P \) along \( Q \):
\[
\bfig
\Vtriangle[({\downarrow}{\bf A})^{op}
`{\bf Set}^{\bf A}`{\bf Set}^{\bf B}\rlap{\ .};
Q`P`{\rm Lan}_Q P = \widetilde{P}]

\morphism(420,300)/=>/<180,0>[`;\eta]

\efig
\]
Indeed,
\[
{\rm Lan}_Q P (\Phi) = \int^{A_1 {\dots} A_n}
P (A_1 {\dots} A_n ; -) \times {\bf Set}^{\bf A}
(Q \langle A_1, {\dots} , A_n \rangle, \Phi)
\]
(see \citet{Mac71}, p. 236) and \( {\bf Set}^{\bf A}
(Q \langle A_1 {\dots} A_n \rangle , \Phi) \cong
\Phi A_1 \times {\dots} \times \Phi A_n \).

\( Q \) may be considered as a profunctor \( {\downarrow}{\bf A}
\longtod {\bf A} \) and we have the following ``softening'' of
Proposition~\ref{Prop-AnalComposite}.

\begin{xproposition}
\label{Prop-PTilde}

\begin{itemize}

	\item[1.] \( \widetilde{P} \) is the composite \( P \otimes (Q \obslash (\ )) \)
\[
{\bf Set}^{\bf A} \to^{Q \obslash (\ )} {\bf Set}^{\downarrow{\bf A}}
\to^{P \otimes (\ )} {\bf Set}^{\bf B} \rlap{\ .}
\]

	\item[2.] \( Q \) satisfies the condition of 2.4.1.

\end{itemize}

\end{xproposition}

\begin{proof}

(1) Same as in \ref{Prop-AnalComposite}.

(2) Again \( \pi_0 Q (A_1, {\dots} , A_n ; -) = n \) for the
same reason (sum of \( n \) representables), but now for
a morphism \( (\sigma, \langle f_1 , {\dots} , f_m \rangle)
\colon \langle A_1, {\dots} , A_n \rangle \to
\langle C_1 , {\dots} , C_m \rangle \)
the morphism
\[
\pi_0 Q (C_1, {\dots} , C_m ; -) \to
\pi_0 Q (A_1 , {\dots} , A_n ; -)
\]
is
\(
\sigma \colon m \to n
\),
which is onto.
\end{proof}

\begin{xcorollary}
\label{Cor-PTilde}

\( \widetilde{P} \) is tense.

\end{xcorollary}

A more elementary understanding of \( \widetilde{P} \)
will be useful. From the coend formula for Kan extension
we see that an element of \( \widetilde{P} (\Phi) (B) \)
is an equivalence class of pairs \( (p, \phi) \)
\begin{center}
\(
\big[p \colon \langle A_1, {\dots} , A_n \rangle \longtod B, 
\phi \colon \sum {\bf A} (A_i, -) \to \Phi\big]
\)
\end{center}
where \( p \in P (A_1, {\dots} , A_n ; B) \) and
\( \sum {\bf A} (A_i, -) \) is short for \( \sum_{i = 1}^n
{\bf A} (A_i, -) \). The equivalence relation is generated
by identifying \( (p, \phi) \) and \( (q, \psi) \) when there
is a morphism \( (\sigma, \langle f_j \rangle) \colon
\langle A_1 , {\dots} , A_n \rangle \to
\langle C_1, {\dots} , C_m \rangle \) in \( {\downarrow} {\bf A} \)
such that
\[
\bfig
\Dtriangle/>`@{>}|{\bb}`@{<-}|{\bb}/<600,350>[\langle A_1 , {\dots} , A_n\rangle
`B
`\langle C_1 , {\dots} , C_m \rangle;
(\sigma, \langle f_j \rangle)`p`q]

\Dtriangle(1400,0)/<-`>`<-/<600,350>[\sum {\bf A} (A_i, -)`\Phi
`\sum {\bf A} (C_j , -);
\sum_\sigma {\bf A} (f_j, -)`\phi`\psi]

\efig
\]
where \( \sum_\sigma {\bf A} (f_j, -) \) represents the natural
transformation taking \( g \colon C_j \to A \) to
\( A_{\sigma (j)} \to^{f_j} C_j \to^g A \).

Functoriality of \( \widetilde{P} \) in \( B \) and \( \Phi \) is
by composition: for \( b \colon B \to B' \)
\[
\widetilde{P} (\Phi) (b) \colon (p, \phi) \longmapsto (b p, \phi)
\]
and for \( \theta \colon \Phi \to \Psi \)
\[
\widetilde{P} (\theta)(B) \colon (p, \phi) \longmapsto (p, \theta \phi) \rlap{\ .}
\]

The universal property of Kan extensions says that for any
functor \( F \colon {\bf Set}^{\bf A} \to {\bf Set}^{\bf B} \)
we have a natural bijection
\[
\widetilde{P} \to^t F \over P \to^u F Q \rlap{\ .}
\]
The correspondence between \( t \) and \( u \) is the
following. \( t \colon \widetilde{P} \to F \) is given by
a family of natural transformations
\[
\langle \widetilde{P} (\Phi) \to F (\Phi)\rangle_\Phi
\]
natural in \( \Phi \in {\bf Set}^{\bf A} \), which further breaks
down into a doubly indexed family of functions
\[
\langle \widetilde{P} (\Phi) (B) \to F (\Phi) (B) \rangle_{\Phi, B}
\]
natural in both \( \Phi \) and \( B \). So for every equivalence
class 
\[ [p \colon \langle A_1, {\dots} , A_n \rangle \longtod B,\allowbreak
\phi \colon \sum {\bf A} (A_i, -) \to \Phi] \] 
we get an element
\( t [p, \phi] \in F (\Phi) (B) \).

On the other hand \( u \colon P \to F Q \) is a doubly
indexed family of functions
\begin{center}
\( \langle P(A_1, {\dots}, A_n ; B) \to F (\sum {\bf A}
(A_i, -)) (B) \rangle \)
\end{center}
natural in \( \langle A_1, {\dots}, A_n\rangle \in\  {\downarrow}{\bf A} \)
and \( B \) in \( {\bf B} \).

Given \( t \) we get \( u \) by restricting to the case
\( \Phi = \sum {\bf A} (A_i, -) \) and \( \phi \) the identity
\[
u (p) = t [p, \id_{\sum {\bf A} (A_i, -)}] \rlap{\ .}
\]
Given \( u \) we get \( t \) by
\[
t [p, \phi] = F (\phi) (u (p)) \rlap{\ .}
\]
There is nothing to check, such as naturality or
well-definedness, as it all follows by the general
theory of Kan extensions. We will use these formulas
in the proof of Theorem~\ref{Thm-NewtonAdj}.
 
Another result that will be useful is the following
fact which, although trivial, is interesting in its own
right and worth pointing out.

\begin{xlemma}
\label{Lem-BimInv}

For a pair \( (p \colon \langle A_i, {\dots}, A_n \rangle
\longtod B, \phi \colon \sum {\bf A} (A_i, -) \to \Phi) \),
the Boolean image of \( \phi \)
\begin{center}
\( \sum {\bf A} (A_i, -) \to {\rm Bim} (\phi) \ \to/^(->/ \Phi \)
\end{center}
is an invariant of the equivalence class \( [p, \phi] \).

\end{xlemma}

\begin{proof}

Suppose \( (p, \phi) \) and \( (q, \psi) \) are related by
a single morphism \( (\sigma, \langle f_j\rangle) \) of
\( {\downarrow}{\bf A} \), i.e.
\[
\bfig
\Dtriangle/>`@{>}|{\bb}`@{<-}|{\bb}/<600,350>[\langle A_1 , {\dots} , A_n\rangle
`B
`\langle C_1 , {\dots} , C_m \rangle;
(\sigma, \langle f_j \rangle)`p`q]

\Dtriangle(1400,0)/<-`>`<-/<600,350>[\sum {\bf A} (A_i, -)`\Phi
`\sum {\bf A} (C_j , -);
\sum_\sigma {\bf A} (f_j, -)`\phi`\psi]

\efig
\]
commute. Because \( (\sigma, \langle f_J\rangle) \) is in
\( {\downarrow}{\bf A} \), \( \sum_\sigma {\bf A} (f_j, -) \)
is \( \pi_0 \)-surjective, so \( {\rm Bim} (\phi) = {\rm Bim} (\psi) \).
\end{proof}

\begin{xdefinition}
\label{Def-SoftAnalytic}

A functor of the form \( \widetilde{P} \colon {\bf Set}^{\bf A}
\to {\bf Set}^{\bf B} \) for \( P \colon {\downarrow} {\bf A}
\longtod {\bf B} \) will be called {\em soft analytic}.

\end{xdefinition}

 It will become clear below that \( P \)
is uniquely determined by \( \widetilde{P} \) (see \ref{Thm-UnitIso}).

\begin{xproposition}

Analytic functors are soft analytic.

\end{xproposition}

\begin{proof}

The category \( ! {\bf A} \) of Section~\ref{SSec-MultVarAndFun} is a
subcategory of \( {\downarrow}{\bf A} \), and the \( Q \)
of \ref{SSec-MultVarAndFun}, the restriction of the one just introduced.
For an \( {\bf A} \)-\( {\bf B} \) symmetric sequence
\( P \colon !{\bf A} \longtod {\bf B} \), \( \widetilde{P} \)
is the left Kan extension
\[
\bfig
\qtriangle/ >->`>`>/<600,500>[(!{\bf A})^{op}`({\downarrow}{\bf A})^{op}
`{\bf Set}^{\bf B};`P`P']

\morphism(350,320)/=>/<150,0>[`;]

\ptriangle(600,0)/>``>/<600,500>[({\downarrow}{\bf A})^{op}`{\bf Set}^{\bf A}
`{\bf Set}^{\bf B};Q``\widetilde{P}]

\morphism(750,320)/=>/<150,0>[`;]

\efig
\]
which can be taken in stages giving, first a soft
\( {\bf A} \)-\( {\bf B} \) symmetric sequence \( P' \)
and then the analytic functor \( \widetilde{P} \) which
is isomorphic to \( \widetilde{P}' \).
\end{proof}

We can describe \( P' \) explicitly. It's the left Kan extension
of \( P \) along the inclusion \( (!{\bf A})^{op} \to/ >->/<220> 
({\downarrow}{\bf A})^{op} \) so
\[
P' (A_1 \dots A_n; B) \cong \int^{\langle C_1 \dots C_m\rangle
\in ! {\bf A}}
P (C_1 \dots C_m ; B) \times {\downarrow}{\bf A} (A_1 \dots A_n;
C_1 \dots C_m) \rlap{\ .}
\]
An element of \( P' (A_1, \dots , A_n ; B) \) is thus an
equivalence class
\[
[\langle A_1 \dots A_n \rangle \to^{(\sigma, \langle f_1 \dots f_n\rangle)}
\langle C_1 \dots C_m \rangle \todd{p}{} B]
\]
where \( \sigma \colon m \to/->>/ n \) is onto, \( f_j \colon A_{\sigma j} \to C_j \)
and \( p \in P (C_1 \dots C_m ; B) \). The equivalence relation is
generated by identifying \( (\sigma, \langle f_j \rangle, p) \) with
\( (\rho, \langle g_j \rangle, q) \) is there exists a morphism
\( (\tau, \langle h_j \rangle) \) in \( !{\bf A} \) such that
\[
\bfig
\square/@{>}|{\bb}`=`<-`>/<1050,500>[\langle A_1, \dots, A_n \rangle
`\langle C_1, \dots , C_m \rangle
`\langle A_1, \dots , A_n \rangle
`\langle D_1, \dots , D_m \rangle;
(\sigma, \langle f_i \dots f_n \rangle)`
`(\tau, \langle h_1 \dots h_m \rangle)
`(\rho, \langle g_1 \dots g_m \rangle)]

\square(1050,0)/@{>}|{\bb}``=`@{>}|{\bb}/<1050,500>[\langle C_1, \dots , C_m \rangle
`B`\langle D_1, \dots , D_m \rangle`B;
p```q]

\efig
\]
i.e.
\[
\bfig
\Ctriangle/->>`>`<<-/<500,250>[m`n`m;\sigma`\tau`\rho]

\square(900,0)/>`=`<-`>/[A_{\sigma j}`C_j`A_{\rho \tau j}`D_{\tau j};
f_j``h_j`g_{\tau j}]

\Dtriangle(2000,0)|lrr|/>`@{>}|{\bb}`@{<-}|{\bb}/<600,250>[\langle C_1 \dots C_m \rangle
`B`\langle D_1 \dots D_m \rangle;
(\tau, \langle h_1, \dots h_m \rangle)`p`q]

\place(2600,0)[.]

\efig
\]
In every equivalence class there are representatives
of the form
\[
\langle A_1 , \dots , A_n \rangle
\to^{(\sigma, \langle 1_{A_{\sigma j}} \rangle)}
\langle A_{\sigma 1}, A_{\sigma 2}, \dots , A_{\sigma m} \rangle
\todd{p}{} B
\]
and, after some calculation, we see that two such are
equivalent if and only if there is a \( \tau \in S_m \) such
that
\[
\bfig
\Ctriangle/->>`>`<<-/<500,250>[m`n`m;\sigma`\tau`\rho]

\Dtriangle(1300,0)|lrr|/<-`@{>}|{\bb}`@{<-}|{\bb}/<600,250>[\langle A_{\sigma 1}, \dots , A_{\sigma m}\rangle
`B`\langle A_{\rho 1}, \dots, A_{\rho m} \rangle;
(\tau, \langle \id_{A_{\sigma m}} \rangle)`p`q]

\efig
\]

We can further nail down the equivalence class by
choosing canonical surjections \( m \to/->>/ n \), the
order preserving ones, and these are determined by
their fibres \( m_i \) which are positive integers. This
gives a relatively simple description of \( P' \)
\[
P' (A_1, \dots A_n ; B) \cong \sum_{m_1, \dots m_n > 0}
P (A_1^{\otimes m_1}, \dots , A_n^{\otimes m_n} ; B) /
S_{m_1} \times {\dots} \times S_{m_n}
\]
where \( A_i^{\otimes m_i} = \langle A_i, A_i, \dots , A_i \rangle
\in {\bf A}^{m_i} \) and the action is by permuting those
entries.

\subsection{The Newton series comonad}
\label{SSec-NewtonCom}

In this section we show that taking iterated differences
is right adjoint to summation of a multivariable symmetric
series. We first combine all the iterated differences
into one soft symmetric sequence.

\begin{xproposition}
\label{Prop-DeltaSymSeq}

Let \( F \colon {\bf Set}^{\bf A} \to {\bf Set}^{\bf B} \)
be tense. Then taking the iterated symmetric differences
 of \( F \) evaluated at \( \Phi \) gives an \( {\bf A} \)-\( {\bf B} \)
 symmetric sequence
 \[
 \Delta_* [F] (\Phi) \colon {\downarrow}{\bf A} \longtod B
 \]
 \[
 \Delta_* [F] (\Phi) (A_1, \dots , A_n ;B) =
 \Delta_{A_1} \dots \Delta_{A_n} [F] (\Phi) (B) \rlap{\ .}
 \]

\end{xproposition}

\begin{proof}

\( \Delta_{A_1} \dots \Delta_{A_n} [F] (\Phi) (B) =
\Delta_{\langle A_i \rangle} [F] (\Phi)(B) \) consists
of the new elements of 
\[ F(\Phi + {\bf A} (A_1, -) +
{\dots} + {\bf A} (A_n, -)) (B) \rlap{\ ,}
\]
i.e.~those elements not
in \( F( \Phi + {\bf A} (A_{\alpha 1}, -) + {\dots} + {\bf A}
(A_{\alpha k}, -)) \) for any proper subsequence
\( \langle A_{\alpha 1}, \dots , A_{\alpha k} \rangle \),
\( \alpha \colon k\  \to/ >>->/ n \) a proper mono. We'll
show that \( \Delta_* [F] \) is a subfunctor of
\( F(\Phi + Q) \). Let \( (\sigma, \langle f_1, \dots , f_m\rangle))
\colon \langle A_1, \dots , A_n \rangle \to \langle C_1, \dots ,
C_m\rangle \) be a morphism in \( {\downarrow}{\bf A} \), and
let \( x \) be an element of 
\[ 
\Delta_{C_1} \dots \Delta_{C_m}
[F] (\Phi) (B) \subseteq F (\Phi + {\bf A} (C_1, -) + {\dots} +
{\bf A} (C_m, -)) (B)\rlap{\ .}
 \]
Then \( y = F (\sigma, \langle f_1, \dots ,
f_m \rangle) (B) (x)  \)
is an element of \( F (\Phi + {\bf A}
(A_1, -) + {\dots} + {\bf A} (A_n, -) (B) \) and suppose it's
not new. There is a proper monomorphism \( \alpha \colon k\ \to/ >>->/ n \)
such that \( y \in F (\Phi _ {\bf A} (A_{\alpha 1}, -) + {\dots} +
{\bf A} (A_{\alpha k}, -)) (B) \).

The pullback of a proper mono along an epi is again proper
so we get
\[
\bfig
\square/ >>->`->>`->>` >>->/[l\ `m`k\ `n;\beta`\rho`\sigma`\alpha]

\place(250,250)[\Pb]

\efig
\]
which, in turn, gives a pullback of complemented subobjects
in \( {\bf Set}^{\bf A} \)
\[
\bfig
\square/^(->`>`>`^(->/<1500,500>[{\bf A} (C_{\beta 1}, -) + {\dots} +
{\bf A} (C_{\beta l}, -)\ 
`\ {\bf A} (C_1, -) + {\dots} + {\bf A} (C_m, -)
`{\bf A} (A_{\alpha 1}, -) + {\dots} + {\bf A} (A_{\alpha k}, -)\ 
`{\bf A} (A_1, -) + {\dots} + {\bf A} (A_n, -)\rlap{\ .};
`(\rho, \langle f_{\beta 1}, \dots , f_{\beta n} \rangle)
`(\sigma, \langle f_1, \dots , f_m \rangle)`]

\place(750,250)[\Pb]

\efig
\]
Adding \( \Phi \) produces another such pullback and
\( F \), being tense, will preserve it
\[
\bfig
\square/^(->`>`>`^(->/<1850,500>[F (\Phi + {\bf A} (C_{\beta 1}, -) + {\dots} +
{\bf A} (C_{\beta l}, -))\ 
`\ F (\Phi + {\bf A} (C_1, -) + {\dots} + {\bf A} (C_m, -))
`F (\Phi + {\bf A} (A_{\alpha 1}, -) + {\dots} + {\bf A} (A_{\alpha k}, -))\ 
`F (\Phi + {\bf A} (A_1, -) + {\dots} + {\bf A} (A_n, -))\rlap{\ .};
```]

\place(925,250)[\Pb]

\efig
\]
Then \( x \) in the upper right corner gets sent to \( y \) which
is in the lower left corner, so \( x \) itself is in the upper
left corner, i.e.~\( x \) wasn't new after all. Thus
\( \Delta_* [F] (\Phi) \) is a subfunctor of \( F (\Phi + Q) \).
\end{proof}

\( \Delta_* [F] (\Phi) \) is functorial in \( F \). Indeed, applying
Proposition~\ref{Prop-PartTense} recursively, we see that
any tense transformation \( t \colon F \to G \) restricts to
\[
\bfig
\square/^(->`-->`>`^(->/<1500,500>[\Delta_{A_1} \dots \Delta_{A_n}{[}F{]}(\Phi)\ 
`F(\Phi + {\bf A} (A_1, -) + {\dots} + {\bf A}(A_n, -))
`\Delta_{A_1} \dots \Delta_{A_n}{[}G{]}(\Phi)\ 
`G(\Phi + {\bf A}(A_1, -) + {\dots} + {\bf A}(A_n, -));
``t(\Phi + {\bf A}(A_1, -) + {\dots} + {\bf A}(A_n, -))`]

\efig
\]
which will be natural and functorial automatically. Thus for
each \( \Phi \) in \( {\bf Set}^{\bf A} \) we get a functor
\[
\Delta [\ \ ] (\Phi) \colon \Tense ({\bf Set}^{\bf A}, {\bf Set}^{\bf B})
\to \Prof ({\downarrow} {\bf A}, {\bf B}) \rlap{\ ,}
\]
i.e.~\( \Delta_* [F] (\Phi) \) is an \( {\bf A} \)-\( {\bf B} \) soft
symmetric sequence.

The main result of this section is the following:

\begin{xtheorem}
\label{Thm-NewtonAdj}

\[
\Delta_* [\ \ ] (0) \mbox{\ is right adjoint to\ } \tilde{(\ )} \rlap{\ .}
\]

\end{xtheorem}

\begin{proof}

\( \widetilde{P} \) is the left Kan extension of \( P \) along
\( Q \)
\[
\bfig
\Vtriangle<400,500>[({\downarrow}{\bf A})^{op}
`{\bf Set}^{\bf A}`{\bf Set}^{\bf B};
Q`P`{\rm Lan}_Q P = \widetilde{P}]

\morphism(350,300)/=>/<150,0>[`;]

\efig
\]
so for any functor \( F \colon {\bf Set}^{\bf A} \to
{\bf Set}^{\bf B} \) we have a bijection
\[
\widetilde{P} \to^t F \over P \to_u FQ
\]
as discussed above. Now \( \Delta_* [F] (0) \) is a
subfunctor of \( F Q \). Indeed
\[
\Delta_* [F] (0) \langle A_1, \dots , A_n \rangle (B) =
\Delta_{A_1} \dots \Delta_{A_n} [F] (0) (B)
\]
consists of the new elements of
\[
F(Q \langle A_1, \dots , A_n \rangle) (B) =
F ({\bf A} (A_1, -) + {\dots} + {\bf A}(A_n, -))(B) \rlap{\ .}
\]

We'll show that \( t \colon \widetilde{P} \to F \) is tense
if and only if \( u \) factors through \( \Delta_* [F] (0) \ 
\to/^(->/ F Q \) which will establish the theorem.

First assume \( t \) is tense. Let \( p \) be in
\( P (A_1, \dots , A_n ; B) \) so \( u (p) \) is in
\( F ({\bf A} (A_1, -) + {\dots} + {\bf A} (A_n, -)) (B) \)
and assume \( u (p) \) is in \( F \) of some subsum
\( F ({\bf A} (A_{\alpha 1}, -) + {\dots} + {\bf A}
(A_{\alpha k} , -)) (B) \) for a subset \( \alpha \colon
k \to/ >->/ n \) of the indices. Tenseness of \( t \)
applied to the complemented subsum
\( \mu \colon \sum {\bf A}(A_{\alpha i}, -)\ \to/^(->/
\sum{\bf A} (A_i, -) \) gives a pullback
\[
\bfig
\square/>`<-_)`<-_)`>/<2000,500>[\int^{\langle C_j\rangle \in {\downarrow}{\bf A}}
P (\langle C_j\rangle;B) \times {\bf Set}^{\bf A}
(\sum {\bf A} (C_j, -), \sum {\bf A} (A_i, -))(B)
`F (\sum {\bf A} (A_i, -))(B)
`\bgr{\int^{\langle C_j\rangle \in {\downarrow}{\bf A}}
P (\langle C_j\rangle;B) \times {\bf Set}^{\bf A}
(\sum {\bf A} (C_j, -), \sum {\bf A} (A_{\alpha i}, -))(B)}
`\bgr{F (\sum {\bf A} (A_{\alpha i}, -))(B)\rlap{\ .}};
```]

\place(1050,250)[\Pb]

\efig
\]
Then \( u (p) = T [p, \id_{\sum {\bf A}(A_i, -)}] \) is in
\( F(\sum{\bf A} (A_{\alpha i}, -)) (B) \) so \( [p, \id] \) is in
the lower left corner which means there are
\( q \colon \langle C_1, \dots , C_m \rangle \longtod B \)
and \( \psi \colon \sum {\bf A} (C_i, -) \to \sum {\bf A}
(A_{\alpha i}, -) \) such that \( [q, \mu \psi] =
[p, \id] \). Thus by Lemma~\ref{Lem-BimInv} we see
that \( {\rm Bim} (\mu \psi) = {\rm Bim} (\id) =
\sum {\bf A}(A_i, -) \). It follows that \( \mu \) is the identity,
so \( u (p) \) is not contained in \( F \) of any proper
subsum, i.e.~is new. This gives our factorization of \( u \)
through \( \Delta_* [F] (0) \).

Conversely, assume that \( u \) factors through
\( \Delta_* [F] (0) \). We'll show that \( t \) is tense.
Let \( \Psi\  \to/^(->/ \Phi \) be a complemented subobject.
We must show that
\begin{equation}\tag{*}
\bfig
\square/>`<-_)`<-_)`>/<550,500>[\widetilde{P} (\Phi)`F(\Phi)
`\bgr{\widetilde{P}(\Psi)}`\bgr{F(\Psi)};
t(\Phi)```t(\Psi)]

\efig
\end{equation}
is a pullback. Take an element \( [p \colon \langle A_1, \dots ,
A_n \rangle \longtod B, \phi \colon \sum {\bf A} (A_i, -) \Phi] \)
of \( \widetilde{P} (\Phi) (B) \) and assume \( t (\Phi) [p, \phi] =
F (\phi) (p) \) is in \( F (\Psi) \). Form the pullback
\[
\bfig
\square/>`<-_)`<-_)`>/<700,500>[\sum {\bf A} (A_i, -)`\Phi
`\bgr{\sum {\bf A} (A_{\alpha j}, -)}`\bgr{\Psi}\rlap{\ .};
\phi```\psi]

\place(350,250)[\Pb]

\efig
\]
It is induced by a monomorphism \( \alpha \colon 
m \to/ >->/<250> n \) because a complemented subobject
of a sum of representables is a subsum. We get a new
pullback now by tenseness of \( F \)
\[
\bfig
\square/>`<-_)`<-_)`>/<750,500>[F(\sum {\bf A} (A_i, -))`F\Phi
`\bgr{F(\sum {\bf A} (A_{\alpha j}, -))}`\bgr{F\Psi\rlap{\ .}};
F(\phi)```F(\psi)]

\place(350,250)[\Pb]

\efig
\]
\( F (\phi) \) takes \( u (p) \) to an element of \( F(\Psi) \) so
\( u (p) \in F (\sum {\bf A}(A_{\alpha j}, -)) \). But \( u (p) \)
was supposed to be a new element of \( F(\sum {\bf A}
(A_i, -)) \) so \( \alpha \) is not a proper subsum which means
that
\[
\bfig
\qtriangle/>`>`<-_)/<600,400>[\sum {\bf A} (A_i, -)`\Phi`\bgr{\Psi\rlap{\ .}};
\phi`\psi`]

\efig
\]
Thus \( [p, \phi] \) is in \( \widetilde{P} (\Psi) \). This shows that
our square \( (*) \) is indeed a pullback.
\end{proof}

The adjoint pair \( \tilde{(\ )} \dashv \Delta_* [\ \ ] (0) \) induces
a comonad on \( \Tense ({\bf Set}^{\bf A}, {\bf Set}^{\bf B}) \)
which we call the {\em Newton series comonad}.

\subsection{Convergence}
\label{SSec-Convergence}

In this section we show that the Newton series for a soft
analytic functor ``converges to it''.

\begin{xtheorem}
\label{Thm-UnitIso}

For every \( {\bf A} \)-\( {\bf B} \) soft symmetric sequence
\( P \colon {\downarrow} {\bf A} \longtod {\bf B} \),
the unit for the adjunction of \ref{Thm-NewtonAdj}
\[
P \to \Delta_* [\widetilde{P}] (0)
\]
is an isomorphism.

\end{xtheorem}

\begin{proof}

An element of \( \Delta_* [\widetilde{P}] (0) \) at
\( \langle A_1, \dots , A_n \rangle, B \) is a new element
of \( \widetilde{P} (\sum {\bf A}(A_i, -)) (B) \), i.e.~of
\begin{center}\(
\int^{C_1, \dots, C_m \in {\downarrow} {\bf A}}
P (C_1, \dots, C_m; B) \times {\bf Set}^{\bf A}
(\sum {\bf A}(C_j, -), \sum {\bf A} (A_i, -))
\)
\end{center}
which is an equivalence class
\begin{center}\(
\big[p \colon \langle C_1, \dots, C_m \rangle \longtod B,
\phi \colon \sum {\bf A} (C_j, -) \to
\sum {\bf A} (A_i, -)\big]
\)
\end{center}
(satisfying the newness condition, of course).

The unit \( P \to \Delta_* [\widetilde{P}] (0) \) takes
\( p \colon \langle A_1, \dots, A_n \rangle \longtod B \)
to the equivalence class
\begin{center}
\(
\big[p \colon \langle A_1, \dots, A_n \rangle \longtod B,
\id \colon \sum {\bf A}  (A_i, -) \to \sum {\bf A} (A_i, -)\big]
\rlap{\ .}
\)
\end{center}
A \( \phi \) as above is, as explained in the discussion
around Proposition~\ref{Prop-TransReps}, of the form
\( \sum_\alpha {\bf A} (f_j, -) \) for \( \alpha \colon m \to n \)
and \( f_j \colon A_{\alpha j} \to C_j \) and we can take
its Boolean factorization by factoring \( \alpha \)
(in \( {\bf Set} \))
\[
\bfig
\Vtriangle/>`->>`<-</<400,350>[m`n`k\ ;\alpha`\sigma`\mu]

\efig
\]
and taking
\[
\sum_{j \in m} {\bf A} (C_j, -)
\to^{\sum_\sigma {\bf A} (f_j, -)}
\sum_{i \in k} {\bf A} (A_{\mu i}, -)\ 
\to/^(->/^{\sum_\mu {\bf A} (1_{\mu i}, -)}
\sum_{i \in n} {\bf A} (A_i, -) \rlap{\ .}
\]
If \( \mu \) were a proper mono, \( [p, \phi] \) wouldn't be
new as it would be in \( \widetilde{P} (\sum_{i \in k}
{\bf A} (A_{\mu i}, -)) \), so \( \mu = \id_n \) and \( \alpha =
\sigma \), a surjection. Thus \( (\sigma, \langle f_1\rangle) \)
is a morphism of \( {\downarrow} {\bf A} \) and we have
\[
\bfig
\Dtriangle/<-`@{>}|{\bb}`@{<-}|{\bb}/<550,350>[\langle C_1, \dots, C_m\rangle
`B`\langle A_1, \dots, A_n \rangle;(\sigma, \langle f_j\rangle)`p`p']

\Dtriangle(1300,0)/>`>`<-/<550,300>[\sum {\bf A} (C_j, -)
`\quad\quad\quad\sum {\bf A} (A_i, -)`\sum {\bf A} (A_i, -);
\sum_\alpha {\bf A} (f_j, -)`\phi`\id]

\efig
\]
so \( [p, \phi] = [p', \id] \), which shows that the unit
\[
P \to \Delta_* [\widetilde{P}] (0)
\]
\[
p \longmapsto [p, \id]
\]
is onto.

To show that the unit is one-one we must show that if
\( [p, \id] = [q, \id] \) then \( p = q \). \( [p, \id] = [q, \id] \) means
there's a zigzag path of
\[
\bfig
\Dtriangle/<-`@{>}|{\bb}`@{<-}|{\bb}/<550,350>[\langle C_1, \dots, C_m\rangle
`B`\langle D_1, \dots, D_r \rangle;(\rho, \langle h_j\rangle)`\ov{p}`\ov{q}]

\Dtriangle(1300,0)/>`>`<-/<550,300>[\sum {\bf A} (C_j, -)
`\quad\quad\quad\sum {\bf A} (A_i, -)`\sum {\bf A} (D_s, -);
\sum_\rho {\bf A} (h_j, -)`\phi`\psi]

\efig
\]
with \( (\rho , \langle h_j\rangle) \) in \( {\downarrow}{\bf A} \)
joining \( [p, \id] \) to \( [q, \id] \). The Boolean image of
\( \phi_i \) (and \( \psi_i \)) is an invariant of the equivalence
class (\ref{Lem-BimInv}) and as \( {\rm Bim} (\id) =
\sum {\bf A} (A_i, -) \), all the \( \phi \) and \( \psi \) also
have \( \sum {\bf A} (A_i, -) \) as their images. That means
that the morphisms \( (\sigma, \langle f_j \rangle) \) and
\( (\tau, \langle g_s \rangle) \) corresponding to \( \phi \)
and \( \psi \) are actually morphisms in \( {\downarrow} {\bf A} \),
i.e.~\( \sigma \colon m \to n \) and \( \tau \colon r \to  n \)
are surjections. Now we have
\[
\bfig
\Ctriangle/<-`<-`>/<600,450>[\langle C_1, \dots, C_m\rangle
`\langle A_i, \dots, A_n\rangle`\langle D_1, \dots, D_r\rangle;
(\sigma, \langle f_j\rangle)`(\rho, \langle h_j\rangle)
`(\tau, \langle g_s\rangle)]

\Dtriangle(600,0)/<-`@{>}|{\bb}`@{<-}|{\bb}/<600,450>[\langle C_1, \dots, C_m\rangle
`B`\langle D_1, \dots, D_r\rangle;`\ov{p}`\ov{q}]

\efig
\]
commuting \( \ov{p} (\sigma, \langle f_j\rangle) =
\ov{q} ( \tau, \langle g_s\rangle) \) at every stage of the path
joining \( (p, \id) \) to \( (q, \id) \), and for these endpoints
we get \( p \) and \( q \) respectively, i.e.~\( p = q \).
\end{proof}

This shows that the Newton series comonad is idempotent.

\begin{xcorollary}

If \( F \colon {\bf Set}^{\bf A} \to {\bf Set}^{\bf B} \) is soft
analytic (in particular analytic) then its Newton series
converges to it, i.e.~the counit
\[
\widetilde{\Delta_*[F](0)} \to F
\]
is an isomorphism.

\end{xcorollary}


\subsection{Concluding remark}

In the previous sections, we touted the functor taking
\( F \) to \( \ov{F} = \widetilde{\Delta_* [F] (0)} \) as a
categorical version of the Newton summation formulas
at the beginning of \ref{SSec-Newton}, but in fact it
looks nothing like them.

Let's consider the first one
\[
\ov{f} (x_1, \dots, x_n) = \sum_{k_1, \dots, k_n}
\frac{\Delta^{k_1}_{x_1} {\dots} \Delta^{k_n}_{x_n} [f] (0, \dots, 0)}
{k_{1!} \dots k_{n!}}
\ x^{\downarrow k_1}_1 \dots x^{\downarrow k_n}_n
\]
where \( f \) is a function \( {\mathbb R}^n \to {\mathbb R} \)
and the sum is taken over all \( n \)-tuples of natural numbers.
We've replaced \( f \) by a (tense) functor \( {\bf Set}^{\bf A}
\to {\bf Set}^{\bf B} \) and the difference operators by our
functorial ones, but it's not clear how to interpret the rest
of the formula. Let's look at it more carefully.

The first thing to note is that, while the \( x_i \) in
\( \Delta_{x_i} \) and in \( x_i^{\downarrow k_i} \) refer to the
same thing, they play different roles. The \( x_i \) in
\( \Delta_{x_i} \) is merely a subscript indicating which
difference operator is used, and we could well have written
\( \Delta_i \) instead, although \( \Delta_{x_i} \) is more
descriptive. The \( x_i \) in \( x_i^{\downarrow k_i} \), on
the other hand, represents a variable which can take values,
\( c_i \). So we have
\[
\ov{f} (c_1, \dots, c_n) = \sum_{k_1, \dots, k_n}
\frac{\Delta^{k_1}_{x_1} {\dots} \Delta^{k_n}_{x_n} [f] (0, \dots, 0)}
{k_{1!} \dots k_{n!}}
\ c^{\downarrow k_1}_1 \dots c^{\downarrow k_n}_n \rlap{\ .}
\]

Here all the like \( \Delta \)'s have been grouped together
which is fine as we have finitely many variables and they're
totally ordered. It would be more natural to sum over all
finite sequences of variables \( \langle x_{\alpha (1)} \dots
x_{\alpha (m)} \rangle \) and group the terms together by the
length \( m \). Of course we get more terms:
\( \Delta^{k_1}_{x_1} \dots \Delta^{k_n}_{x_n} \) gets counted
\[
\binom{k_1 + {\dots} + k_n}{k_1, \dots , k_n} \ =\ 
\frac{(k_1 + {\dots} + k_n)!}{k_1 ! \dots k_n !} \ =\ 
\frac{m !}{k_1 ! \dots k_n !}
\]
times, so now we have
\[
\ov{f} (c_1, \dots, c_n) =
\sum_{\alpha \colon m \to<100> n}
\frac{\Delta_{\alpha (1)} \dots \Delta_{\alpha (m)} [f] (0, 0)}
{m !}
\ c^{\downarrow k_1}_1 \dots c^{\downarrow k_n}_n \rlap{\ .}
\]

In fact this takes care of the finiteness and total ordering
of the variables, as far as the \( \Delta \) part of the formula
is concerned. We take a set of variables \( {\rm Var} \) and
consider the free monoid on it \( {\rm Var}^* \), over which
the sum is to be taken. The \( c_i \) are a choice of value
for each variable \( \phi \colon {\rm Var} \to {\mathbb R} \)
but we still have to deal with the \( k_i \) in this setup.

The \( k \)'s count the number of occurrences of a given
variable \( y \) in a sequence \( \langle x_1, \dots, x_n\rangle \).
Let\[
\delta \colon {\rm Var} \times {\rm Var} \to {\mathbb N}
\]
be the Kronecker delta, i.e.~\( \delta (x, y) = 1 \) if \( x = y \)
and \( 0 \) otherwise. For each \( y \), extend \( \delta (-, y) \)
to a function \( \delta (-, y) \colon {\rm Var}^* \to {\mathbb N} \)
using the additive structure of \( {\mathbb N} \), so
\[
\delta (x_1, \dots, x_n; y) =
\sum^m_{i = 1} \delta (x_i, y)
\]
is exactly the number of \( y \)'s in \( \langle x_1, \dots, x_n\rangle \).
Thus we end up with the Newton series in the form we want
\[
\ov{f} (\phi) = \sum_{\langle x_1, \dots, x_n\rangle \in {\rm Var}^*}
\frac{\Delta_{x_1} \dots \Delta_{x_m} [f](0)}{m !}
\prod_{y \in {\rm Var}} \phi (y)^{\downarrow \delta (x_1, \dots x_n ;y)}
\]
which, admittedly, looks more complicated than the original but
it's the closest we can get to the categorical version.

Now the Newton series comonad of Section~\ref{SSec-NewtonCom}
\[
\widetilde{\Delta_*[F] (0)} =
\int^{\langle A_1, \dots, A_m\rangle \in \downarrow {\bf A}}
\Delta_{A_1} \dots \Delta_{A_m} [F] (0) \times {\bf Set}^{\bf A}
({\bf A} (A_1, -) + {\dots} + {\bf A} (A_m, -), \Phi) 
\]
looks similar to the above, with the following
correspondences:
\[
\begin{array}{rcl}
f \colon {\mathbb R}^n \to {\mathbb R}  & \leftrightarrow  & F \colon {\bf Set}^{\bf A} \to {\bf Set}^{\bf B}\\
\mbox{variables  } x                               & \leftrightarrow  & \mbox{objects } A \mbox{ of } {\bf A}\\
{\rm Var}                                                 & \leftrightarrow  & {\bf A}\\
{\rm Var}^*                                              & \leftrightarrow  & {\downarrow}{\bf A}\\
\phi \colon {\rm Var} \to {\mathbb R}       & \leftrightarrow  & \Phi \colon {\bf A} \to {\bf Set}\\
\delta (x, y)                                              & \leftrightarrow  & {\bf A} (A', A)\\
\delta (x_1, \dots, x_n, y)                        & \leftrightarrow  & {\bf A} (A_1, A) + {\dots} + {\bf A} (A_m, A)\\
\prod {\phi (y)}^{\downarrow \delta (x_1 \dots x_n;y)} & \leftrightarrow  & {\bf Set}^{\bf A} ({\bf A} (A_1, -)+ {\dots} + {\bf A} (A_m, -), \Phi)
\end{array}
\]
The correspondence is not perfect, of course. \( {\rm Var}^* \) might
rightly be said to correspond to \( ! {\bf A} \) rather than
 \( {\downarrow} {\bf A} \). Then the \( m ! \) in the sum is incroporated
 in the coend via the symmetric groups.
 
 Also \( \prod \phi (y)^{\downarrow \delta (x_1, \dots, x_m, y)} \)
 should correspond to monomorphisms 
 \[
  {\bf A} (A_1, -) + {\dots} +
 {\bf A} (A_m, -) \to \Phi 
 \]
 rather than arbitrary natural
 transformations. That's what the extra morphisms in
 \( {\downarrow}{\bf A} \) (involving surjections \( \sigma \))
 take care of. We need a bit more theory to explain this.

 \begin{xdefinition}
 \label{Def-Diverse}
 
 Let \( \Phi \colon {\bf A} \to {\bf Set} \) and \( x \in \Phi A \).
 An {\em ancestor} of \( x \) is a \( y \in \Phi A' \) for which
 there is a morphism \( f \colon A' \to A \) such that
 \( \Phi (f)(y) = x \). Two elements \( x_1 \in \Phi A_1 \) and
 \( x_2 \in \Phi A_2 \) are {\em relatives} if they have
 a common ancestor. A sequence \( \langle x_1 \in \Phi A_1,
 \dots, x_n \in \Phi A_n \rangle \) is called {\em diverse}
 if no two elements are relatives. A natural transformation
 \( \phi \colon \sum {\bf A} (A_i, -) \to \Phi \) is {\em diverse}
 if the corresponding sequence of elements
 \( \langle \phi (A_i) (1_{A_i}) \rangle \) is.
 
 \end{xdefinition}
 
 All the elements of a diverse sequence are different
 and more, but not enough more to make the corresponding
 transformation monic. One could have \( i \neq j \)
 and \( f \colon A_i \ \to A, g \colon A_j \to A \) with
 \( \Phi (f) (x_i) = \Phi (g) (x_j) \). But if \( {\bf A} \) is
 a groupoid, then \( \phi \) is monic if and only if it is
 diverse. The variables \( x_1, \dots, x_n \) in the
 formula we're abstracting from form a finite discrete
 set so diverse restricts to one-one in that case.

 \begin{xproposition}
 \label{Prop-DiverseFact}
 
 \begin{itemize}
 
 	\item[(1)] \( \phi \) as below is diverse if and only if
	for every factorization of \( \Phi \)
	\[
	\bfig
	\Dtriangle/>`>`<-/<550,400>[\sum {\bf A}(A_i, -)`\Phi`\sum {\bf A}(C_j, -);
	\sum_\sigma {\bf A}(f_i, -)`\phi`\psi]
	
	\efig
	\]
	with \( (\sigma, \langle f_i\rangle) \colon \langle C_1, \dots , C_m \rangle
	\to \langle A_1, \dots, A_n \rangle \) in \( {\downarrow} {\bf A} \),
	we have that \( \sigma \) is a bijection, i.e.~\( (\sigma, \langle f_i\rangle)
	\in ! {\bf A} \).
	
	\item[(2)] Every \( \phi \) factors as \( \psi \sum_\sigma {\bf A}
	(f_i, -) \) with \( (\sigma, \langle f_i\rangle) \in {\downarrow}{\bf A} \)
	and \( \psi \) diverse.
 
 \end{itemize}
 
 \end{xproposition}

 \begin{proof}
 
 (1) \( \phi \) and \( \psi \) as in the statement correspond
 to an \( n \)-tuple \( x_1 \in A_1, \dots, x_n \in \Phi A_n \)
 and an \( m \)-tuple \( y_1 \in \Phi C_1, \dots, y_m \in \Phi C_m \),
 respectively. The \( x \)'s and \( y \)'s are related by
 \[
 x_i = \Phi (f_i) (y_{\sigma i}) \rlap{\ .}
 \]
 
 If \( \sigma \) is not one-to-one, say \( \sigma (i_1) =
 \sigma (i_2) \), then \( x_{i_1} \) and \( x_{i_2} \) are
 relatives as they have the common ancestor
 \( y_{\sigma (i_1)} = y_{\sigma (i_2)} \). So the \( x_i \)
 are not diverse nor is \( \phi \).
 
 Conversely, if the \( x_i \) are not diverse, then there
 are two \( x \)'s that are relatives. Assume, for
 simplicity of notation, that they are \( x_{n - 1} \)
 and \( x_n \). So we have \( f \colon C \to A_{n - 1},
 g \colon C \to A_n \) and \( y \in \Phi C \) such that
 \( \Phi (f) (y) = x_{n - 1} \) and \( \Phi (g) (y) = x_n \).
 Then we get a morphism
 \[
 (\sigma, \langle f_i\rangle) \colon \langle A_1, \dots, A_{n - 2}, C
 \rangle \to \langle A_1, \dots, A_n \rangle
 \]
 \[
 \sigma (i) = \left\{ \begin{array}{lcl}
 i        & \mbox{if}    &  i < n\\
 n -1  &  \mbox{if}    &  i = n \rlap{\ ,}
 \end{array}
 \right.
 \]
 
 \[
  \langle f_i \rangle = \langle 1_{A_1}, \dots, 1_{A_{n - 2}}, f, g\rangle \rlap{\ .}
\]
Let \( \langle y_1, \dots, y_{n - 1} \rangle = \langle x_1, \dots,
x_{n - 2}, y \rangle \). Then
\[
x_i = \Phi (f_i) (y_{\sigma i})
\]
so the \( y \) determine a \( \psi \) giving a factorization as
above, and \( \sigma \) is not a bijection.

This proves (1).

(2) If \( \phi \) is not diverse, there exists a factorization
as in (1) with \( \sigma \) onto but not one-to-one, so
\( \sum {\bf A} (C_j, -) \) has fewer terms than \( \sum
{\bf A} (A_i, -) \). If we take, among all factorizations,
one with the minimal number of terms, the \( \psi \) must
be diverse, otherwise we could factor it again and get
a smaller one.
 \end{proof}

 \begin{xcorollary}
 
 Every equivalence class
 \begin{center}\(
 \big[ x \in F (\sum {\bf A} (A_i, -)) (B),
 \phi \colon \sum {\bf A} (A_i, -) \to \Phi \big] \)
 \end{center}
 in
 \[
 \int^{\langle A_1, \dots, A\rangle \in {\downarrow}{\bf A}}
 \Delta_{A_1} \dots \Delta_A [F] (0) (B) \times
 {\bf Set}^{\bf A} \big(\sum {\bf A} (A_i, -), \Phi\big)
 \]
 has a representative in which \( \phi \) is diverse.
 
 \end{xcorollary}
 
 \begin{proof}
 
 Factor \( \phi \) as in \ref{Prop-DiverseFact} (2) above.
 Then
 \[
 \bfig
 \node a(0,0)[y]
 \node b(150,0)[\in]
 \node c(700,0)[F(\sum {\bf A} (C_j, -)) (B)]
 
 \node g(700,500)[]
 
 \node d(0,500)[x]
 \node e(150,500)[\in]
 \node f(700,500)[F(\sum {\bf A} (A_i, -))(B)]
 
 \arrow/|->/[d`a;]
 \arrow|r|/->/[g`c;F(\sigma, \langle f_i\rangle)]

 \square(2000,0)/>`>`=`>/<600,500>[\sum{\bf A} (A_i, -) `\Phi
 `\sum {\bf A} (C_j, -)`\Phi;\phi
 `\sum_\sigma {\bf A} (f_i, -)``\psi]
 
 \efig
 \]
 so \( [x, \phi] = [y, \psi] \) and \( \psi \) is diverse.
 \end{proof}

 The diverse transformations are our categorified
 set injections so
\begin{center}\(
 \mbox{Diverse}\  (\sum {\bf A} (A_i, -), \Phi)
 \)
 \end{center}
 is our version of falling power. Note, however, that
 it is not functorial, and we need all transformations
 to make it so.

\bibliography{Multivariate}
\bibliographystyle{msclike}

\end{document}